\documentclass[a4paper,reqno,12pt,oneside]{amsart}

\usepackage{amsmath,geometry}
\usepackage{xcolor}        
\usepackage{latexsym}
\usepackage{amsmath}
\usepackage{amsfonts}
\usepackage{amssymb}
\usepackage{graphicx}
\usepackage{esint}
\usepackage{slashed}
\usepackage{enumerate}   
\usepackage[scr=boondoxo]{mathalpha}

\usepackage{mathtools}

\DeclareMathOperator*{\osc}{osc}

\numberwithin{equation}{section}

\newtheorem{thm}{Theorem}[section]

\newtheorem{pro}[thm]{Proposition}
\newtheorem{lm}[thm]{Lemma}
\newtheorem{oss}[thm]{Remark}
\newtheorem{cor}[thm]{Corollary}

\usepackage{setspace}

\title[Blow-up profile decomposition and optimal Harnack inequalities]{Blow-up profile decomposition and optimal Harnack inequalities for the $n$-Liouville equation}

\author{Pierpaolo Esposito}
\address{Pierpaolo Esposito, Dipartimento di Matematica e Fisica, Universit\`a degli Studi Roma Tre, Largo S.~Leonardo Murialdo 1, Roma 00146, Italy.}
\email{pierpaolo.esposito@uniroma3.it}

\author{Marcello Lucia}
\address{Marcello Lucia, Department of Mathematics, College of Staten Island $\&$ The Graduate Center, The City University of New York, 365 Fifth Ave., New York, NY 10016, USA.}
\email{marcello.lucia@csi.cuny.edu}

\begin{document}

\begin{abstract}
We consider the problem 
$- \Delta_n u = h(x) e^u$ in a bounded domain $\Omega$ of $\mathbb R^n$, where $\Delta_n$ stands for the $n$-Laplacian and $h$ is a continuous function satisfying 
$0 < a \leq h \leq b$.
Several new a priori estimates are obtained for the set of solutions whose boundary oscillation is uniformly bounded. 
Those results are exploited to provide a profile decomposition in $L^\infty$-sense for a sequence of blowing-up solutions.
We also present an optimal Harnack type inequality, which states that on a compact set $K \subset \Omega$ we have a universal estimate
$\displaystyle \sup_K u + (n-1)\displaystyle \inf_{\Omega} u \leq C$. 
\end{abstract}

\maketitle

\section{Introduction}
\medskip

In the present paper we are concerned with solutions to 
\begin{equation} \label{Liouville-Quasilinear}
      -\Delta_n u = h(x) e^u  \hbox{ in } \Omega ,
\end{equation}
where $\Omega \subset \mathbb{R}^n$, $n \geq 2$, is a bounded open set, $h \in C(\Omega)$  and $\Delta_n u = \hbox{div} ( |\nabla u|^{n-2} \nabla u) $ stands for the $n$-Laplace operator. Solutions are meant in a weak sense and by elliptic estimates \cite{Dib,Ser1,Tol} such solutions are in $C^{1,\alpha}(\Omega)$ for some $\alpha \in (0,1)$. 

\medskip

When $n=2$ and $h$ is a constant, the PDE~\eqref{Liouville-Quasilinear} goes back to Liouville \cite{Liou}, and is the simplest case of the problem of prescribing the Gauss curvature on a two-dimensional surface, 
a question that has sparked a renewed interest since the '70s through the works of Moser \cite{Mos}, Kazdan and Warner \cite{KW}.
As observed in  \cite{MaQing}, 
such a question admits  a generalization in $\mathbb R^n$ with $n \geq 3$, since after a conformal change of the Euclidean metric $g = e^{2u} {\rm id}$ the Ricci curvature satisfies 
${\rm Ric} (\frac{\nabla^g u}{|\nabla^g u|_g}, \frac{\nabla^g u}{|\nabla^g u|_g})  |\nabla^g u|_g^{n-2} = - \frac{ \Delta_n u }{e^{nu}}$. 
Moreover, the $n$-Liouville equation~\eqref{Liouville-Quasilinear} represents a simplified version of a quasilinear fourth-order problem arising in the theory of log-determinant functionals  \cite{ChYa,EsMa}, 
which are relevant in the study of the conformal geometry of a $4$-dimensional closed manifold. In order to understand some of the bubbling phenomena that may occur in such geometrical settings, we are naturally led to study the simplest situation given by~\eqref{Liouville-Quasilinear}.

\medskip 

Of particular relevance is the question of compactness for a solutions sequence to \eqref{Liouville-Quasilinear}. 
In that direction, the seminal work of Brezis and Merle~\cite{BrMerle} in dimension two and its generalization by Aguilar and Peral in \cite{AgPe}  to the case $n \geq 3$ (see also~\cite{ReWe}) have shown that 
the asymptotic behavior of a solution sequence $u_k$ to
\begin{equation} \label{758}
-\Delta_n u_k=  h_k (x)  e^{u_k} \hbox{ in }\Omega,
\end{equation} 
with
\begin{equation} \label{758bis}
\sup_k \int_\Omega e^{u_k}<+\infty ,
\end{equation} 
can be  described by the following  ``concentration-compactness" alternative: 

\medskip
\goodbreak

{\bf Concentration-Compactness Principle:}
{\it Consider a sequence of solutions $u_k$ such that \eqref{758}-\eqref{758bis} hold with measurable functions $h_k$ satisfying
\begin{equation} \label{eq:Apperitivo}
    0 \leq h_k \leq C
    \, \, \text{ in } \Omega.
\end{equation} 
Then, up to a subsequence, the following alternative holds: 

 \begin{enumerate} \itemsep=5pt
\item[{\rm \bf (i)}] $u_k$ is bounded in $L^{\infty}_{loc} (\Omega)$;
\item[{\rm \bf (ii)}] 
 $u_k \to - \infty$ locally uniformly in $\Omega$ as $k\to +\infty$;
\item[{\rm \bf (iii)}] 
the blow-up set $\mathcal{S}$ of the sequence $u_k$, defined as
\begin{equation} \label{eq:Blowup-Set}
   \mathcal{S} = \{ p \in \Omega: \hbox{ there exists }  x_k \in \Omega \hbox{ s.t. } \, x_k \to p , u_k (x_k) \to \infty  \hbox{ as } k \to +\infty \},
\end{equation}
is nonempty, finite  with $u_k \to - \infty$ locally uniformly in $\Omega \setminus S$ and  
\begin{equation} \label{253}
   h_k e^{u_k} \xrightharpoonup{\, \, k \to \infty \, \,}
   \sum_{p \in {\mathcal S}} \beta_p \delta_p    
 \end{equation}
weakly in the sense of measures for some coefficients $\beta_p \geq  n^{n} \omega_n$, 
where $\omega_n$ stands for the volume of the unit ball in $\mathbb R^n$.
\end{enumerate} 
}

\medskip 

The compact case, in which the sequence $e^{u_k}$ does converge locally uniformly in $\Omega$, is expressed by alternatives (i) and (ii), thanks to elliptic estimates \cite{Dib,Tol}. 
Alternative (iii) describes the non-compact case, 
and the characterization of the possible values for the Dirac masses $\beta_p$ is crucial towards an accurate description of the blow-up mechanism.
In this matter a central role is played by the following invariance property of the PDE: starting from a function $u$ that solves \eqref{Liouville-Quasilinear}, then by introducing for a given 
$(\mu, x_0) \in (0,\infty) \times \mathbb R^n$ the rescaling functions 
\begin{equation} \label{eq:Zwei-Expressi}
    U (y) = u(\mu y + x_0) + n \log \mu
    \quad \hbox{ and } \quad
    H(y) := h(\mu y + x_0)
\end{equation}
defined in the domain $\Omega_{\mu} := \{ \mu^{-1} (x - x_0) \, : \, x \in \Omega\}$, we find that
\begin{equation} \label{eq:Rescaled}
   - \Delta_n U  = H(y) e^U
    \, \, \hbox{ in }  \Omega_{\mu}  ,
   \quad \hbox{ and } \quad 
   \int_{\Omega_{\mu}} e^{U}  = \int_{\Omega} e^u .
\end{equation} 
Now consider a sequence $(h_k,u_k)$ such that \eqref{758}, \eqref{758bis}, \eqref{eq:Apperitivo} hold. If the set  \eqref{eq:Blowup-Set} is not empty we can assume that
\begin{equation} \label{eq:Pastel1} 
  \qquad
  \overline{B_2 (p)} \subset \Omega, 
  \qquad 
  \overline{B_2 (p)} \cap {\mathcal S} = \{ p \},
  \quad \forall p \in {\mathcal S},
\end{equation}
and we always have a sequence $x_k \in B_2 (p)$ such that, as $k \to \infty$
\begin{equation} \label{eq:Pastel2} 
   \max_{\overline{B_2 (p)}} u_k = u_k(x_k) \to \infty, 
   \quad
   x_k \to p.
 \end{equation}
Hence, the rescaling function $U_k$ defined in \eqref{eq:Zwei-Expressi} with $\mu_k := e^{- \frac{u_k (x_k)}{n}}$ and center $x_k$ has the property
$\displaystyle \max_{B_{\mu_k^{-1}} (0)} U_k = U_k (0)  \leq 0$, and so by the concentration-compactness principle $U_k$ converges in $C_{loc} (\mathbb R^n)$ 
to some function $U: \mathbb R^n \to \mathbb R$. If in addition, as $k \to \infty$, it holds  
\begin{equation} \label{eq:Peso1}
        h_k \xrightarrow{\,  C_{loc} (\Omega) \,}  h, 
\end{equation}
the limiting function $U \leq U(0) = 0$ is then a solution of the problem
\begin{equation} \label{limitpb}
   -\Delta_n U=h(p) e^U \quad \hbox{ in }\mathbb{R}^n, \qquad \int_{\mathbb{R}^n} e^U<+\infty ,
\end{equation}
from which one can derive a more precise asymptotic behavior of the original blowing-up sequence $u_k$. To achieve this, one must classify the solutions to \eqref{limitpb} 
and take into consideration the possibility of several sharp peaks (local maximum points) collapsing to~$p$.

\medskip
\subsection{Two-dimensional results}

When the problem~\eqref{Liouville-Quasilinear} is considered with $n=2$,  the possible values for $\beta_p$ in \eqref{253} have been found by Y.Y.~Li and Shafrir \cite{LiSh}, who showed that under the 
convergence condition \eqref{eq:Peso1} then, for any $p \in \mathcal S$, $h(p)>0$ and the concentration mass $\beta_p$ is quantized as follows:
\begin{equation} \label{1011}
\beta_p \in 8\pi \mathbb{N}.
\end{equation}

\medskip  

The meaning of the value $8\pi$ in \eqref{1011} can be roughly understood follows: The blowing-up sequence $u_k$ develops several sharp peaks collapsing into $p$, each  of which looks like, after a proper rescaling, a solution $U$ of
\eqref{limitpb} with $h(p)>0$. Using the complex representation formula obtained by Liouville \cite{Liou} or the more recent PDE approach by Chen-Li~\cite{ChenLi} (see also \cite{CEL} for an alternative $P-$function approach), the solutions of~\eqref{limitpb} are explicitly known 
and they all have the same mass: $\int_{\mathbb R^2} h(p) e^U = 8 \pi$. Therefore the value of $\beta_p$ in \eqref{1011}  just represents the sum of the masses $8\pi$ carried by each of such sharp peaks collapsing in $p$. 

\medskip

However, to make this argument precise, one needs to ensure that
the mass contribution coming from the neck regions between the sharp peaks vanishes in the limit.  
This crucial property relies on a $\sup + \inf$ inequality,  an inequality of Harnack type that has been first derived by Shafrir~\cite{Shafrir}. It states that
 if one considers the set
$$\Lambda_{a,b} :=\{h  \in L^{\infty} (\Omega): \ a \leq h  \leq b \hbox{ in } \Omega\},$$
and a compact subset $K \subset \Omega$, then there are constants $C_1 = C_1 (\frac{b}{a})$ and $C_2 = C_2 (a,b, K) $ so that 
\begin{equation} \label{0317}
\sup_K u + C_1 \inf_{\Omega} u \leq C_2
\end{equation}
holds for any solution $u$ to \eqref{Liouville-Quasilinear}$_{n=2}$ with $h \in \Lambda_{a,b}$, where hereafter $0 < a \leq b$ are fixed.
This result was established using an isoperimetric argument together with an analysis of the mean average for a solution $u$ to \eqref{Liouville-Quasilinear}$_{n=2}$.
In addition, in the case $h\equiv 1$, by making use of the Liouville representation formula, Shafrir also showed that the inequality \eqref{0317} holds with $C_1=1$.  
Later, based on a blowup argument combined with the moving plane/sphere technique, Brezis, Li and Shafrir \cite{BrLiSh} showed the validity of \eqref{0317} 
in its sharp form with $C_1=1$ provided $h$ belongs to the set:
\begin{equation} \label{1100}
  \Lambda_{a,b}'=\{h  \in {\rm Lip} (\Omega) \, : \,  a\leq h \leq b,\ |\nabla h|  \leq b  \hbox{ in }\Omega\}.
\end{equation}
Relying again on some isoperimetric arguments, in  \cite{ChenLin}  Chen and Lin  showed 
the validity of this result for a much wider class of functions $h$. In particular, in lieu of \eqref{1100},  one can consider a set of H\"older continuous functions with a bound on the H\"older semi-norm. 
An alternative proof of this type of optimal $\sup + \inf$ inequality can also be obtained through Green's representation formula  by following some arguments in \cite{Rob}, and we also mention the contribution of 
\cite{ChenLi1}.

\medskip

Under some bound requirements on the boundary oscillation, Y.Y. Li obtained in \cite{YYLi-Harnack}  
a more precise asymptotic description of a blowing-up solutions sequence $u_k$ to \eqref{758}, which as a matter of fact strengthens the optimal $\sup + \inf$ inequality. More specifically,  assuming $h_k \in \Lambda_{a,b}'$ along with \eqref{253} and
\begin{equation} \label{eq:TaschenTuch-Mouchoir2}
  {\rm osc_{\partial \Omega} }u_k  := \sup_{x,y \in \partial \Omega} |u_k(x)-u_k(y)|\leq C  ,
\end{equation}
then in a ball $B_2 (p)$ as in \eqref{eq:Pastel1} and along a sequence of maximum point $x_k \in B_2 (p)$ as in \eqref{eq:Pastel2}, one has  (see \cite[Theorem 0.3]{YYLi-Harnack}) 
\begin{equation} \label{eq:DasDreieck}
   \Big| u_k - \log \Big( \frac{e^{-u_k(x_k)}}{(e^{-u_k(x_k) } + \frac{ h_k(x_k)}{8}   |x - x_k|^2 )^2} \Big)  \Big|  \,  \leq \, C ,
   \quad \hbox{ in } B_1 (p) .
\end{equation}
From \eqref{eq:DasDreieck} one can easily derive the following estimate
\begin{equation} \label{eq:Zeit-Ist-Geld}
    - C \, \leq  \, u_k (x_k) +  \inf_{B_2 (p) } u_k  \, \leq \, C ,
\end{equation}
which is more accurate than the optimal estimate coming from \eqref{0317}. This was the initial step of a larger program that in the two-dimensional case has been accomplished by Chen and Lin \cite{ChenLin1,ChenLin2}.

\medskip

The goal of this paper will be to extend the profile decomposition \eqref{eq:DasDreieck} and the optimal $\sup + \inf$ inequality when the operator in \eqref{Liouville-Quasilinear} is 
replaced by the $n-$Laplacian with $n \geq 3$.  However, due to the nonlinearity of the differential operator, we face many obstacles.
For instance, the Green representation formula is not longer available in the quasilinear context, and 
its alternative  in terms of the Wolff potential (see \cite{KiMa}) fails to provide sharp constants as needed to derive precise asymptotic estimates on blowing-up solutions to~\eqref{758}-\eqref{758bis}.  
We refer the interested reader to \cite{HKM,KuMi} for an overview on the nonlinear potential theory. In the present paper we overcome those difficulties through a radically new approach.

\medskip
\subsection{Blow-up results for the quasilinear problem}

Like in dimension two the blowup mechanism  for the solution to  \eqref{Liouville-Quasilinear} is tightly related to the scaling and translation invariance described in  \eqref{eq:Zwei-Expressi} leading 
to the limiting problem \eqref{limitpb}, for which the classification of the solutions has been recently provided by the first author in \cite{Esp}, see also \cite{CEL,Esp1}. 
In particular, the unique solution to \eqref{limitpb} with $U\leq U(0) = 0$  is given by the ``standard bubble"
\begin{equation} \label{eq:Bubble}
   U(y)=-n \log \left(1+c_n^{-\frac{1}{n-1}} h(p)^{\frac{1}{n-1}}|y|^{\frac{n}{n-1}}\right) ,
\end{equation} 
and from the invariance described by \eqref{eq:Zwei-Expressi}, \eqref{eq:Rescaled}, all the solutions to \eqref{limitpb} are
\begin{equation} \label{eq:Classification-PP}
   \qquad 
   U_{\mu,x_0}  (y) =  \log   \frac{ \mu^\frac{n}{n-1}}{ ( \mu^\frac{n}{n-1} +  (\frac{h(p)}{c_n})^{\frac{1}{n-1}}  |y - x_0 |^{\frac{n}{n-1}} )^n }  ,\quad (\mu,x_0)\in (0,+\infty)\times \mathbb{R}^n,
\end{equation}
an $(n+1)$-parameter family which satisfies
\begin{equation} \label{quantization}
    \qquad
    \int_{\mathbb{R}^n} h(p) e^{U_{\mu,x_0}}  =c_n \omega_n , \quad c_n :=n \big( \frac{n^2}{n-1} \big)^{n-1}.
\end{equation}

As a by-product, under the condition \eqref{eq:Peso1} we necessarily have in \eqref{253} that  $\beta_p \geq c_n \omega_n $, 
a bigger value than the one appearing in the alternative (iii) of the Concentration-Compactness Principle.

\medskip

When a boundary control on $u_k$ is assumed, the possible value of $\beta_p$ in \eqref{253} are easier to be determined.
More specifically, in \cite{EsMo} it has been shown 
that if the oscillation of $u_k$ on $\partial B_{\delta} (p)$, $p \in \mathcal S$, is uniformly bounded for some $\delta>0$, then an application of a Pohozaev identity leads
to $\beta_p = c_n \omega_n $, provided $h_k \in C^1 (\Omega) \cap \Lambda_{a,b}'$ (defined in \eqref{1100}). In the general case, one may take advantage of the scaling invariance property 
\eqref{eq:Zwei-Expressi}-\eqref{eq:Rescaled}   which in dimension two is the key reason why each peak of a blowingup sequence carries the same mass and leads to \eqref{1011}. Since this invariance holds for any $n \geq 2$, a 
similar quantization property was expected for the $n$-Liouville problem \eqref{Liouville-Quasilinear}, and indeed the authors of the present paper established in \cite{EsLu}
\begin{equation} \label{10111}
    \beta_p \in c_n\omega_n \mathbb{N},
\end{equation}
whenever \eqref{253}, \eqref{eq:Peso1} and $h_k \in \Lambda_{a,b}$ hold.
To reach this result we needed a decay estimate of the sequence 
$h_k e^{u_k}$ on the neck regions between the peaks of the sequence. This was made possible thanks to 
a $\sup+\inf$ similar to  \eqref{0317} that we 
proved  in \cite{EsLu}, using a new blow-up approach, and it reads as follows  

\medskip

\begin{thm} \label{thm:Non-Optimal} 
Given $0<a\leq b <\infty$, consider a family $\Lambda \subset \Lambda_{a,b}$ which is equicontinuous at each point $x \in \Omega$ and let 
\begin{equation} \label{eq:Admissible-Sol}
   {\mathcal U} := 
      \big\{ u \in C^{1,\alpha} (\Omega) \cap C(\overline{\Omega}) \, : \, u \hbox{ satisfies } \eqref {Liouville-Quasilinear} \hbox{ with }  h \in \Lambda  \big\}.
\end{equation}

Given a compact set $K \subset \Omega$, then for each $C_1>n-1$ there exists $C_2 =C_2( \Lambda, K, C_1) >0$  such that 
\begin{equation} \label{eq:NonOptimal}
         \qquad
         \sup_K u + C_1 \inf_{\Omega} u \leq C_2,
         \quad \forall u \in {\mathcal U} .
\end{equation}
\end{thm}

\medskip

In the present paper, we push further the analysis undertaken in \cite{EsLu} and aim to have a profile description of a solutions sequence  $u_k$ to the $n$-Liouville equation that blows up at a point $p \in \Omega$ 
chosen without loss of generality to be the origin $0$. Our set up will be a sequence $(h_k,u_k)$ such that
\begin{equation}  \label{1553} 
    \left.
    \begin{array}{c}
        (h_k, u_k) \hbox{ satisfies }       \eqref{758} \hbox{ and }\eqref{758bis}  ,
        \vspace{2mm} \\
                h_k \in \Lambda_{a,b} \hbox{ satisfies } \eqref{eq:Peso1},
        \vspace{2mm} \\
         \displaystyle  \max_{\overline{B_2 (0)}} u_k = u_k (x_k) \to \infty, 
    \quad   
    x_k \to 0 \in \overline{B_2(0)} \subset \Omega, 
            \vspace{2mm} \\
               u_k \hbox{ is bounded in } L^{\infty}_{loc} (\overline{B_2(0)}) \setminus \{0 \} ,
             \vspace{2mm} \\
 \end{array}
\right\}
\end{equation}
where this last condition can be assumed without loss of generality thanks to the Concentration-Compactness Principle. In addition to \eqref{1553},  we will also assume \eqref{eq:TaschenTuch-Mouchoir2}.

To derive our results, we will first combine BMO bounds with Serrin's elliptic estimate to derive some general apriori estimates for the $n$-Poisson equation $-\Delta_n u =f$. 
By applying those results to the $n$-Liouville equation, we will obtain oscillation estimates outside the blowing-up regions of a bubbling sequence. One main difficulty arises from the fact that the sequence $u_k$ may develop, 
aside the maximum points $x_k$, other peaks given by local maxima along which $u_k$ tends to infinity.  Focusing first on a bubbling solutions sequence that develops only one peak, a main asymptotic estimate of the present paper 
is given by 
\begin{equation} \label{eq:MusicaTunis}
   u_k  - \inf_{\partial B_{\mu_k} (x_k)} u_k \, = \,  
   -  \Big( \frac{1}{n \omega_n} \int_{B_{\frac{1}{2}}(x_k) } h_k e^{u_k} \Big)^{\frac{1}{n-1}} \log \Big( \frac{|x-x_k|}{\mu_k} \Big)  +O(1) ,
\end{equation} 
where $\mu_k := e^{- \frac{u(x_k)}{n} }$.  By combining this result with a Pohozaev identity,  we are able to extend the profile decomposition \eqref{eq:DasDreieck} to our quasilinear setting as follows:

\begin{thm} \label{thm:Fundamental-Estimate}
Assume that \eqref{eq:TaschenTuch-Mouchoir2}, \eqref{1553} hold and $\| \nabla h_k \| \leq C$. Then, 
\begin{equation} \label{eq:Variable2}
   \Big| u_k - \log   \frac{ e^{-\frac{u_k(x_k)}{n-1}} }{ ( e^{-\frac{u_k(x_k)}{n-1}}  +   (\frac{h(p)}{c_n})^{\frac{1}{n-1}}  |x - x_k |^{\frac{n}{n-1}} )^n }  \Big| \leq
  C,
   \quad \hbox{ in } B_1 (p).
\end{equation}
\end{thm}

\medskip

Notice that our result covers the two-dimensional formulation \eqref{eq:DasDreieck}, since
\eqref{eq:Variable2} can also be written with $h(p)$ in place of $h_k(x_k)$, a change that would only add a bounded term in the inequality.
From \eqref{eq:Variable2} we obtain an analogue result to \eqref{eq:Zeit-Ist-Geld}, namely we can see that the blowing up sequence $u_k$ satisfies the following inequality of $\sup + \inf$ type
\begin{equation} \label{eq:Zeit-Ist-Geld-2}
    - C  \, \leq \, u_k (x_k) +  (n-1) \inf_{B_2 (p) } u_k  \, \leq \, C .
\end{equation}
In light of this latest estimate, we will look into the borderline case $C_1 = n-1$ of the inequality~\eqref{eq:NonOptimal} and exhibit a class of  pairs $(h,u)$ for which it holds.  
Our result will not require any $C^1$-bound on the class of admissible functions $h$.

\medskip
\subsection{Optimal quasilinear $\sup+\inf$ inequality}

The inequality \eqref{eq:NonOptimal} was established in \cite{EsLu} and the question of its validity with sharp constant $C_1 = n-1$ was left unanswered.  
That this should be the case is supported by the fact that  the functions $U_{0,\lambda}$ given by \eqref{eq:Classification-PP} solve
\eqref{Liouville-Quasilinear} in $B_{1} (0)$ with $h \equiv 1$,  and satisfy
$$
  \max_{\{ 0 \}} U_{0,\lambda} + \alpha \inf_{B_1 (0)} U_{0,\lambda}
   \, = \, 
   \log ( \lambda^n) + \alpha \log \big(  \frac{\lambda^n}{\lambda^{\frac{n^2}{n-1}} } \big)  + O(1)
   \, = \, 
   n ( 1 - \frac{\alpha}{n-1}) \log \lambda + O(1)
$$
which is bounded from above for $\lambda$ large if and only if $\alpha \geq n -1$.  

\medskip

The second main result of our paper shows that indeed inequality \eqref{eq:NonOptimal} holds with $C_1 = n-1$ if we work as in \cite{ChenLin} with the class of functions $h$ that have a uniform logarithmic modulus of continuity at each point:
$$
   \Lambda^{{\rm Log}}_{a,b} := 
   \Big\{ h \in L^{\infty} (\Omega) \, : \, a \leq h \leq b, \quad |h(x) - h(y)| \leq \frac{b}{\big| \log |x - y | \big| } \Big\}, 
$$
and by considering for each $M >0$ the set of solutions to \eqref {Liouville-Quasilinear}  that satisfy for some $M >0$ 
\begin{equation} \label{1228} 
   \hbox{ osc}_{\partial \Omega} (u)\leq M,
   \qquad
   \int_{\Omega} e^u \leq M .
 \end{equation}

Observe that 
$ \Lambda_{a,b}' \subset   \Lambda^{{\rm Log}}_{a,b} \subset \Lambda_{a,b}$ and that $ \Lambda^{{\rm Log}}_{a,b} $ is a family 
which is equicontinuous at each point of~$\Omega$ (thus Theorem \ref{thm:Non-Optimal} can be applied with this family). 
Our second main result is the following optimal $\sup + \inf$ inequality, to be compared with \cite{BCLT}:

\begin{thm} \label{711intro} 
Given $0<a\leq b <\infty$ and $M>0$, consider
$${\mathcal U}_M := 
   \big\{ u \in C^{1,\alpha} (\Omega) \cap C(\overline{\Omega}) \, : \, u \hbox{ satisfies } \eqref {Liouville-Quasilinear}  \hbox{ and } \eqref{1228}  \hbox{ with }  h \in   \Lambda^{{\rm Log}}_{a,b}  \big\}.$$
Then, for each  given compact set $K \subset \Omega$  there exists $C_2 =C_2(a,b,K,M) >0$ so that
\begin{equation}\label{317bis}
        \max_K u+(n-1) \inf_\Omega u \leq C_2
        \quad
        \forall u \in {\mathcal U}_M.
\end{equation}
\end{thm}
To undertake this task, we will argue by contradiction and consider a sequence of solutions that fails to satisfy \eqref{317bis}. 
Such a sequence must necessarily blow-up, and in full generality  might develop several peaks of local maxima.  
Since after proper rescaling the limiting profile of such a solutions sequence $u_k$ has the form \eqref{eq:Bubble}, our strategy will be 
to use the sharp constants available in the radial case, by comparing effectively $u_k$ with the bubble \eqref{eq:Bubble} through  
the new concept of {\it nonlinear spherical average} that will be presented in Section \ref{sec:Spherical-Mean}. 
In dimension $n=2$, this notion coincides with the usual  spherical average $\bar u = \fint_{\partial B_r } u$, which satisfies $\Delta \bar u=\fint_{\partial B_r } \Delta u$, a fact that can be exploited 
through an ODE analysis  to obtain information on $u$ as long as $u-\bar u$ is controlled, as this is usually the case in the blow-up regime. 
When $n \geq 3$, we extend this by  considering the unique radial solution $\tilde u$ solving
$$  
       \Delta_n \tilde u=\fint_{\partial B_r}(\Delta_n u) \quad \hbox{ for }0<r<R, \quad \tilde u (R)=\bar u(R) ,
$$
a function that will be called  ``nonlinear spherical average of $u$".  In contrast with $n=2$, this defines an operator that does not commute with the $n$-Laplacian.
Nevertheless, since $g=\Delta_n u-\fint_{\partial B_r}(\Delta_n u)$ has vanishing radial mode, the one that is typically responsible for the wilder behavior, 
one expects better estimates for $u -\tilde u$ through elliptic estimates on the difference equation $\Delta_n u-\Delta_n \tilde u=g$. 
Though we cannot achieve this in a full generality, we have successfully obtained some type of improved estimates in a blow-up regime, where good estimates on $g$ and $u-\bar u$ are available. 

\medskip 

Deriving a general regularity theory for $u-\tilde u$, or proving the inequality \eqref{317bis} without assumption \eqref{1228} are open and challenging problems.

\medskip 

{\bf Acknowledgments.}
Pierpaolo Esposito acknowledges support by INDAM-GNAMPA and by the project
PNRR-M4C2-I1.1-PRIN 2022-PE1-Variational and Analytical aspects of Geometric
PDEs-F53D23002690006 - Funded by the U.E.-NextGenerationEU. 
 Marcello Lucia has been supported by the Simons Foundation MP-TSM-00002452.

\smallskip

The authors confirm that no artificial intelligence tools or machine learning models were used to elaborate the proofs
or write this manuscript.

\bigskip \goodbreak
\section{A priori estimates for the $n$-Poisson equation} 
\medskip 

Throughout this section we aim to give a priori estimates for the set of weak solutions to  
\begin{equation} \label{eq:XiAn}
   - \Delta_n u = f ,
   \qquad
   u \in W^{1,n} (\Omega), 
\end{equation}
where  $f \in L^{1}(\Omega)$ is defined on a bounded domain $\Omega \subset \mathbb R^n$ of class $C^1$. Here, ``weak solution"  refers to $u \in W^{1,n} (\Omega)$ 
satisfying
\begin{equation} \label{1116}
      \qquad
       \int_{\Omega} |\nabla u|^{n-2} \nabla u \nabla \varphi = \int_{\Omega} f \varphi ,
            \qquad
       \forall \varphi \in W^{1, n}_0 (\Omega) \cap L^{\infty}(\Omega),
\end{equation}
and we  observe that by density \eqref{1116} holds for all $\varphi \in W^{1, n}_0 (\Omega)$ whenever  $f \in L^{\infty} (U)$.

\medskip

In the coming subsections, we will present a Caccioppoli type inequality together with some $W^{1,n}$ estimates, and  by following the arguments in \cite{DHM} we will establish 
uniform BMO a priori bounds for the set of solutions of~\eqref{eq:XiAn}, without assuming any boundary conditions. We will conclude with several useful oscillation estimates.

\medskip
\subsection{Interior gradient estimates}

In many instances we will derive a priori bounds for solutions to \eqref{eq:XiAn} by considering an appropriate comparison problem $-\Delta_n u_0 = f_0$ combined with suitable test functions.
More specifically, if on a bounded open set $U \subset \mathbb R^n$  we look at the two Poisson problems
\begin{equation} \label{eq:Comparison-Prob}
   - \Delta_n u  \, = \, f,
   \qquad
   - \Delta_n u_0 \, = \, f_0,
\end{equation}
with $f, f_0 \in L^{\infty} (U)$, by taking the difference we get
\begin{equation} \label{eq:Ein-Kleiner-Espresso}
   - \Delta_n u  + \Delta_n u_0 \, = \,  f - f_0  .
\end{equation} 

\medskip

Given $\eta \in W_0^{1,\infty}(U)$ and $g \in {\rm Lip}  (\mathbb R)$, 
 we can choose  $\eta^n g (u  - u_0 ) \in  W_0^{1,n}(U)$ as a legitimate test function in the weak formulation of \eqref{eq:Ein-Kleiner-Espresso} (since $f - f_0 \in L^{\infty} (U)$). This leads to
$$
   \int_{U} \big\langle a(\nabla u) - a(\nabla u_0) , \nabla [\eta^n g(u - u_0) ] \big\rangle \, = \, \int_{U} (f - f_0) \, \eta^n g(u - u_0)  ,
$$
where $a(X)$ denotes the vector field $a (X) := |X|^{n-2} X$, $X \in \mathbb R^n$.We then get

\begin{eqnarray}
&& 
     \int_{U} \eta^n g'(u-u_0) \big\langle a(\nabla u) - a(\nabla u_0) , \nabla (u - u_0)  \big\rangle
     \nonumber \\
     && \qquad \leq   \, 
     n  \Big| \int_{U}  \eta^{n-1} g(u - u_0) \big\langle a(\nabla u) - a(\nabla u_0) , \nabla \eta \big\rangle \Big|  \, + \, \int_{U} (f - f_0) \, \eta^n g(u  - u_0)  .
    \qquad
    \label{eq:General-Caccioppoli}
\end{eqnarray} 

\medskip

If in the comparison problem \eqref{eq:Comparison-Prob} we consider the simplest case where 
$u_0 \equiv cst $ and $f_0 \equiv 0$, the inequality \eqref{eq:General-Caccioppoli} leads to some Caccioppoli type estimates:
\begin{pro} \label{stimaC}
For any weak solution $u$ of \eqref{eq:XiAn} we have
\begin{equation} \label{Caccio}
        \int_{\{|u-c|\leq M\} \cap B_{r_1}(y)}|\nabla u|^n \, \leq \, (\frac{4n^2}{r_2-r_1})^n \int_{B_{r_2}(y)\setminus B_{r_1}(y)}|u-c|^n+2n M \int_{B_{r_2}(y)}|f|
\end{equation}
for all $M>0$, $c \in \mathbb{R}$, $0<r_1<r_2 <\hbox{dist}(y,\partial \Omega)$.
\end{pro} 
\begin{proof} 
Consider a cut-off function $\eta \in C_0^{\infty}(B_{r_2}(y))$ so that 
\begin{equation} \label{eq:CutOff}
   0  \, \leq \, \eta  \, \leq \, 1,  \quad \eta=1 \hbox{ in }  B_{r_1}(y), 
   \quad 
   (r_2-r_1)|\nabla \eta|\leq 2 ,
 \end{equation}
 and consider also the $C^1$-function
$$
   g_M(s)=\left\{ \begin{array}{ll} s &\hbox{ if } |s| \leq M 
   \vspace{2mm} \\ 
   \left[n \frac{M}{|s|}-(n-1) (\frac{M}{|s|})^{\frac{n}{n-1}}\right] s & \hbox{ if }|s|>M \end{array} \right.
$$
which satisfies 
\begin{equation} \label{1745}
    |g_M| \leq n M, \quad 0\leq g_M'\leq 1, \quad 0\leq \frac{g_M(s)}{n s} \leq [g_M'(s)]^{\frac{n-1}{n}}.
\end{equation}

\medskip

From \eqref{eq:General-Caccioppoli}  with $u_0 \equiv c $ and $f_0 \equiv 0$, we get
\begin{eqnarray*}
&& 
\int_\Omega   |\nabla u|^n  \eta^n g_M'(u-c)  
\leq n  \int_\Omega   |\nabla u|^{n-1}   \eta^{n-1} | \nabla \eta| |g_M(u-c)| +  \int_\Omega  |f|  \, \eta^n |g_M (u -c)|
 \\  [9pt]
&&
\stackrel{\eqref{1745}}{\leq}
n^2  \int_\Omega   |\nabla u|^{n-1}  \eta^{n-1} |  \nabla \eta |  \, | |u-c| [g'_M(u-c)]^{\frac{n-1}{n}} 
+  n M  \int_\Omega |f| \,  \eta^n  
  \\  [9pt]
&& \leq \frac{1}{2}    \int_\Omega    |\nabla u|^n \eta^n  g'_M(u-c) +2^{n-1} (n-1)^{n-1} n^n     \int_\Omega   |\nabla \eta|^n |u-c|^n +n M  \int_\Omega  |f| \, \eta^n 
\end{eqnarray*}
where the last step follows from $ab \leq \frac{n-1}{n}a^{\frac{n}{n-1}}+\frac{b^n}{n}$, $a,b>0$. 
Hence, by rearranging the latter inequality and observing that $g'(s) = 1$ for $|s| \leq M$, we obtain 
$$
   \int_{\{|u-c|\leq M\} } |\nabla u|^n \eta^n 
\leq
2^{n} n^{2 n}     \int_\Omega  |u-c|^n  |\nabla \eta|^n  + 2 n M  \int_\Omega   |f| \, \eta^n .
$$
Finally, by using the properties of $\eta$ given in \eqref{eq:CutOff}, we obtain the inequality claimed in \eqref{Caccio}.
\end{proof}

\medskip

We will later exploit \eqref{eq:General-Caccioppoli} with some non-constant $u_0$, and in that case 
it is convenient to recall that the vector field 
$a(X) $ satisfies for all $X,Y \in \mathbb R^n$ the following:

\begin{eqnarray} \label{eq:KokosMilch-Or}
   | a(X) - a(Y) | &\leq & C \big(|X|^{n-2} + |Y|^{n-2} \big) |X - Y|  ,\\
  \langle a(X) - a(Y), X - Y \rangle 
  &=& 
\frac{1}{2} \big(|X|^{n-2} + |Y|^{n-2} \big) |X - Y|^2  \nonumber \\
 && \quad + \,  \frac{1}{2} \big( |X|^{n-2} -  |Y|^{n-2}  \big) \big( |X|^{2} -  |Y|^{2} \big)
  \label{eq:Coercivity-Or-Noise} \\ 
  &\geq&   
  \frac{1}{2} \big(|X|^{n-2} + |Y|^{n-2} \big) |X - Y|^2 .  \nonumber \end{eqnarray}
Indeed, inequality \eqref{eq:KokosMilch-Or} follows by using the fundamental theorem of calculus, while equality \eqref{eq:Coercivity-Or-Noise} is readily obtained after expanding each term. 
We will use a direct variant of these inequalities that reads as follows

\begin{lm} \label{lem:Meeting-Monday}
For some $C,\delta  >0$, the following hold for all $X,Y \in \mathbb R^n$
\begin{equation} \label{eq:KokosMilch}
   | a(X) - a(Y) | \, \leq \, C \big(|X - Y|^{n-2} + |Y|^{n-2} \big) |X - Y|  ,
\end{equation}
\begin{equation} \label{eq:Coercivity}
  \langle a(X) - a(Y), X - Y \rangle \, \geq \,  \delta \big(|X -Y |^{n-2} + |Y|^{n-2} \big) |X - Y|^2 .
\end{equation}
\end{lm}

\medskip

With the lower bound \eqref{eq:Coercivity}, the general inequality \eqref{eq:General-Caccioppoli} in the case $g(s) = s$ leads to
\begin{eqnarray}
&& 
      \delta \int_U  \underbrace{\eta^n  \Big\{|\nabla (u - u_0)|^{n-2}+|\nabla u_0|^{n-2} \Big\} |\nabla (u - u_0) |^2}_{I_1}      \nonumber \\
     && \qquad \leq   \, 
     \Big|  \int_{U}  
      \underbrace{n \eta^{n-1}  \big\langle a(\nabla u) - a(\nabla u_0) , \nabla \eta \big\rangle}_{I_2}   (u - u_0)  \Big|  \, + \, \int_{U} (f - f_0) \, \eta^n (u  - u_0)  .
    \qquad
    \label{eq:General-Caccioppoli-Bis}
\end{eqnarray} 

We can get useful estimates by combining suitably in \eqref{eq:General-Caccioppoli-Bis} the term $I_2$ with $I_1$. We first stress two properties of $I_2$:

\begin{lm} \label{lem:Scusa}
Consider $I_1$, $I_2$ as in \eqref{eq:General-Caccioppoli-Bis}, then 
for any $w: U \to \mathbb R$, we have 
\begin{equation} \label{eq:Heiss} 
   w I_2  \,  \leq \,  \frac{\delta}{2} \, I_1 + \,   C_2 \Big\{  |\nabla \eta|^n |w|^n +  |\nabla u_0|^{n-2} | \nabla \eta|^2  |w|^2  \Big\}   .
\end{equation}

Furthermore, given $B_r (p_0) \Subset U$ and $\eta \in C^1 (U)$ which is radial $\eta(x) = \eta (|x - p_0|)$ we have 
\begin{equation} \label{eq:Heiss2}
    \int_{\partial B_r (p_0) } I_2  =   (\eta^n)' (r)   \int_{B_r (p_0)}   ( f_0 - f) 
\end{equation}
\end{lm} 
\begin{proof}
From inequality \eqref{eq:KokosMilch} of Lemma \ref{lem:Meeting-Monday} we obtain
\begin{eqnarray}
   w I_2 &=& 
   w  n \eta^{n-1} \, \big\langle a(\nabla u) - a (\nabla u_0) , \nabla \eta \big\rangle  
   \nonumber \\
   &\leq& 
   C_1 w  \eta^{n-1} \,  \big(  |\nabla (u - u_0)|^{n-2} +   |\nabla u_0 |^{n-2}  \big) |\nabla (u-u_0)| |\nabla \eta|   .
   \label{eq:Der-Verzehr}
\end{eqnarray}

Since by Young's inequality there holds
\begin{equation} \label{eq:Gitarre2}
  C_1  |\nabla (u - u_0)|^{n-1} \eta^{n-1}  |\nabla \eta|  w \, \leq \,  \frac{\delta}{2} \, \eta^n   |\nabla (u - u_0)|^{n} + C_2    |\nabla \eta|^n  w^n  ,
\end{equation}
\begin{equation} \label{eq:Gitarre}
   C_1  |\nabla u_0 |^{n-2}  |\nabla (u-u_0)|   \eta^{n-1}  |\nabla \eta|  w \, \leq \,  |\nabla u_0 |^{n-2}  \Big\{  \frac{\delta}{2} \eta^n |\nabla (u-u_0)|^2 + C_2 \eta^{n-2}  |\nabla \eta|^2  w^2  \Big\} ,
\end{equation}

\medskip

by plugging \eqref{eq:Gitarre2} and \eqref{eq:Gitarre} into \eqref{eq:Der-Verzehr}, we obtain the conclusion \eqref{eq:Heiss}.

\medskip

Concerning \eqref{eq:Heiss2}, 
denoting by $\nu$ the exterior unit normal vector field to $\partial B_r (p_0)$, Stokes' Theorem gives
$$
     \int_{\partial B_r (p_0) } I_2  \, = \, 
    (\eta^n)' (r) \int_{\partial B_r (p_0)}  \big\langle a(\nabla u) -     a( \nabla u_0) , \nu \big\rangle
    \, = \,
    (\eta^n)' (r) \int_{ B_r (p_0)}  ( \Delta_n u -    \Delta_n u_0 ) $$
and therefore we get \eqref{eq:Heiss2}.
\end{proof} 

\bigskip

By substituting the inequality \eqref{eq:Heiss} with $w = u - u_0$ in \eqref{eq:General-Caccioppoli-Bis}, we immediately obtain the following type of Caccioppoli inequality. 

\begin{pro} \label{prop:Benutzen}
For some universal constant $C >0$, the following holds
\begin{eqnarray*}
  && 
  C  \int_U  \eta^n  \Big\{|\nabla (u - u_0)|^{n-2}+|\nabla u_0|^{n-2} \Big\} |\nabla (u - u_0) |^2  \\
  && \qquad\, \leq \,
\Big| \int_{U} \eta^n (f_0 - f) (u - u_0)  \Big| + \int_U \Big\{ |\nabla \eta|^n | u-u_0|^n + |\nabla u_0|^{n-2} |\nabla \eta|^2 (u-u_0)^ 2 \Big\}
\end{eqnarray*}
for all $(u,f)$, $(u_0, f_0)$ satisfying  \eqref{eq:Comparison-Prob} and for all $\eta \in W^{1,\infty}_0 (U)$.
\end{pro} 

\medskip

Previous proposition admits the following generalization, which adds some flexibility to the choice of test functions. 
\begin{pro} \label{prop:Benutzen2}
 In a ball $B_{R_2} (p_0)$ consider
$(u,f)$, $(u_0, f_0)$ satisfying~ \eqref{eq:Comparison-Prob} and
any pair of functions $(\eta,  {\mathcal R} ) \in C^1_0 (B_{R_2} (p_0)) \times C^1 (B_{R_2} (p_0)) $ which are radial with respect to $p_0$.  
Then, 
by considering $A := B_{R_2} (p_0) \setminus \overline{B_{R_1} (p_0)}$ ($0 \leq R_1 < R_2$) and letting $w := u - u_0 - {\mathcal R}$, we have
\begin{eqnarray}
  && C  \int_A   \eta^n  \Big\{|\nabla (u - u_0)|^{n-2} + |\nabla u_0|^{n-2} \Big\} |\nabla (u - u_0) |^2  
   \leq  \int_A \Big\{ |\nabla \eta|^n | w|^n + |\nabla u_0|^{n-2} |\nabla \eta|^2 w^ 2 \Big\}
\nonumber \\
  &&
  \qquad \qquad \qquad \quad 
   \, + \,
  \Big|  \int_{R_1}^{R_2}  \underbrace{\int_{B_s (p_0) }  {\mathcal R}  (\eta^n)'  ( f_0 - f) dx \, ds }_{I_0} \Big|+   \Big|\int_{A} ( f - f_0) \eta^n  (u-u_0) \Big|
  \quad
  \label{eq:Ich-Stoere-Sie}
\end{eqnarray}
for some universal constant  $C >0$.
\end{pro}
\begin{proof}
In  \eqref{eq:General-Caccioppoli-Bis}, we write 
$\int_{A} I_2  (u-u_0) =      \int_{A} I_2 w  + \int_{A} I_2  {\mathcal R} $.  
The term $ \int_{A} I_2 w$ can be estimated using \eqref{eq:Heiss}, and by exploiting the radial symmetry of both functions 
$\eta, {\mathcal R}$  we have 
$$
  \int_{A} I_2  {\mathcal R}   =  \int_{R_1}^{R_2}{\mathcal R}  \int_{\partial B_s (p_0)} I_2  d\sigma ds 
 \stackrel{\eqref{eq:Heiss2}}{=}   \int_{R_1}^{R_2}{\mathcal R}  (\eta^n)'\int_{B_s (p_0)} (f_0 -f) ds. 
$$
Using those facts in \eqref{eq:General-Caccioppoli-Bis}, we then deduce the claimed estimate \eqref{eq:Ich-Stoere-Sie}.
\end{proof}

\medskip

Notice that when $( f - f_0)$ has zero spherical average, by using the fact that $\eta^n {\mathcal R}$ is a radial function, we see that the integral $I_0$ in \eqref{eq:Ich-Stoere-Sie}
vanishes. This will be exploited in Section~\ref{sec:Spherical-Mean}.

\medskip

The Caccioppoli estimates obtained in Prop. \ref{prop:Benutzen} and Prop.~\ref{prop:Benutzen2} will be applied to get gradient bounds on some annuli $A' \subset A \subset \Omega$ 
centered at $p_0$ with a term $|\nabla u_0|  \leq \frac{C}{|x - p_0|}$, and a class of function $\| w \|_{L^{\infty} (\Omega)} \leq C$. In that situation, we will choose 
a function $\eta$ that satisfies
$$
    0 \leq \eta \leq 1, \quad \eta \in C^1_0 (A), 
    \quad \eta \equiv 1 \hbox{ on } A',  
$$
which  minimizes the problems
$$
    \inf \Big\{ \int_{A^{\prime}} \big| \nabla \xi \big|^n \, : \, \xi \in W^{1,n}_0 (A) , \, \xi \equiv 1 \, \, {\rm on } \, \, A^{\prime} \Big\} ,
$$
$$
    \inf \Big\{ \int_{A^{\prime}} \frac{\big| \nabla \xi \big|^2}{|x - p_0|^{n-2}} \, : \, \xi \in W^{1,n}_0 (A) , \, \xi \equiv 1 \, \, {\rm on } \, \, A^{\prime}\Big\} .
$$
By applying the direct method in calculus of variation, we can see that such a function exists and 
is given by the radial function $\eta (x) :=\eta(|x- p_0|)$  that solves on $A \setminus {\overline A^{\prime}}$ the respective Euler-Lagrange equation
$$
 \begin{array}{c}
    \Delta_n \eta = 0,
\end{array}
\qquad
   \begin{array}{cc}
   {\rm div} (\frac{ \nabla \eta }{|x - p_0|^{n-2}}) = 0,
\end{array}
$$
subjected to the boundary condition $\eta|_{\partial A} = 0$ and $\eta |_{\partial A^{\prime}} = 1$. 
The solutions, that can be computed explicitly, turn out to be the same and is given by

\begin{equation} \label{eq:Choice-Eta}
   \eta (x) =\left\{ \begin{array}{cl} 
\frac{\log \big( \frac{|x- p_0|}{r_1} \big )}{ \log ( \frac{r_1^{\prime}}{r_1} )} &
             \hbox{if } \, r_1  < |x- p_0|  < r_1^{\prime}
 \vspace{3mm} \\
1 &\hbox{if } \, r_1^{\prime}   \leq |x-p_0| \leq r_2^{\prime}
 \vspace{3mm} \\
 \frac{\log \big( \frac{|x- p_0|}{ r_2} \big )}{ \log ( \frac{r_2^{\prime}}{r_2} )} &
   \hbox{if }  r_2^{\prime}  < |x- p_0| <  r_2,
\end{array} \right. 
\end{equation}
where $A'=B_{r_2'}(p_0) \setminus \overline{B_{r_1'}(p_0)}$ and $A=B_{r_2}(p_0) \setminus \overline{ B_{r_1}(p_0)}$. For later use, we mention that after explicit computations one gets the following
\begin{equation} \label{eq:Regeln1}
   \frac{1}{n \omega_n}  \int_A  \frac{|\nabla \eta|^2}{|x - p_0|^{n-2}} = 
   \frac{1}{\log (\frac{r^{\prime}_1}{r_1} ) }  + \frac{1}{\log (\frac{r_2}{r^{\prime}_2} ) } =:W,
\end{equation}
\begin{equation} \label{eq:Regeln2}
   \frac{1}{n \omega_n} \int_A |\nabla \eta|^n =   \Big( \frac{1}{\log (\frac{r^{\prime}_1}{r_1} ) } \Big)^{n-1}  + \Big(  \frac{1}{\log (\frac{r_2}{r^{\prime}_2} ) }  \Big)^{n-1}
   \, \leq \,
W^{n-1}.
\end{equation}

\bigskip 

For solutions to \eqref{eq:XiAn},  several types of $W^{1,s}$ estimates with $s \in [1,n)$ can be found in the literature. 
The following is borrowed from  \cite[Lemma 2.3]{DHM} and will be useful for us.
\begin{pro} \label{pro:DolzHung-Mul}
Consider a  bounded $C^1$-domain  $U \subset \mathbb R^n$, $v \in W^{1,n} (U)$, $c \in \mathbb R$ and $\alpha_0 >0$ such that 
\begin{equation} \label{eq:SleepingDogs}
              \int_{ \{ |v - c| \leq \alpha \} }  |\nabla v|^n  \leq c_0 (1 + \alpha),  \quad \forall \alpha \geq \alpha_0.
 \end{equation}             
Then,  for each $s \in [1,n)$  one can find $q=q(s) \geq 1$ and $C=C(s,  c_0, \alpha_0,U)>0$ such that 
\begin{equation} \label{eq:SleepingDogs2}
    \int_{U} |\nabla v|^{s}  \, \leq \,  
   C \Big[ 1 +  \Big(\fint_U |v - c|\Big)^q  \Big] .
\end{equation}
\end{pro} 
\begin{proof}
 Denoting by $\mu_{f} (t) :=  | \{ |f|  > t \} |$ the distribution function of a measurable function~$f$, the claim will follow by playing with $\mu_{|\nabla v|}$
using the technics introduced by Talenti~\cite{Talenti:1979}. This approach has also been exploited in \cite{boccardogallouet1} to derive a priori bounds similar to \eqref{eq:SleepingDogs2},
and for the sake of completeness we give the arguments. 

\medskip

Writing  $\big\{  | \nabla v|  >  \tau \big\} =  \big\{ | \nabla v|  >  \tau, \, |v-c| \leq \alpha \big\} \cup \big\{  |\nabla v| >  \tau , |v-c| > \alpha \big\} $,  
we get 
 \begin{equation} \label{eq:Russell} 
 \mu_{|\nabla v|} (\tau)  \, \leq \, 
     \frac{1}{\tau^n} \, \int_{ \{ |v-c| \leq \alpha \} } |\nabla v|^n \, + \, \mu_{| v-c|} (\alpha) .
 \end{equation}
 For $1 \leq q<n$, by Sobolev-Poincar\'e inequality on $v_\alpha=\min\{|v-c|,\alpha\}$ we deduce
$$
    \mu_{|v-c|} (\alpha)  \, \leq \, \frac{\|v_\alpha\|^q_q}{\alpha^q} \, \leq \, 
    \frac{c_1}{\alpha^q}\Big[\|\nabla v_\alpha\|^q_n + \big( \fint_U v_\alpha \big)^q\Big] \, \leq \,  
\frac{c_1}{\alpha^q}\Big[ \big(\int_{ \{ |v-c| \leq \alpha \} } |\nabla v|^n \big)^{\frac{q}{n}} + \big( \fint_U |v-c| \big)^q\Big] $$
with a constant $c_1:= c_1 (q, U)$.  
Hence, with the hypothesis \eqref{eq:SleepingDogs}, this last inequality combined with  \eqref{eq:Russell} give for all $\tau >0$ and $\alpha \geq \alpha_0$ the following
\begin{equation} \label{eq:Ceramica}
    \mu_{|\nabla v|} (\tau)   \, \leq \,   
   \frac{ c_0 [1 + \alpha]}{\tau^n}    \, + \,  \frac{ c_1}{\alpha^q}   [ c_0  (1+\alpha)]^\frac{q}{n}
\, + \,  \frac{c_1}{\alpha^q} \big( \fint_U |v-c|) \big)^q.
\end{equation}
Fix $0 <\varepsilon < n-s $, $q>\max\{\frac{n}{n-1},\frac{s}{\varepsilon}\}$ and set $t_0 := \alpha_0^{\frac{1}{\varepsilon}}$. Choosing $\alpha := \tau^{\varepsilon}$ in \eqref{eq:Ceramica}
we obtain 
 \begin{eqnarray*}
        \int_U  |\nabla v|^{s} &\leq&   2 t_0^s |U|   \, + \,   \int_{\{ |\nabla v| > t_0 \}  }(|\nabla v|^{s}-t_0^s) 
        \, = \, 
      2 t_0^s  |U|  \, + \,  s  \int_{t_0}^{\infty}  \tau^{s-1}  \mu_{| \nabla v|} (\tau)d \tau \\
       &\stackrel{\eqref{eq:Ceramica}}{\leq}&
        2t_0^s   |U|   \, + \,   c_0 s  \int_{t_0}^{\infty} \big( 1 + \tau^{\varepsilon} \big)  \tau^{s -1 - n}  +  
s c_1 c_0^\frac{q}{n}  \int_{t_0}^{\infty}  (1+\tau^\varepsilon)^\frac{q}{n} \tau^{s-1-\varepsilon q} d\tau\\
&&\, + \,  s c_1(\fint_U |v-c|)^q  \int_{t_0}^{\infty}   \tau^{s-1-\varepsilon q} d\tau  \\
        &\leq&c_2 \Big[1+ \big( \fint_U |v -c|\big)^q\Big] 
        \end{eqnarray*}
with a constant $c_2 := c_2 (s, \varepsilon, c_0, \alpha_0,q,U)$.  Hence, the conclusion \eqref{eq:SleepingDogs2} follows by taking for instance $\varepsilon=\frac{n-s}{2}$ and $q=\frac{n}{n-1}+\frac{2s}{n-s}$.
\end{proof}

\medskip
\subsection{BMO estimates}

Before discussing the space of BMO functions, let us first mention the following preliminary estimates.

\begin{lm}  \label{lm:stimaLuoYang}
For each $p \geq 1$, there exists a constant $C_p >0$ such that 
$$      \big\| u - \fint_E u   \big\|_{L^p (E)}   \leq C_p \, |\Omega|^{\frac{1}{p} }  \, \| f \|_{L^1 (\Omega)}^{\frac{1}{n-1}}   + 2   \,  |E|^{\frac{1}{p}}  \, {\rm osc}_{\partial \Omega}  (u)$$
holds for all  solutions $u \in C(\overline{\Omega})$ of \eqref{eq:XiAn}, and  all  subsets $E \subset \Omega$ of positive measure.
\end{lm} 
\begin{proof} 
Consider the $n-$harmonic function $\Phi$ in $\Omega$ so that $\Phi =u$ on $\partial \Omega$. On the one hand,
by the weak maximum principle we have that 
\begin{equation} \label{eq:YiChi}
      \min_{\partial \Omega} u   \leq \Phi \leq \max_{\partial \Omega} u  
     \quad \hbox{in }\Omega.
\end{equation}
On the other hand, from \cite[Lemma 1.8]{AgPe}, there exists a non-empty interval $(0, \alpha^*)$ such that: 
$$
        \fint_{\Omega} e^{\alpha \frac{|u - \Phi |}{\| f \|_{L^1(\Omega)}^{1/(n-1)}}} \leq  C(\alpha) 
$$
for all $ \alpha \in (0, \alpha^*)$. In particular, for some constant $\tilde C_p >0$ we have 
\begin{equation} \label{1251} 
   \frac{1}{|\Omega|^{1/p}}  \frac{ \| u - \Phi \|_{L^p (\Omega)}}{\| f \|_{L^1(\Omega)}^{\frac{1}{n-1} }}  \leq  \tilde C_p .
\end{equation} 
Using \eqref{eq:YiChi}  we have $\displaystyle \big| \Phi  - \fint_{\Omega}  \Phi  \big| \leq  {\rm osc}_{\partial \Omega}  (u)$  in $\Omega$, and as a consequence,
\begin{equation} \label{eq:LuoYang2}
   \Big|\Phi -\fint_{E} \Phi  \Big |= \Big|\Phi - \fint_{\Omega} \Phi-\fint_{E} [\Phi- \fint_{\Omega} \Phi]\Big|\leq 2 \,  {\rm osc}_{\partial \Omega}  (u) .
\end{equation}
Furthermore,
\begin{equation} \label{eq:LuoYang3}
         \Big|\fint_{E} u  - \fint_{E} \Phi \Big|\leq \fint_{E} |u  - \Phi| \leq   |E|^{-\frac{1}{p}} \| u - \Phi \|_{L^p (E)}.
\end{equation}

\medskip

Hence, using \eqref{eq:LuoYang2} and \eqref{eq:LuoYang3} we deduce
\begin{eqnarray*}
  \big \|  u - \fint_{E}  u  \big\|_{L^p (E)}   &\leq&  \big\|  u - \Phi  \big\|_{L^p (E)}  +  \big\|  \Phi - \fint_{E}  \Phi \big\|_{L^p (E)}  + \big\| \fint_E \Phi  - \fint_{E}  u \big\|_{L^p (E)}   
   \\
   &\leq& 
  2   \| u - \Phi \|_{L^p (\Omega)} + 2  \,  |E|^{\frac{1}{p}}  \,  {\rm osc}_{\partial \Omega}  (u) \\
  &\stackrel{\eqref{1251}}{\leq}&
   2 \tilde C_p  |\Omega|^{\frac{1}{p}}  \| f \|_{L^1 (\Omega)}^{\frac{1}{n-1} }  + 2  \,  |E|^{\frac{1}{p}} \, {\rm osc}_{\partial \Omega}  (u)  ,
\end{eqnarray*} 
and the validity of the lemma follows.
\end{proof}

Let us start by reviewing a few facts about functions of {\it bounded mean oscillation} or for short BMO-functions.
In the seminal paper by John and Nirenberg \cite{JN:1961} those are elements $u \in L^1_{loc} (\mathbb R^n)$ whose 
mean oscillation  $\fint_Q | u - \fint_Q u|$ is uniformly bounded on the set of cubes $Q \subset \mathbb R^n$. 
From the ``John-Nirenberg inequality" one can equivalently consider the $p$-mean oscillation $\fint_Q | u - \fint_Q u|^p$ ($p \geq 1$), and it is known that 
one may use balls instead of cubes. 

\medskip

For the need of the present paper, we say that $u \in L^n (\Omega)$  is a BMO-function on $\Omega \subset \mathbb R^n$ if
$$
    [u]_{BMO,\Omega}^n :=\sup \left\{ \fint_B  \Big|u-\fint_B u \Big|^n:  \ B \subset \Omega , \, B \, {\rm ball }    \right \} < \infty .
$$
We easily check that $[\cdot]_{BMO,\Omega}$ defines a semi-norm on the vector space of BMO-functions, and modulo constants it becomes a norm. 
For later purpose, we note that given a function $u \in L^n (\Omega)$ and balls $B  \subset \tilde B  \subset \Omega$ we have
\begin{equation} \label{eq:Holiday}
   |\fint_{B } u -\fint_{\tilde B} u | \leq \fint_{B} |u  -\fint_{\tilde B} u|
\leq  \frac{| \tilde B | }{|B|}  \fint_{\tilde B} | u  -\fint_{\tilde B} u | \leq  \frac{| \tilde B | }{|B|}   [u]_{BMO,\Omega} .
\end{equation}
 
As it is well known, this estimate can be refined as follows:

\begin{lm} \label{lem:LiTi}
For any ball $B \subset \tilde B  \subset \Omega$ and any function $u \in L^n (\Omega)$, we have
 \begin{equation}      \label{eq:LiTi}
       \big| \fint_{B } u -\fint_{\tilde B} u \big| \leq e^n \Big\{  1 +  \frac{1}{n} \log  \Big(\frac{| \tilde B |}{|B|} \Big)   \Big\} [u]_{BMO, \Omega}
\end{equation}
\end{lm}
\begin{proof}
For the given ball $B_r (p) = B$, consider the largest $k \in \mathbb N$ such that  $B_{ e^k r} (p) \subset \tilde B$. 
A telescopic argument combined with~\eqref{eq:Holiday} and the observation that   $e^{n k} \leq \frac{|\tilde B|}{|B|}$ lead to
\begin{eqnarray*}
    \big| \fint_{B} u -\fint_{\tilde B} u \big| 
    &\leq&  
    \sum_{j=0}^{k-1}  \big| \fint_{B_{e^{j+1} r} (p)  } u -\fint_{B_{e^j r} (p) } u \big| +\big| \fint_{\tilde B} u -\fint_{B_{e^k r} (p) } u \big| \leq e^n (k+1) [u]_{BMO,\Omega} \nonumber \\
    &\leq& 
    e^n \Big\{ 1 +    \frac{1}{n}   \log \Big(  \frac{|\tilde B|}{|B|} \Big)  \Big\} [u]_{BMO,\Omega}
\end{eqnarray*} 
which gives the claim \eqref{eq:LiTi}.
\end{proof}

\bigskip

We now derive  a priori BMO estimates for the solutions to the $n$-Poisson equation~\eqref{eq:XiAn}. It is based on a blowup argument
 that has been provided for quasilinear problem subjected to zero boundary conditions in \cite[Lemma 2.5]{DHM} (see also \cite{EsMa} and \cite{Simon:1996}). 
We show how similar a priori bounds can be obtained without prescribing boundary data.
\begin{thm} \label{thmBMO} 
Let $\Omega \Subset \mathbb R^n$, $ f \in L^1 (\Omega)$
 and $K_\delta=\{x \in \overline{\Omega}: \ \hbox{dist}(x,\partial \Omega)\leq \delta\}$ for some $\delta >0$ fixed. 
Then, there exists a constant $C := C(\delta) >0$ such that 
\begin{equation}  \label{BMOest}
         [u]_{BMO,\Omega} \leq C \big( \| f \|^{\frac{1}{n-1}}_{L^1 (\Omega)} + {\rm osc}_{K_\delta} (u) \big) 
\end{equation}
holds for all solution $u \in W^{1,n}_{loc}(\Omega)\cap C(\overline{\Omega})$  of  \eqref{eq:XiAn}.
\end{thm}
\begin{proof} 
Assuming by contradiction that \eqref{BMOest} is false, we get sequences $\tilde u_k \in W^{1,n}_{loc}(\Omega)\cap C(\overline{\Omega})$ and $\tilde f_k \in L^1 (\Omega)$ satisfying 
\begin{equation} \label{eq:Contradiction2}
    - \Delta_n \tilde u_k = \tilde f_k, 
    \qquad
    \frac{[\tilde u_k]_{BMO,\Omega}}{  \| \tilde f_k \|^{\frac{1}{n-1}}_{L^1 (\Omega)} + {\rm osc}_{K_\delta}  (\tilde u_k)  }  \longrightarrow  \infty .
\end{equation} 
Setting then
$$
  u_k := \frac{\tilde u_k}{ \| \tilde f_k \|^{\frac{1}{n-1}}_{L^1 (\Omega)} + {\rm osc}_{K_\delta}  (\tilde u_k)  } 
  \quad \hbox{ and } \quad
  f_k := \frac{\tilde f_k}{ \big( \| \tilde f_k \|^{\frac{1}{n-1}}_{L^1 (\Omega)} + {\rm osc}_{K_\delta}  (\tilde u_k)  \Big)^{n-1}},
$$
we obtain a sequence $(u_k, f_k)$ that satisfies
\begin{equation} \label{eq:Montesarchio}
   - \Delta_n u_k = f_k , 
     \qquad
 [u_k]_{BMO,\Omega} \to \infty,
\end{equation}
with 
\begin{equation} \label{1022} 
   {\rm osc}_{K_\delta} (u_k) \leq 1
   \quad  \hbox{ and } \qquad \|f_k\|_{L^1(\Omega)} \leq 1  .
    \end{equation}

\medskip

Since $ [u_k]_{BMO,\Omega} \to \infty$, we can find a sequence of balls $B_k :=B_{r_k}(y_k) \subset \Omega$ such that
\begin{equation} \label{contrad}
\fint_{B_k} \Big|u_k-\fint_{B_k} u_k \Big|^n \geq \frac{1}{2} [u_k]_{BMO,\Omega} \xrightarrow{k \to \infty} +\infty
\end{equation}
We claim that 
\begin{equation} \label{0646}
 \lim_{k \to +\infty} r_k=0 ,
\quad  \hbox{dist}(y_k,\partial \Omega)\geq \frac{\delta}{2}.
\end{equation}

Indeed, from \eqref{1022}, \eqref{contrad} and applying Lemma~\ref {lm:stimaLuoYang} we deduce that $|B_k| \to 0$ and then $r_k \to 0$ as $k \to +\infty$.   Furthermore, 
noting that for each ball $B  \subset K_\delta$ we have
$$
    \left(\fint_{B} \Big|u_k-\fint_{B} u_k \Big|^n \right)^{\frac{1}{n}} \leq {\rm osc}_{K_\delta} (u_k) \stackrel{\eqref{1022}}{\leq}  1 ,
$$
assumption \eqref{contrad} with $r_k \to 0$ gives $r_k < \frac{\delta}{2}$ and $\hbox{dist}(y_k,\partial \Omega)\geq \frac{\delta}{2}$. So the claim \eqref{0646} is fully established.

\medskip 

Consider then the rescaled function 
$$
    W_k(z)=\frac{u_k(r_k z+y_k)-\fint_{B_k} u_k}{[u_k]_{BMO,\Omega}},
\quad
   z \in  \Omega_k :=\frac{\Omega - y_k}{r_k},
$$
which, by taking also into account \eqref{contrad},  satisfies
\begin{equation} \label{1052}
 \left\{ 
     \begin{array}{c}
     - \Delta_n W_k (z) = r_k^{n} \frac{ f_k (r_k z+y_k) }{[u_k]^{n-1}_{BMO, \Omega}} ,  \quad   
      [W_k]_{BMO,\Omega_k}=1, 
    \\ \\
     \fint_{B_1(0)} W_k=0,  \quad 
    \fint_{B_1(0)}|W_k|^n \geq \frac{1}{2}.
    \end{array}
    \right.
\end{equation}

\bigskip

Fix $r >1$.  
By applying the Caccioppoli inequality \eqref{Caccio} 
to the solutions $W_k$  of \eqref{1052}, with the balls 
$B_{r} (0) \subset B_{2r} (0) \Subset \Omega_k$ ,
we obtain for each $M>0$ and $c \in \mathbb R$ the estimate
\begin{eqnarray} \label{Caccioresc}
      \int_{\{| W_k - c |\leq M\} \cap B_{r}(0)}|\nabla W_k|^n 
      &\leq& \big( \frac{4n^2}{r} \big)^n\int_{B_{2r }(0)\setminus B_{r}(0)}|W_k - c|^n
      +   \frac{2n M }{[u_k]_{BMO,\Omega}^{n-1}} 
\end{eqnarray}
where we have used
$$
    \Big\| r_k^{n} \frac{ f_k (r_k z+y_k) }{[u_k]^{n-1}_{BMO, \Omega}}\Big\|_{L^1(B_{2r}(0))} \, \leq \,  \frac{ \|f_k\|_{L^1(\Omega)}}{[u_k]^{n-1}_{BMO, \Omega}} \, \stackrel{\eqref{1022}}{\leq} \, \frac{1}{[u_k]^{n-1}_{BMO, \Omega}}.
$$

\medskip

Letting $c_k=\fint_{B_{2r}(0)} W_k$ in \eqref{Caccioresc} and using $ [W_k]_{BMO,\Omega_k} \, = \,1$ we obtain
\begin{equation} \label{eq:Ralf}
    \int_{\{| W_k - c_k |\leq M\} \cap B_{r}(0)}|\nabla W_k|^n  \leq (8n^2)^n \omega_n 
      +   \frac{2n M }{[u_k]_{BMO,\Omega}^{n-1}} .
\end{equation}

The estimate \eqref{eq:Ralf} allows to apply Prop. \ref{pro:DolzHung-Mul} in the ball $B_{r} (0)$, and by taking into account that $ [W_k]_{BMO,\Omega_k}=1$, 
for each $s \in [1,n)$ we  get for some  $q\geq 1$ and $C>0$ (depending on $r$):
\begin{equation} \label{eq:EdisonStation}
       \int_{B_{r} (0)} |\nabla W_k|^s  \, \leq \, C \Big[1 +  \big( \fint_{B_{r}(0)} |W_k-c_k|\big)^q \Big] \, \leq \, C 
\Big[1 +  \Big(1+|c_k-\fint_{B_{r}(0)} W_k|\Big)^q \Big].
\end{equation}
Furthermore, since $\int_{B_1 (0)} W_k = 0$, from \eqref{eq:Holiday} we have for each $r >1$,
\begin{equation} \label{1918}
        \big|\fint_{B_{r}(0)} W_k \big| \leq r^n [ W_k]_{BMO, \Omega_k}  =  r^n  
\end{equation}
and so the  Poincar\'e-Wirtinger inequality with \eqref{eq:EdisonStation} and \eqref{1918} give
   $$\| W_k \|_{W^{1,s} (B_r(0)} \, \leq \, C,
   $$ 
 for some  $C>0$ (depending on $r$). Hence, up to a subsequence, as $k \to +\infty$ we have
\begin{equation}  \label{eq:ZaoShang}
    \left\{
    \begin{array}{c}
      W_k  \rightharpoonup W \, \, \hbox{ in } W^{1,s} (B_{r} (0) ) ,
      \quad \forall s \in [1, n) ,
      \vspace{2mm} \\
      W_k \longrightarrow W \, \, \hbox{ in } L^{n} (B_{r} (0) ) .
   \end{array}
   \right.
\end{equation} 

\medskip By considering the truncation function  
$T^M \! \!: \mathbb R \to \mathbb R$ defined as  
$$
   T^M (s):= \min \{ s, M\}  \, \hbox{ for } s >0, 
   \quad 
   T^M(-s) = - T^M(s) ,
$$
the left-hand side of \eqref{Caccioresc} can be rewritten as $\int_{B_r (0)}   \big|  \nabla [ T^{M} ( W_k - c  )] \big|^n$.
By doing so, the estimate \eqref{Caccioresc} with \eqref{eq:ZaoShang} and the fact that
$ |T^{M} ( W_k - c )| \leq M$ imply
$$
    T^{M} ( W_k - c )  \xrightharpoonup{\quad}  \bar W \quad \hbox{ weakly in } W^{1,n} (B_r(0)) .
$$
Since $T^M$ is a Lipschitz continuous function, by \eqref{eq:ZaoShang} we have that  $T^{M}  ( W_k - c )  \to T^{M} (W -c)$ strongly  in $L^n (B_r(0))$, and therefore
$\bar W = T^{M} ( W -c)$. Hence, 
\begin{equation}  \label{eq:ZaoShang2}
          T^{M} ( W_k - c )  \xrightharpoonup{\quad}  T^{M} ( W -c)  \hbox{  weakly in } W^{1,n} (B_{r}(0) ) .
\end{equation} 

 So, by letting $k \to \infty$ in inequality~\eqref{Caccioresc}, and using \eqref{eq:Montesarchio} with \eqref{eq:ZaoShang} we arrive at
   $$\int_{  B_{r}(0)}   \big|  \nabla [T^{M} ( W - c  )] \big|^n 
    \leq 
     \big( \frac{4n^2}{r} \big)^n \int_{B_{2r }(0)\setminus B_{r}(0)}|W - c|^n.$$
Letting then $M \to \infty$, we conclude 
\begin{equation} \label{eq:Gelato}
      \int_{  B_{r}(0)}   | \nabla  W  |^n 
    \leq 
     \big( \frac{4n^2}{r} \big)^n \int_{B_{2r }(0) \setminus B_{r}(0) } |W - c|^n,
\end{equation} 
and furthermore by \eqref{1052} and \eqref{eq:ZaoShang} we get that
\begin{equation} \label{1253}
  [W]_{BMO,\mathbb{R}^n}\leq 1, \quad \fint_{B_1(0)} W=0, \quad  \fint_{B_1(0)}|W|^n \geq \frac{1}{2}.
\end{equation}
Choosing in  \eqref{eq:Gelato} the value  $c := \fint_{B_{2r}(0)} W$, from \eqref{1253}  we first have
$$    \int_{  B_{r}(0)}   |  \nabla W  |^n  \leq   ( 8n^2)^n  \omega_n   \fint_{B_{2r}(0)} | W - \fint _{B_{2r}(0)} W |^n \leq ( 8n^2)^n  \omega_n ,$$
which implies $ \int_{ \mathbb R^n}   |  \nabla W  |^n < \infty$. Choosing then $c := \fint_{B_{2r}(0) \setminus B_r(0)} W$, from \eqref{eq:Gelato} and 
the Poincar\'e-Wirtinger inequality we have 
\begin{equation} \label{eq:Gelato2}
    \int_{  B_{r}(0)}   | \nabla  W  |^n  \leq \frac{(4 n^2)^n}{r^n \Lambda\Big(B_{2r}(0) \setminus \overline{B_r(0)} \Big)} \int_{B_{2r }(0) \setminus B_{r}(0) } |\nabla W|^n,
\end{equation} 
where for each open set $U \subset \subset \mathbb R^n$ 
$$
     \Lambda (U) :=  \inf \Big\{  \int_{U} |\nabla w|^n \, : \,  \int_U w^n =1, \,  \int_U w = 0 \Big\} >0 .
$$
Noting that $r^n \Lambda\Big(B_{2r}(0) \setminus \overline{B_r(0)} \Big)$ is scaling invariant, and thus independent of $r$, 
by letting $r \to \infty$ in the estimate~\eqref{eq:Gelato2} 
we obtain 
$ \int_{ \mathbb R^n}   | \nabla  W  |^n = 0 $. 
But,  by using \eqref{1253} we end up with the contradiction
\begin{eqnarray*}
         0 = 
         \int_{B_{1}(0)} |\nabla W|^n 
  &\geq&  
  \Lambda\Big(B_1(0)\Big) \int_{ B_{1} (0)}  \big |  W  -  \fint_{B_1(0) }  W    \big|^n
  \geq
  \Lambda \Big(B_1(0)\Big) \frac{\omega_n}{2} .
  \end{eqnarray*}
So \eqref{eq:Contradiction2} cannot hold and this concludes the proof of the theorem. 
\end{proof}

\medskip
\subsection{Oscillation estimates} \label{subsec:Oscillation}

 We now collect some oscillation estimates for the set of solutions to the $n$-Poisson equation \eqref{eq:XiAn}. 
First, we highlight that Serrin's elliptic estimates for quasilinear equations \cite{Ser1} lead to the following ``$L^{\infty}$-version" of Lemma \ref{lm:stimaLuoYang}:

\begin{pro} \label{prop:Musica}
Consider the problem  \eqref{eq:XiAn} with  $ f \in L^{\infty} (\Omega)$. Then, for some constant $C >0$ 
the following a priori estimates hold for any solution $u \in C(\overline{\Omega})$ of \eqref{eq:XiAn}:
\begin{equation} \label{eq:Faust1}
   C \Big\| u - \fint_{B_r (p)} u \Big\|_{L^{\infty} (B_r (p) ) }   \, \leq \, 
   \min \Big\{ [u]_{BMO, \Omega}, \frac{ 1}{r}  \| f \|_{L^1 (\Omega)}^{\frac{1}{n-1}}  + {\rm osc}_{\partial \Omega} (u) \Big\} + r^\frac{n}{n-1}  \| f \|^{\frac{1}{n-1}}_{L^{\infty} (B_{2r} (p)) } 
\end{equation}
on any balls $B_{2 r} (p)  \subset \Omega$; and furthermore
\begin{equation} \label{eq:Faust}
   C  \Big\| u - \fint_{B_R (p_0)} u \Big\|_{L^{\infty} (B_r (p) ) }   \, \leq  \,  \log ( \frac{R}{r}) \,  [u]_{BMO, \Omega} + r^\frac{n}{n-1}  \| f \|^{\frac{1}{n-1}}_{L^{\infty} (B_{2r} (p)) }
\end{equation} 
on any balls $B_{2 r} (p)  \subset  B_{R} (p_0) \subset \Omega$.
\end{pro}
\goodbreak
\begin{proof}
We set $B_r:=B_r(p)$ and $B_{2r}:=B_{2r}(p)$, and we also note that we have by assumption $f \in L^{\frac{n}{n- \varepsilon}} (\Omega)$ for all $\varepsilon \in (0,n)$.  

\medskip

{\it Proof of \eqref{eq:Faust1}:} 
Since $u - \fint_{B_{2r}} u $ solves \eqref{eq:XiAn}, 
the elliptic estimates in \cite[Theorem 2]{Ser1}  give a constant $C := C(\varepsilon)$ such that
\begin{eqnarray}
   \big\| u -  \fint_{B_{r}} u   \big\|_{L^{\infty} ( B_r )} &\leq&   2 \big\| u -  \fint_{B_{2r}} u \big\|_{L^{\infty} ( B_{r} )} \nonumber  \\ 
  &\leq& 
  C  \Big\{  \frac{1}{r} \big\| u - \fint_{B_{2r}} u  \big\|_{L^n (B_{2r}) } +   r^\frac{\varepsilon}{n-1}  \| f \|^{\frac{1}{n-1}}_{L^{\frac{n}{n-\varepsilon}} (B_{2r} ) }  \Big\}  \nonumber  \\  
  &\leq& 
  C \Big\{ \frac{1}{r} \big\| u -\fint_{B_{2r}} u \big\|_{L^n (B_{2r}) }  + r^\frac{n}{n-1}  \| f \|^{\frac{1}{n-1}}_{L^{\infty} (B_{2r} ) }   \Big\} .
    \label{eq:Ospedale}
\end{eqnarray} 
The conclusion follows from
the definition of BMO-norm and by applying
Lemma \ref{lm:stimaLuoYang}  to \eqref{eq:Ospedale}.

\bigskip

{\it Proof of \eqref{eq:Faust}:} 
Given any $x \in B_r := B_r (p)$, by applying \eqref{eq:Faust1} with Lemma~\ref{lem:LiTi}, we obtain
\begin{eqnarray*}
    \Big | u (x)  - \fint_{B_R (p_0 )} u \Big| 
    &\leq&   
    \Big| u (x) - \fint_{B_r } u  \Big|  \, +  \, \Big|  \fint_{B_R (p_0)} u - \fint_{B_r } u   \Big|  \\
    &\leq&  
    C \left( [u]_{BMO, \Omega }   + r^\frac{n}{n-1}  \| f \|^{\frac{1}{n-1}}_{L^{\infty} (B_{2r}) }  \right) + e^n \big( 1 + \log (\frac{R}{r}) \big) [u]_{BMO, \Omega}
    \\
    &\leq&
    C  \big( 1 + \log (\frac{R}{r}) \big) [u]_{BMO, \Omega}   + r^\frac{n}{n-1}  \| f \|^{\frac{1}{n-1}}_{L^{\infty} (B_{2r}(p)) }  
\end{eqnarray*}
and therefore \eqref{eq:Faust} follows (since $2r \leq R$). 
\end{proof} 

\medskip

Obviously, from  \eqref{eq:Faust1} we immediately get a similar priori bound for $\osc_{B_{r} (p)} u$, which will be useful later. 
Furthermore, the claim~\eqref{eq:Faust} leads to the following oscillation estimates on annuli subsets: 

\begin{pro} \label{prop:Birch}
Consider the problem \eqref{eq:XiAn} with $ f \in L^{\infty} (\Omega)$. Given $p_0 \in \Omega$  and $0<\delta <R$, set $A_{R, \delta} := B_{R+\delta} (p_0) \setminus \overline{B_{R- \delta} (p_0)} $.
Then, for some constant $C >0$, we have
\begin{equation} \label{eq:Faust-Liszt1}
   {\rm osc}_{A_{R, \delta} } (u) \leq C \Big\{  \log ( \frac{R}{\delta} ) [u]_{BMO, \Omega} + \delta^\frac{n}{n-1}  \| f \|^{\frac{1}{n-1}}_{L^{\infty} (A_{R, 2 \delta} ) }  \Big\} ,  
 \end{equation}
 whenever $B_{R+2 \delta}(p_0) \subset \Omega$, and for any solution $u \in C(\overline{\Omega})$ to \eqref{eq:XiAn}.

\medskip

In particular,  for all $x \in A_{R,\delta}$ we have
\begin{equation}  \label{eq:Faust-Liszt2}
    \Big | u (x)  - \fint_{\partial B_R (p_0 )} u d \sigma \Big|  \leq 
    C \Big\{  \log ( \frac{R}{\delta} ) [u]_{BMO, \Omega} + \delta^\frac{n}{n-1}  \| f \|^{\frac{1}{n-1}}_{L^{\infty} (A_{R, 2 \delta} ) }  \Big\} .
\end{equation}
\end{pro}
\begin{proof}
Given $x_1, x_2 \in A_{R, \delta}$,  consider balls $B_{2\delta} (p_i) \subset A_{R,2 \delta}$ with $p_i \in \partial B_R(p_0)$ containing $x_i$ ($i=1,2$). 
Applying the estimate~\eqref{eq:Faust} with the two balls  $B_{2\delta} (p_i) \subset B_{R + 2\delta } (p_0) \subset \Omega$ we get
\begin{eqnarray*}
    |u (x_1) - u (x_2)|
    &\leq&   
    \Big| u (x_1) - \fint_{B_{R + 2 \delta} (p_0) } u dx \Big|  \, +  \, \Big| u (x_2) - \fint_{B_{R + 2 \delta} (p_0) } u dx \Big|  \\
    &\leq&  C \left(  \log ( 2 + \frac{R}{\delta} ) \,   [u]_{BMO, \Omega }   + \delta^\frac{n}{n-1}  \| f \|^{\frac{1}{n-1}}_{L^{\infty} (B_{2 \delta} (p_1)) }  + \delta^\frac{n}{n-1}  \| f \|^{\frac{1}{n-1}}_{L^{\infty} (B_{2 \delta} (p_2)) }  \right),
\end{eqnarray*}
and the claim \eqref{eq:Faust-Liszt1} follows. 

\medskip

Noting that $\fint_{\partial B_R (p_0 )}  |u(x) - u(y)| d \sigma (y) \leq {\rm osc}_{A_{R, \delta} } (u)  $ for all $x \in A_{R, \delta}$, 
the estimate \eqref{eq:Faust-Liszt2} follows then from \eqref{eq:Faust-Liszt1}.
\end{proof} 

\bigskip
\section{Nonlinear spherical average for $n$-Poisson problems} 
\label{sec:Spherical-Mean}
\medskip

For each $u \in C (B_{R}(p_0) )$ defined on a ball  $B_R(p_0) \subset \mathbb R^n$ one can consider its spherical average $\bar u(r) := \fint_{\partial B_r (p_0)} u$ for $0<r<R$, and this generates an operator that commutes with the Laplacian.
In particular, if one considers a solution $u$ to the linear Poisson equation $-\Delta u = f$ in $B_R (p_0)$, then $\bar u (r)$ gives a radial function that satisfies the ODE 
$$
             -\frac{ ( r^{n-1} \bar u')'}{r^{n-1}}= \fint_{\partial B_r(p_0)}  f,
             \quad r \in (0, R).
 $$
This is extremely useful since information for $\bar u$ can be derived by an ODE analysis and estimates on $u$ are reduced to have a good control on $u - \bar u$.

 \medskip
 
 If one replaces the Laplacian by the operator $\Delta_n$ with $n \geq 3$ in $\mathbb R^n$, its quasilinearity prevents to use the function $\bar u$ in the way previously described.
 Given $(u,f)$ with 
 \begin{equation} \label{eq:Poisson-Ball}
      -\Delta_n u = f \hbox{ in }   B_{R}(p_0)  \subset \mathbb R^n ,
 \end{equation}
we  consider the  ``nonlinear spherical average" of $u$, defined as the unique radial solution $\tilde u$  to: 
\begin{equation} \label{eq:1209} 
\left\{ \begin{array}{rll}
    -\Delta_n \tilde u &= \fint_{\partial B_{r}(p_0)} f  & \hbox{in }B_R(p_0),
    \vspace{2mm} \\
   \tilde u   &=\fint_{\partial B_R(p_0)} u  &\hbox{on }\partial B_R(p_0).
\end{array}\right. 
\end{equation}
In the case of the Laplacian we have $\tilde u (r) = \bar u (r) $ for all $r \in (0, R)$, but in the quasilinear case we need to make a distinction between those two radial functions, and 
we are led to consider separately $u - \tilde u$ and $u - \bar u$. Those two functions coincide on $\partial B_R (p_0)$, but contrary to the linear case it is not obvious how to derive a priori bound on those two functions from their boundary value. 
A priori bounds on $u - \bar u$ have been obtained in the Section \ref{subsec:Oscillation} 
and here we focus on the  function $u - \tilde u$.  
 
\medskip

Observe first that since in radial coordinates the function $\tilde u (r)$ (with $r := |x - p_0|$)  solves the ODE problem
$$
\left\{
\begin{array}{c}
   - \frac{ \big( r^{n-1} |\tilde u'|^{n-2}\tilde u' \big)'}{r^{n-1}} = \fint_{\partial B_{r}(p_0)} f  ,
   \vspace{2mm} \\
   \tilde u' (0) = 0 ,
   \qquad
   \tilde u (R) := \fint_{\partial B_{R}(p_0)} u ,
\end{array} \right.
$$
its derivative $\tilde u' (r)$ can be written explicitly as
 $$
     - |\tilde u' |^{n-2} \tilde u' (r) =  \frac{1}{n\omega_n r^{n-1} } \int_0^r  \int_{\partial B_{t}(p_0)} f d \sigma  dt =  \frac{1}{n\omega_n r^{n-1} }   \int_{B_r (p_0)} f .
 $$
Thus, we have 
 \begin{equation}  \label{eq:00943}
       |\tilde u' (r) | = \frac{1}{r} \Big|\frac{1}{n\omega_n} \int_{B_{r}(p_0)} f \Big|^{\frac{1}{n-1}} .
\end{equation} 

\medskip

When $f \geq 0$, the relation  \eqref{eq:00943} can be rephrased as 
\begin{equation}  \label{eq:Cornell-Hospital1} 
     \big( | \nabla \tilde u|  |x - p_0|  \big)^{n-1} 
     \, = \,      \big( \frac{1}{n\omega_n}  \big)    \| f \|_{L^1 (B_{|x - p_0|} (p_0) ) } 
\end{equation} 
for any $x \in B_R (p_0)$, and it also gives the explicit representation
\begin{eqnarray}
  \tilde u (r)=
   \fint_{\partial B_{R}(p_0)} u
   \, + \, \int_r^{R} (\frac{1}{n \omega_n} \int_{B_t(p_0)} f \big)^{\frac{1}{n-1}} \frac{dt}{t} \label{eq:Eq-Spherical-Mean} .
\end{eqnarray}

 \bigskip
  
 A first relevant result is given by he following estimate.

\begin{thm} \label{thm:Spherical-Mean}
Let $u$ be a solution to \eqref{eq:Poisson-Ball}, $\tilde u$  be the corresponding solution to \eqref{eq:1209} and 
consider two annuli $A^{\prime}  \subset A \subset B_R (p_0)$ centered at $p_0$: 
$$
    A^{\prime}  := \{ r_1^{\prime}  < |x - p_0|  < r_2^{\prime}   \} ,
   \qquad
   A  := \{ r_1  < |x - p_0| < r_2   \} .
$$
Then, there holds
\begin{eqnarray} 
       && 
       \int_{A^{\prime}}  \Big\{  | \nabla  ( u - \tilde u) |^{n} +    |\nabla \tilde u|^{n-2} |\nabla (u - \tilde u )|^2  \Big\} \nonumber \\ 
       && 
       \quad \leq  \,
            C  \| u - \bar u \|_{L^{\infty} (A)}   
            \Big\{
            \| f  \|_{L^1(A)}  
       \, + \, 
        \| f \|^{\frac{n-2}{n-1} }_{L^1 (B_R(p_0))}  \|u - \bar u \|_{L^{\infty} (A)} W 
        \, + \, 
         \|u - \bar u \|_{L^{\infty} (A)}^{n-1} W^{n-1} 
        \Big\} 
        \, \qquad
        \label{eq:1647}
\end{eqnarray}
where $ \displaystyle W :=  \frac{1}{\log \big( \frac{r^{\prime}_1}{r_1} \big) }  + \frac{1}{\log \big(\frac{r_2}{r^{\prime}_2} \big) } $.
\end{thm}
\begin{proof}
By considering in \eqref{eq:Comparison-Prob} the comparison problem $-\Delta_n \tilde u =  \fint_{\partial B_{r}(p_0)} f$, the result will be 
a consequence of \eqref{eq:Ich-Stoere-Sie} applied with the radial function ${\mathcal R} :=\bar u- \tilde u $ and a radial cut-off function $\eta$ defined as in \eqref{eq:Choice-Eta}.

\medskip

More specifically, observe first that the function $g (x)  := f (x) - \fint_{\partial B_{r}(p_0)} f$ with $r := |x - p_0|$ satisfies
\begin{equation} \label{eq:312}
        \int_{\partial B_r (p_0) } g  \, d \sigma =0 \qquad \forall \ 0 <  r \leq R, 
\end{equation}
and consequently on any ball $B_{r} (p_0) \subset B_R (p_0)$ and any radial function $\Psi (|x - p_0|)$ we have
   $$     \int_{B_{r}  (p_0)} g (x)  \Psi (| x - p_0| ) \, dx = \,  \int_0^{r} ds \, \Psi (s)  \int_{\partial B_{s} (p_0)} g    d\sigma   \, \stackrel{\eqref{eq:312}}{=} \, 0 .$$
In particular, we have $ \int_{B_{r}  (p_0)} g (x)  \eta^n (\bar u-\tilde u )  \, dx = 0$, and so the integral $I_0$ in \eqref{eq:Ich-Stoere-Sie} is zero.

Hence, applying Prop. \ref{prop:Benutzen2} and using then \eqref{eq:Cornell-Hospital1}  we get for some constant $C_1>0$ the following 
\begin{eqnarray*}
 && C_1 \int_A \eta^n  \Big\{ |\nabla (u - \tilde u )|^{n-2}+|\nabla \tilde u|^{n-2} \Big\} |\nabla (u - \tilde u )|^2  
       \nonumber \\
      && \, \leq \,   
    2 \|  f  \|_{L^1(A)} \,  \| u -\bar u \|_{L^{\infty} (A)}   \, + \,  
       \| u - \bar u \|^n_{L^{\infty} (A)} \int_A |\nabla \eta|^n +     \| u - \bar u \|^2_{L^{\infty} (A)}  \int_A  |\nabla \tilde u|^{n-2} |\nabla \eta|^2
      \nonumber \\
  &&
    \stackrel{\eqref{eq:Cornell-Hospital1}}{\leq} 
   2  \| u -\bar u \|_{L^{\infty} (A)}   \| f  \|_{L^1(A)} 
    \nonumber  \\
  && 
  \qquad \, + \, 
      \| u - \bar u \|^{n}_{L^{\infty} (A)} \int_A |\nabla \eta|^n +   
      (n \omega_n)^{-\frac{n-2}{n-1}} \| u - \bar u \|^{2}_{L^{\infty} (A)}   \| f \|_{L^1 (B_R(p_0)) }^{\frac{n-2}{n-1}} \int_A  \frac{|\nabla \eta|^2}{|x - p_0|^{n-2}}  .
 \qquad 
\end{eqnarray*}

\medskip

Using the test function $\eta$ defined by \eqref{eq:Choice-Eta} with the related estimates \eqref{eq:Regeln1}-\eqref{eq:Regeln2} we then reach the conclusion~\eqref{eq:1647}.
\end{proof}

\bigskip

When $n=2$, the differential operator in \eqref{eq:Poisson-Ball}  reduces to the Laplacian and $\bar u \equiv \tilde u $ since they both satisfy equation \eqref{eq:1209}, as previously observed.  
In the quasilinear case $n \geq 3$, the oscillation of $\bar u - \tilde u$ does not vanish and can be estimated using the following two propositions:

\begin{pro} \label{prop:Capodanno}
Assume $n \geq 3$. Let $u$ be any solution to \eqref{eq:Poisson-Ball} with $f \geq 0$, $\tilde u$  be the corresponding solution to \eqref{eq:1209} and 
consider an annulus $A^{\prime} := \{ r_1^{\prime} < |x - p_0 | < r_2^{\prime}\} \subset B_R (p_0)$. Then,  
for some constant $C >0$ (independent of $u$), there holds
\begin{eqnarray*}
 [ {\rm osc}_{A^{\prime}}  \, ( \bar u - \tilde u )]^2     \leq C
\| f \|_{L^1 (B_{r^{\prime}_1}(p_0))  } ^{-\frac{n-2}{n-1}}
  \log \big( \frac{r^{\prime}_2}{r^{\prime}_1}\big) \,   \int_{A^{\prime}} |\nabla \tilde u|^{n-2} |\nabla (u - \tilde u) |^2
\end{eqnarray*}
\end{pro}
\begin{proof}
Given $s , t \in [r_1^{\prime},  r_2^{\prime}]$, we have 
\begin{eqnarray} \label{eq:0919}
            \big| (\bar u - \tilde u )(t)-(\bar u -\tilde u)(s)       \big|
          &=& \Big| \int_{s}^{t} \frac{d}{dr}  (\bar u -\tilde u)dr \Big| 
= 
    \Big| \int_{s}^{t} dr \fint_{\partial B_r(p_0)}  \langle \nabla (u - \tilde u) , \nu \rangle \Big| \nonumber \\
&\leq & \frac{1}{n \omega_n}  \int_{A^{\prime}} \frac{|\nabla (u-\tilde u)|}{|x-p_0|^{n-1}}. 
\end{eqnarray}
By the Cauchy-Schwarz inequality we have the estimate 
\begin{eqnarray} \label{eq:0919bis}
\Big(\int_{A^{\prime}} \frac{|\nabla (u-\tilde u)|}{|x-p_0|^{n-1}}\Big)^2&\leq&
\Big( \int_{A^{\prime}}  \frac{1}{|\nabla \tilde u|^{n-2} |x - p_0|^{2(n-1)} } \Big)
 \Big( \int_{A^{\prime}} |\nabla \tilde u|^{n-2} |\nabla (u - \tilde u) |^2 \Big) \nonumber \\
&\leq & (n\omega_n)^{\frac{2n-3}{n-1}}  \| f \|_{L^1 (B_{r^{\prime}_1}(p_0))}^{-\frac{n-2}{n-1}}
  \log \big( \frac{r^{\prime}_2}{r^{\prime}_1} \big) \int_{A^{\prime}} |\nabla \tilde u|^{n-2} |\nabla (u - \tilde u) |^2 \end{eqnarray}
in view of \eqref{eq:Cornell-Hospital1}. Inserting \eqref{eq:0919bis} into \eqref{eq:0919} we conclude the proof.
\end{proof}

\bigskip

Above proposition immediately gives 
 $$\big[  {\rm osc}_{A^{\prime}} (u - \tilde u)  \big]^2 
  \leq 
  C \Big(
   \big[ {\rm osc}_{A^{\prime}} (u - \bar u)  \big]^2 
\, + \,   
   \| f \|^{-\frac{n-2}{n-1}}_{L^1 (B_{r^{\prime}_1}(p_0))  } 
  \log ( \frac{r^{\prime}_2}{r^{\prime}_1}) \int_{A^{\prime}}  |\nabla \tilde u|^{n-2} |\nabla (u - \tilde u) |^2 
  \Big)$$ 
and in our application we will combine it with the estimate~\eqref{eq:1647}.

\bigskip

\begin{pro}  \label{prop:Storia-General}
Assume $n \geq 3$. As in Prop.~\ref{prop:Capodanno} consider any solution $u$ to \eqref{eq:Poisson-Ball} with $f \geq 0$, and also two annuli 
$$
   A' :=  \Big\{r'_1 =  r_1  \big( \frac{r_2}{r_1} \big)^{\theta_0} <  |x - p_0| <  r_2 \big( \frac{r_1}{r_2} \big)^{\theta_0}  = r'_2 \Big\}   \subset  A = \{ r_1 < |x - p_0| < r_2 \} \subset B_{1} (p_0)
$$
defined for some fixed $\theta_0 \in (0,\frac{1}{2} )$. Then, for some constant $C = C(\theta)$ we have
 \begin{eqnarray*}
 && [ {\rm osc}_{A^{\prime}} (\bar u - \tilde u) ]^2 \\
 &\leq& 
 C  \frac{\| u - \bar u\|_{L^{\infty} (A)}}{ \| f \|_{L^1( B_{r^{\prime}_1}(p_0) )}^{\frac{n-2}{n-1}}}
   \Big\{ \log \big( \frac{r_2}{r_1} \big) \| f \|_{L^1(A)}    + \| u - \bar u \|_{L^{\infty} (A)}  \| f \|_{L^1( B_{1}(p_0) )}^{\frac{n-2}{n-1}}
           +  \,  \frac{ \| u - \bar u \|^{n-1}_{L^{\infty} (A) } }{ \big(  \log \frac{r_2}{r_1} \big)^{n-2}  } \Big\} .
\end{eqnarray*}
\end{pro}
\begin{proof}
The result follows from Proposition \ref{prop:Capodanno} together with Theorem \ref{thm:Spherical-Mean} by also noting that in the inequality \eqref{eq:1647} we have
$$
   \log \big( \frac{r'_2}{r'_1}  \big) \, W =  \frac{2}{\theta} (1 - 2 \theta) ,
   \qquad
      \log \big( \frac{r'_2}{r'_1}  \big) \, W^{n-1}  = (1 - 2 \theta)    \big(\frac{2}{\theta} \big)^{n-1}  \Big( \frac{1}{\log \frac{r_2}{r_1} } \Big)^{n-2}   .
$$
\end{proof}

The previous estimates on $u - \tilde u$ and $\tilde u - \bar u$ will have strong implications for the $n$-Liouville equation. However,  the interest for such result goes largely beyond 
the problem we are considering, because it identifies for the first time (at least to our knowledge) the correct ``spherical mean'' substitute 
to be used for PDEs involving $\Delta_n$ in $n-$dimensions. It would be interesting to understand if similar techniques work also for the quasilinear operator $\Delta_p$ when $p<n$.

\goodbreak
\section{Decay with inverse-power source terms} 

In order to exploit in Section~\ref{sec:AsProfile} some decay bounds obtained in \cite{EsLu} for blowing-up sequences of solutions to the $n$-Liouville problem, we provide here some results about the logarithmic decay of radial solutions to the 
$n$-Poisson problem $-\Delta_n u = \frac{1}{|y|^{n+\alpha}}$.

\begin{pro} \label{lem:Safoud}
Given $D \in \mathbb R$, $\alpha >0$ then the radial solutions to 
\begin{equation} \label{eq:Saima}
\left\{
\begin{array}{c}
    -\Delta_n \Phi \, = \, \frac{1}{|y|^{n+\alpha}}  \quad 
    \hbox{ in } \mathbb{R}^n \setminus \overline{B_1(0) }
   \vspace{2mm} \\
    \Phi \big|_{\partial B_{1}} = D
\end{array}
\right.
\end{equation}
are given by the one-parameter family of functions
\begin{equation} \label{eq:Solution-Phi}
   \Phi_A  = D + \int_1^r \Big| A + \frac{1}{\alpha s^{\alpha}} \Big|^{\frac{2-n}{n-1}}  \Big( A + \frac{1}{\alpha s^{\alpha}}\ \Big)  \frac{ds}{s} , \quad A \in \mathbb R.
\end{equation}
Furthermore,  for some $C_{\alpha} > 0$ we have the following asymptotic decay
\begin{equation}  \label{eq:Expansion-Phi}
   D  \, \leq \, \Phi_A  - {\rm sgn} (A) |A|^{\frac{1}{n-1}} \log |y| \, \leq \,  D  +C_{\alpha}  
   \quad
   \hbox{ on } \mathbb R^n \setminus B_{1} (0) .
\end{equation}
\end{pro}
\begin{proof}
Expressing the problem \eqref{eq:Saima} in radial coordinates, the derivative of any radial solution $\Phi$ is given by
\begin{equation} \label{eq:Safoud}
  \Phi'  (r) = \frac{1}{r} \Big| A + \frac{1}{\alpha r^{\alpha}} \Big|^{\frac{1}{n-1} -1}  \Big( A + \frac{1}{\alpha r^{\alpha}}\ \Big) 
\end{equation}
where $A :=  |\Phi'(1)|^{n-2} \Phi'(1) - \frac{1}{\alpha}$ is an arbitrary constant. By integrating \eqref{eq:Safoud} we arrive to the expression~\eqref{eq:Solution-Phi}.

\medskip

The explicit solution $\Phi_A$ can be rewritten as follows
\begin{eqnarray*}  
       \Phi_A (r) &=&    D +{\rm sgn} (A) |A|^{\frac{1}{n-1}} \log r  \nonumber \\
       &&   +  \int_1^r  \underbrace{\Big\{ - {\rm sgn} (A)  |A|^{\frac{1}{n-1}}  +  \Big| A  + \frac{1}{ \alpha t^\alpha} \Big|^{\frac{2-n}{n-1}} 
       \Big( A + \frac{1}{\alpha t^\alpha} \Big) \Big\}}_{I(t)} \frac{dt}{t} .
\end{eqnarray*}

The integrand $I$ can be written as follows:
$$
I (t) = 
\left\{ 
  \begin{array}{cl}
   -  A ^{\frac{1}{n-1}}  +  \Big(  A   + \frac{1}{\alpha t^\alpha} \Big)^{\frac{1}{n-1} } 
   &\hbox{ if } A \geq 0 \vspace{1mm}  \\
    |A|^{\frac{1}{n-1}}  -  \Big| |A|  - \frac{1}{\alpha t^\alpha} \Big|^{\frac{2-n}{n-1}} \Big(  |A|  - \frac{1}{ \alpha t^\alpha}  \Big)
    &\hbox{ if } A < 0 .
  \end{array}
\right.
$$

Noting that $I \geq 0$ we infer 
$$
    \Phi_A  \geq D + {\rm sgn} (A)   |A|^{\frac{1}{n-1}}  \log r
    $$
which gives the left-hand side inequality in \eqref{eq:Expansion-Phi}. 

\medskip

To get the upper estimate in \eqref{eq:Expansion-Phi}, we first observe the following inequalities for some $C_0 \geq 1$:
\begin{equation} \label{eq:FirstAlgebra}
   1 - x^{\frac{1}{n-1}}  \leq (1 - x)^{\frac{1}{n-1}}, 
   \quad \forall x \in [0,1],
\end{equation}
\begin{equation} \label{eq:SecondAlgebra}
  1+x^{\frac{1}{n-1}} \leq C_0 (1 + x)^{\frac{1}{n-1}} \quad \forall x \geq 0 .
\end{equation}

\medskip

When $A \geq 0$, the inequality \eqref{eq:FirstAlgebra} applied with $x = \Big(  1   +\frac{1}{ \alpha A t^\alpha} \Big)^{-1}$ gives

\begin{equation} \label{eq:SeminarioMartedi2}
    I (t) =   -  A ^{\frac{1}{n-1}}  +  \Big(  A   + \frac{1}{ \alpha t^\alpha} \Big)^{\frac{1}{n-1} } 
     \leq  (\alpha t^\alpha)^{-\frac{1}{n-1}}.
\end{equation}

\medskip

When $A \leq 0$, applying \eqref{eq:FirstAlgebra} with $x =1   -\frac{1}{ \alpha |A| t^\alpha}$ if $|A| -  \frac{1}{\alpha t^\alpha}  \geq 0$ and \eqref{eq:SecondAlgebra} if the other sign holds, we derive

\begin{equation} \label{eq:SeminarioMartedi}
      I (t) = |A|^{\frac{1}{n-1}}  -  \Big| |A|  - \frac{1}{\alpha t^\alpha} \Big|^{\frac{2-n}{n-1}} 
      \Big( |A| - \frac{1}{\alpha t^\alpha} \Big) \, \leq \, C_0  (\alpha t^\alpha)^{-\frac{1}{n-1}} ,
\end{equation}
independently of the sign of $ |A| - \frac{1}{\alpha t^\alpha}$.  

\medskip

Therefore, for any $A \in \mathbb R$, from \eqref{eq:SeminarioMartedi2} and \eqref{eq:SeminarioMartedi} we obtain
\begin{eqnarray*}
  \Phi_A - {\rm sgn} (A) |A|^{\frac{1}{n-1}} \log r  = 
   D +  \int_1^r  I(t) \frac{dt}{t}   \leq  D +
    \frac{C_0}{\alpha^{\frac{1}{n-1}} }
    \int_{1}^{\infty} \frac{1}{t^{ \frac{\alpha}{n-1} +1}} dt  
=D + C_{\alpha}. 
\end{eqnarray*}
This concludes the uniform expansion stated in \eqref{eq:Expansion-Phi}.
\end{proof} 

\bigskip

The above result will provide a comparison function that will be used to get an  upper bound for subsolutions to the $n$-Poisson equation in \eqref{eq:Saima}. 
\begin{pro} \label{prop:Galois}
Consider $\varphi$ satisfying
   $$        0 \leq -\Delta_n   \varphi  \, \leq \, \frac{1}{|y|^{n+\alpha}} \quad \hbox{ in } B_{R}(0) \setminus \overline{B_1(0) }.$$
Setting $\displaystyle m_r  :=  \min_{\partial B_r (0)} \varphi $ with $r \in [1,R]$, and defining   
$$
\displaystyle \alpha_R  :=  \frac{m_1 - m_R}{\log (R)}  , 
\quad
D = \max \big\{   \mathop{\rm{osc}}_{\partial B_1 (0) } \varphi,   \mathop{\rm{osc}}_{\partial B_R (0) } \varphi \big\}, 
$$ 
then, for some $C_{\alpha} > 0$  the following estimate holds in $\overline{ B_{R}(0) } \setminus B_1(0)$.
\begin{equation} \label{eq:LogEstimate2}
   -  \alpha_R \log |y|   \,  \leq \, \varphi  - m_1   \, \leq \,  D  + C_{\alpha}  -  \alpha_R \log |y|   .
\end{equation} 
\end{pro} 
\begin{proof}
Let us set $\varphi_0 := \varphi - m_1$. We start by proving the left-hand side inequality in \eqref{eq:LogEstimate2}. 
On $\partial B_1 (0)$ we have 
\begin{equation} \label{eq:SuSfera1}
      - \alpha_R \log |y|  \, =  \, 0  \, \leq \,   \varphi_0 (y) ,
     \quad \forall y \in \partial B_1 (0),
\end{equation} 
whereas on $\partial B_R(0)$,  we notice  that 
\begin{equation} \label{eq:SuSferaR}
    \quad
    - \alpha_R\log |y|  \, =  \,   m_R - m_1  \, \leq \, \varphi_0 (y) ,
     \quad \forall y \in \partial B_R (0).
\end{equation} 
Furthermore, since on the domain $B_R(0) \setminus \overline{B_1 (0) }$ we have 
$$  
         - \Delta_n  \big( -  \alpha_R  \log (|y|)   \big) \, = 0 \, \leq \, -\Delta_n \varphi_0 ,
         $$
from \eqref{eq:SuSfera1},  \eqref{eq:SuSferaR} and the weak maximum principle we deduce 
\begin{equation} \label{eq:Nevicata3}
  \qquad \quad
      - \alpha_R  \log |y|    \, \leq \,  \varphi_0 = \varphi - m_1
    \quad \hbox{ in } \overline{B_{R}(0) }\setminus B_1(0) .
\end{equation}

\bigskip

Let us now prove the right-hand side inequality in \eqref{eq:LogEstimate2}.  
Referring to Proposition~\ref{lem:Safoud}, the function $\Phi$ defined in \eqref{eq:Solution-Phi} with 
$$
      {\rm sgn A} |A|^{\frac{1}{n-1}} \, = \,  - \alpha_R ,
      \qquad
      D = \max \big\{  {\rm osc}_{\partial B_1 (0)} \varphi,  {\rm osc}_{\partial B_R (0)} \varphi \big\},
 $$
 satisfies for some constant $C_{\alpha} >0$ the following
\begin{equation} \label{eq:OspedaleCornell}
\left\{
\begin{array}{c}
    -\Delta_n \Phi =\frac{1}{|y|^{n+\alpha}}  \quad 
    \hbox{ in } \mathbb{R}^n \setminus B_1(0)  ,
   \vspace{2mm} \\
    \Phi  = D
    \quad \hbox{ on } \partial B_{1} ,
     \vspace{2mm} \\
   D  - \alpha_R   \log |y| \leq \Phi  \leq D + C_{\alpha} -  \alpha_R  \log |y| ,
   \quad
   |y|  \geq 1.
\end{array}
\right.
\end{equation}

\medskip

On $\partial B_1 (0) $ we have 
\begin{equation} \label{eq:MusicaSchubert}
   \qquad
    \varphi_0 (y) + \alpha_R   \log |y|  \, \leq \, 
   \max_{\partial B_1 (0)} \varphi_0 \, = \, {\rm osc}_{\partial B_1 (0)} \varphi \leq D,
   \quad
   \forall y \in \partial B_1 (0) ,
\end{equation} 
while for  $y \in \partial B_R (0) $ we have 
\begin{equation} \label{eq:MusicaSchubert2}
  \varphi_0 (y) + \alpha_R   \log |y|  \, = \,  
  \varphi (y) - \min_{\partial B_R (0) } \varphi \, \leq \,  {\rm osc}_{\partial B_R (0)}  \varphi \, \leq \, D .
 \end{equation}

\medskip

Thus, the inequalities \eqref{eq:MusicaSchubert}-\eqref{eq:MusicaSchubert2} with the lower bound on $\Phi$ stated in \eqref{eq:OspedaleCornell} give
$$ 
    \qquad
     \varphi_0 (y)  \, \leq \, D -  \alpha_R   \log |y|  \, \leq \, \Phi (y) 
     \quad
     \forall y \in \partial B_1 (0) \cup \partial B_R (0) ,
$$     
and since by assumption  
$ - \Delta_n \varphi_0 \leq  - \Delta_n \Phi $, the weak comparison principle implies
\begin{equation} \label{eq:TestiScritti}
          \qquad  \varphi_0   \, \leq  \,   \Phi  \quad \hbox{ in }   \overline{ B_{R}(0) }\setminus B_1(0) .
\end{equation}

\medskip

Therefore,  from \eqref{eq:TestiScritti} and the upper bound on $\Phi$ stated in \eqref{eq:OspedaleCornell} we obtain 
\begin{equation} \label{eq:Nevicata4}
\qquad \quad
           \varphi_0 (y) \, \leq \,  \Phi (y)  \, \leq \,  D + C_{\alpha} - \alpha_R  \log |y| , 
           \quad
           \forall y \in  \overline{ B_{R}(0) } \setminus B_1(0).
\end{equation} 
From \eqref{eq:Nevicata3} and \eqref{eq:Nevicata4}, we deduce that $\displaystyle \varphi_0 (y) :=  \varphi (y)  - m_1$ satisfies the conclusion  \eqref{eq:LogEstimate2}.
\end{proof}

\bigskip \goodbreak
\section{Oscillation estimates for an isolated blow-up sequence}
\label{sec:Oscillation-Liouville}
\medskip

The results of the  previous sections will be applied to get uniform estimates for a sequence of solutions to \eqref{Liouville-Quasilinear} that blows-up at $0 \in \Omega$. 

\subsection{Some general results} 
\label{sec:BlowupMechanism}

We first point out the following result which is the analogue of \cite[Lemma 2.2]{YYLi-Harnack} in  the case of the $n$-Laplacian.
\begin{lm} \label{lem:Arrival}
Consider a sequence $(h_k, u_k)$ such that \eqref{eq:TaschenTuch-Mouchoir2} and \eqref{1553} holds, and let $K$ be a compact subset of $B_2 (0) \setminus \{ 0 \} $. 
Then,  $\sup_k \osc_{K} (u_k)<+\infty$.
\end{lm}
\begin{proof}
With our assumptions we have $\| h_k e^{u_k}\|_{L^{\infty} (K)} = O(1)$, and note that without loss of generality we can assume $K$ to be connected.
We then cover the compact set $K$ 
with a chain of finitely many balls $ \{ B_{r_j}(p_j) \}_{ 1 \leq j \leq J}$ satisfying  $B_{r_j}(p_j) \cap B_{r_{j+1}}(p_{j+1}) \not= \emptyset$ for $1 \leq j \leq J-1$ 
and $ \overline{B_{2r_j}(p_j) } \subset B_2(0) \setminus \{ 0 \}$ for $ 1 \leq j \leq  J$.  
An application of the estimate \eqref{eq:Faust1} with $f:=  h_k e^{u_k}$ on each $B_{r_j}(p_j)$ gives a uniform bound $\osc_{B_{r_j}(p_j)} u \leq C_j$, and 
the conclusions follows by choosing an arbitrary point in each $B_{r_j}(p_j) \cap B_{r_{j+1}}(p_{j+1})$.
\end{proof} 

\bigskip

In \cite{EsLu}, we derived some decay estimates for a blow-up solutions sequence to \eqref{758} which rely on some bounds that can be phrased as follows.

\begin{pro} \label{thm:Giovedi}
Given $0 < a \leq b$, consider a family $\Lambda \subset \{ H \in L^{\infty} (B_2 (0) ) \, : \,   a \leq H \leq b \}$ which is equicontinuous at each point, and denote by $\mathcal U$ the set of solutions to 
\begin{equation} \label{eq:Giovedi}
   \qquad 
   - \Delta_n V = H e^V 
   \quad
   \hbox{ in } B_2 (0) ,
   \quad
   \hbox{ with } H \in \Lambda.
\end{equation}
Then, setting  $A_0 := B_{2} (0) \setminus \overline{B_{\frac{1}{2}} (0)}$, there exist some constants $c_1, c_2, c_3  > 0$ such that
\begin{equation} \label{eq:AmaliaBis}
   \qquad
   \sup_{\partial B_1 (0)} V  \, \leq \,  c_1 -   c_2  V(0)  + c_3 \max_{A_0} \{ e^\frac{V}{n-1} \},
   \quad \forall V \in {\mathcal U}.
\end{equation}
\end{pro} 
\begin{proof}
Given $V \in {\mathcal U}$ and setting $M = \sup_{A_0} V$, consider  in the annulus $A_0$ the non-negative function $M-V$ which satisfies $\Delta_n ( M -V) =  H e^V$.
By applying the Harnack inequality (see \cite[Theorem 7]{Ser1}) we get for some $C_0 > 1$
$$
   \sup_{\partial B_1 (0)} (M- V )   \, \leq \,  C_0 \big\{ \inf_{\partial B_1 (0)} (M- V)  \, + \, \| H e^V \|_{L^{\infty} (A_0)}^{\frac{1}{n-1}}  \big\}   ,
   \quad 
   \forall V \in {\mathcal U}.
$$
By rearranging and noting that  $\displaystyle \inf_{\partial B_1 (0)} V = \inf_{B_1 (0)} V$ (by the minimum principle) we get for some constant $\tilde C_0 >0$
\begin{eqnarray}
   C_0 \sup_{ \partial B_1 (0)} V    
   &\leq&
   \inf_{ \partial B_1 (0)} V   \, + \, ( C_0 -1) \sup_{A_0} V  \, + \, C_0 \| H e^V \|_{L^{\infty} (A_0)}^{\frac{1}{n-1}}   
   \nonumber \\
     &\leq& 
   \inf_{B_1 (0)} V  + \, \tilde C_0  \| e^V \|_{L^{\infty} (A_0)}^{\frac{1}{n-1}}  ,
    \quad  \forall V \in {\mathcal U} .
   \label{eq:Tatsachlich-Princeton2}
\end{eqnarray}
Furthermore, by Theorem~\ref{thm:Non-Optimal} we have for some constants $C_1, C_2 >0$ 
\begin{equation} \label{eq:Tatsachlich-Princeton3}
 V(0) + C_1 \inf_{B_1 (0) } V \leq C_2 .
\end{equation} 

Using \eqref{eq:Tatsachlich-Princeton3} in \eqref{eq:Tatsachlich-Princeton2} we arrive at
$$
  C_0  \sup_{\partial B_1 (0)} V   \, \leq \,  \frac{C_2}{C_1} -  \frac{V(0)}{C_1} +  \tilde C_0 \| e^ \frac{V}{n-1} \|_{L^{\infty}  (A_0)}  ,
$$
and inequality \eqref{eq:AmaliaBis} follows.
\end{proof} 

\bigskip

The following is a slightly more general reformulation of \cite[Corollary 2.4]{EsLu}.

\smallskip

\begin{cor} \label{cor:Amalia}
Given $0 < a \leq b$, consider a family $\Lambda \subset \{ h \in L^{\infty} (\Omega) \, : \,   a \leq h \leq b \}$ which is equicontinuous at each point, and denote by $\mathcal U$ the set of solutions to 
$$
   \qquad 
   - \Delta_n u = h e^u 
   \quad
   \hbox{ in } \Omega ,
   \quad
   \hbox{ with } h \in \Lambda.
$$
Then, for some constants $c_1, c_2, c_3 >0$ the following holds
\begin{equation} \label{eq:Harn-Sup-Inf}
   \sup_{\partial B_r(p_0)} u \, \leq\,   c_1 - c_2 u(p_0) - (1 + c_2) \log r^n + c_3 \max_{A_r} \big( |x - p_0|^n e^u \big)^{\frac{1}{n-1}},
\end{equation}
for any  $u \in {\mathcal U}$ and any $A_r := B_{2r } (p_0) \setminus \overline{ B_{r/2 } (p_0) }  \subset \Omega$
with $0<r \leq \frac{1}{4} \hbox{dist}(p_0,\partial \Omega)$.
\end{cor} 
\begin{proof}
Let $u \in {\mathcal U}$. By considering the rescaled functions 
$$ V (y) := u (r y + p_0) + \log r^n,
     \qquad 
     H (y) :=  h (r y + p_0),
     $$
for which Problem \eqref{eq:Giovedi} holds. On the one hand, we note that the set of functions
$$
   \tilde \Lambda=\Big\{y \in B_2(0) \to h(ry+p_0): \, h \in \Lambda, \, 0<r \leq \frac{1}{4} \hbox{dist}(p_0, \partial \Omega)\Big\}
   $$
 is equicontinuous at each point in $B_2(0)$ thanks to the to the uniform equicontinuity of $\Lambda$ on $K$ and 
 compactness of the set $K:=\{ry+p_0: \, y \in \overline{B_2(0)},\, 0<r \leq \frac{1}{4} \hbox{dist}(p_0,\Omega^c)\}\Subset \Omega$.
On the other hand, 
 for $ y \in B_2 (0) \setminus \overline{B_{\frac{1}{2}} (0)}$ by setting $x = ry + p_0$ we have
  $$e^{V(y)} =r^n e^{u(x)} = \big( \frac{|x-p_0|}{|y|} \big)^ne^{u(x)} \leq  2  |x-p_0|^n e^{u(x)} .$$
Applying then the inequality \eqref{eq:AmaliaBis} with the class of functions $\tilde \Lambda$, and writing $V$ in terms of $u$ we finally arrive at  \eqref{eq:Harn-Sup-Inf}. 
\end{proof} 

\medskip

\begin{oss}
In dimension two, the $\sup + \inf$ inequality (in the non-optimal form) has been proved by considering the entire set of functions $h \in \Lambda_{a,b}$ ($0 <a \leq b$), 
and therefore the estimate \eqref{eq:Harn-Sup-Inf} remains valid  within this larger class of weights (see  \cite{Shafrir}). However, in both the linear and quasilinear cases the arguments are applied 
along a sequence $(h_k)$ that converges locally uniformly, and therefore considering a family $\Lambda$ as in Cor.~\ref{cor:Amalia} is not quite a restriction (thanks to Ascoli-Arzela Theorem).
\end{oss} 

\bigskip

In full generality, under the requirements \eqref{1553} a blowing-up sequences $u_k$ of solution to the problem~\eqref{758}  may exhibit several ``peaks". More precisely, 
as proved in \cite[Theorem 3.1]{EsLu} (an extension of~\cite[Lemma 4]{LiSh}) to the quasilinear case), 
one can find finitely many local maximum $x_k^1,\dots,x_k^N \in B_2 (0)$  so that: 
\begin{equation} \label{eq:Blowup1}
   x^i_k \to 0, 
   \qquad 
   \mu_k^i := e^{- \frac{u_k (x^i_k)}{n}} \to 0 ,
   \qquad
u_k(x_k^i)=\max_{B_{\mu_k^i}(x_k^i)} u_k, 
\end{equation}
\begin{equation} \label{eq:Blowup3}
 U_k (y) := u_k(\mu_k^i y+x_k^i)+n \log \mu_k^i \to U(y) \hbox{ in }C^1_{loc}(\mathbb{R}^n),
\end{equation}
\begin{equation} \label{eq:Blowup2}
 \frac{\mu_k^i}{|x_k^i-x_k^j|}\to 0,\quad \forall i \not = j ,
 \end{equation}
\begin{equation} \label{eq:Blowup4}
\min\{ |x-x_k^1|^n,\dots,|x-x_k^N|^n \}  e^{u_k} \leq C \hbox{ in }B_2(0) .
\end{equation} 

\bigskip

In the present section we focus on the case $N= 1$, which we refer to as an  {\it ``isolated blow-up solutions sequence"}, in accordance with the terminology used in some of the literature. 
Under this circumstance, we are dealing with a sequence $u_k$ that has a single peak at a maximum point $x_k$ as in \eqref{eq:Pastel2} and satisfies
\begin{equation}  \label{1317-Isolated}
   \left.
    \begin{array}{c}
     \displaystyle \lim_{k \to \infty}  x_k  = 0 ,  \quad
     \vspace{2mm}\\ 
 |x-x_k|^n e^{u_k} \leq C \quad  \hbox{ in } B_2(0)\subset \Omega 
\end{array}
\right\}
\end{equation}
for some constant $C >0$. 

\medskip

Setting $\mu_k := e^{- \frac{u_k (x_k)}{n} } $, the invariance stated in \eqref{eq:Zwei-Expressi} implies that the rescaled function \eqref{eq:Blowup3} satisfies
\begin{equation} \label{1044-Isolated}
  -\Delta_n U_k=h_k(\mu_ky+x_k) e^{U_k} \quad \hbox{in } \Omega_k :=\frac{\Omega -x_k}{\mu_k}, 
  \qquad
  \int_{\Omega_k} e^{U_k} \leq C.
\end{equation}
Since the limiting problem of \eqref{1044-Isolated} is given by \eqref{limitpb},  in each ball $B_{\mu_k} (x_k)$ the solutions $u_k$ for the $n$-Liouville equation looks like the standard bubble $U$ properly rescaled. 
To complete the asymptotic estimate of $u_k$ we must study the behavior of $u_k$ in the ``neck"  $B_1 (x_k) \setminus B_{\mu_k} (x_k)$.  

\medskip
\subsection{Neck estimate for an isolated blow-up sequence}
\label{sub:DecayNecks-Isolated}

We first mention the following decay estimate.

\begin{pro} \label{thm:VENTO-Isolated}
Consider an isolated blow-up sequence $u_k$ of solutions to \eqref{758}, satisfying  \eqref{1553} and~\eqref{1317-Isolated}. 
Then, in the annulus $A_k := B_{1} (0) \setminus \overline{B_{\mu_k} (0)}$ the following holds for some uniform  constants $C, \alpha>0$
\begin{equation} \label{00954-Isolated}
    \qquad
        h_k e^{u_k} \, \leq \, C  \frac{ \mu_k^{\alpha}}{|x-x_k|^{n+\alpha}}  \quad \hbox{ in } A_k .
\end{equation}
\end{pro} 
\begin{proof}
Given any $r \in (\mu_k, 1)$, using  \eqref{1317-Isolated}  in the estimate \eqref{eq:Harn-Sup-Inf}  give for some uniform constant $C>0$
$$    
      \qquad
       u_k (x)   \, \leq \, C  + n c_2 \log \mu_k - \log r^{(n + n c_2)},
       \quad \forall x \in  \partial B_{r} (x_k) .
$$
Hence,  \eqref{00954-Isolated} holds with $\alpha = c_2 n$.
\end{proof} 

\medskip

Obviously, from \eqref{1317-Isolated} we readily derive that \eqref{00954-Isolated} also holds in $B_{\mu_k} (x_k)$, a fact that is not quite useful.
In fact, the relevance of \eqref{00954-Isolated}  is a decay property that  for later purposes is stressed in the following observation.
\begin{oss} \label{oss:Mass-Annulus-Decay}
Whenever \eqref{00954-Isolated} holds in an annulus $A := B_{r_2} (0) \setminus \overline{B_{r_1} (0)} \subset B_1 (0)$, then 
\begin{equation} \label{eq:TrenoTunisia-Isolated} 
  \| h_k e^{u_k} \|_{L^1 (A)}    \stackrel{\eqref{00954-Isolated}}{\leq}
  C \int_{A}
 \frac{ \mu_k^\alpha}{|x-x_k|^{n+\alpha}} \, \leq \, 
C \big( \frac{\mu_k}{r_1} \big)^{\alpha} .
\end{equation}
\end{oss} 

\bigskip

We give now some oscillation bounds for an isolated blow-up sequence (i.e. satisfying \eqref{1317-Isolated})
that rely on Serrin's elliptic estimates stated in the subsection~\ref{subsec:Oscillation}.
\begin{pro} \label{coro1135-Isolated} 
Consider an isolated blow-up solution sequence $u_k$  satisfying \eqref{eq:TaschenTuch-Mouchoir2} and  \eqref{1553}. 
Then, on each annulus~\footnote{We use the notation $A_{R, \delta} := \{ R - \delta < |x - x_k| < R + \delta \}$ }
 $A_{r, 3 \delta} \subset  B_1 (0) \setminus \{ x_k \}$, centered at $x_k$,  we have 
\begin{equation} \label{eq:WhyNot-Isolated}
   {\rm osc}_{A_{r, \delta} } (u_k)  \leq C  \Big\{  \log ( \frac{r}{\delta} ) + 1 \Big\} , 
\end{equation}
\begin{equation} \label{eq:Lucevan-Isolated}
    \Big|  u_k (x) -\fint_{\partial B_r(x_k)}u_k \Big| \, \leq \, C   \Big\{  \log ( \frac{r}{\delta} ) + 1 \Big\} ,  \quad \forall x  \in A_{r, \delta}  .
\end{equation}
\end{pro}    
\begin{proof}
We first stress that   Lemma~\ref{lem:Arrival} and the assumption \eqref{1553} allow to apply Theorem~\ref{thmBMO} (with $f = h_k e^{u_k}$), and we derive
\begin{equation} \label{1615-Isolated}
\sup_k [u_k]_{BMO,B_{\frac{3}{2}}(0)}<+\infty.
\end{equation}

\medskip

{\bf Proof of \eqref{eq:WhyNot-Isolated}:} 
For each $p  \in \partial B_{r} (x_k) $,  consider $B_{3 \delta} (p) \subset A_{r,3 \delta} \subset \Omega \setminus \{ x_k \}$. 
Since  $|x_k - x| > \delta$ forall $x \in B_{2 \delta} (p)$,
in the annulus $A_{r, 2 \delta}$ we in particular get
\begin{equation} \label{eq:Foglie-Isolated}
   \delta^n h_k e^{u_k} \, \leq \,  |x - x_k|^n h_k  e^{u_k(x)} \, \leq  \, C,
   \quad \forall x \in A_{r, 2 \delta} .
\end{equation}

\medskip

By applying  Prop. \ref{prop:Birch} on $B_\frac{3}{2}(0)$  with $f:= h_k e^{u_k}$
(see \eqref{eq:Faust-Liszt1}), the BMO-estimate \eqref{1615-Isolated} 
and inequality \eqref{eq:Foglie-Isolated}, we readily deduce the estimate \eqref{eq:WhyNot-Isolated}  in view of $B_{r+2\delta}(x_k) \subset B_\frac{3}{2}(0)$.

\medskip

{\bf Proof of \eqref{eq:Lucevan-Isolated}:}
Based on \eqref{eq:WhyNot-Isolated}, for all $x \in A_{r, \delta}$ we have 
$$ 
   \Big| \fint_{\partial B_r (x_k)} \{ u_k(x) - u_k(y) \}  d \sigma(y) \Big| \, \leq \,   {\rm osc}_{A_{r, \delta} } (u_k) \, \leq \, C \Big\{  \log ( \frac{r}{\delta} ) + 1 \Big\} .
  $$
 \end{proof}

\medskip
\subsection{Comparison with the nonlinear spherical average}

For an isolated blow-up sequence $u_k$ of solutions to \eqref{758}, we proved in Proposition \ref{coro1135-Isolated} that
$u_k=\fint_{\partial B_{|x-x_k|}(x_k)}u_k+O(1)$.  But, the nonlinearity of $n$-Laplacian operator makes it difficult to get further bounds on the function $\fint_{\partial B_{|x-x_k|}(x_k)}u_k$, and 
a better comparison radial function is provided by the  ``nonlinear spherical average" proposed in Section~\ref{sec:Spherical-Mean}. 
More specifically, for each solution $u_k$ we look at the corresponding nonlinear average function
 $\tilde u_k = \tilde u_k(|x-x_k|)$ defined through the Poisson problem~\eqref{eq:1209} with $f := h_k e^{u_k}$ and $R := \frac{1}{2}$, namely 
\begin{equation} \label{eq:1209Bis-Isolated}
\left\{ 
       \begin{array}{cl}
       -\Delta_n \tilde u_k= \fint_{\partial B_{|x-x_k|}(x_k)} h_k e^{u_k} &\hbox{in }B_{R} (x_k),
       \\
      \tilde u_k=\fint_{\partial B_R(x_k)}  u_k &\hbox{on }\partial B_{R} (x_k) ,
 \end{array}\right.
 \quad R :=  \frac{1}{2} .
\end{equation}

\bigskip

A relevant feature of the ``nonlinear spherical average"~$\tilde u_k$ is the fact that its oscillation can be precisely estimated:
\begin{pro} \label{lem:Tobin-Isolated}
Consider an isolated blow-up solution sequence $u_k$ to \eqref{758}, satisfying  \eqref{1553} and \eqref{1317-Isolated}. 
Then, the associated sequence $\tilde u_k$ solving \eqref{eq:1209Bis-Isolated} satisfies
\begin{equation} \label{eq:Fulbright-Isolated}
   \tilde u_k (x) = 
 \fint_{\partial B_{R}(x_k) } u_k
 - \Big( \frac{1}{n \omega_n} \int_{B_{R} (x_k) } h_k e^{u_k} \Big)^{\frac{1}{n-1}} \log \Big( \frac{|x-x_k|}{R} \Big)\, + \, 
 O (1) ,
\end{equation}
\begin{equation} \label{0942Bis-Isolated} 
   \tilde u (x) - \tilde u (y) = 
   - \Big( \frac{1}{n \omega_n} \int_{B_{R} (x_k)} h_k e^{u_k} \Big)^{\frac{1}{n-1}} \log \Big( \frac{|x - x_k|}{ |y - x_k|} \Big) \, + \,
   O (1) ,
\end{equation}
uniformly for all $x,y \in B_{R} (x_k)  \setminus  B_{ \mu_k} (x_k) $. 
\end{pro}
\begin{proof}
Set $f_k :=  \frac{1}{n \omega_n} h_k e^{u_k}$. From \eqref{eq:Eq-Spherical-Mean} applied in the ball $B_R(x_k)$ with $R := \frac{1}{2}$ we have 
\begin{equation} \label{eq:Fulbright1-Isolated}
   \tilde u_k (r) = \fint_{\partial B_{R} (x_k) }  \! u_k  \, d \sigma +   \int_r^{R}   \big( \int_{B_t(x_k)} f_k  \big)^{\frac{1}{n-1}}   \frac{dt}{t},
\end{equation}
for all $r \in (0, R]$. 
Based on the decay given by \eqref{00954-Isolated} and the inequality $1 - x^{\frac{1}{n-1} } \leq (1 - x)^{\frac{1}{n-1} }$  whenever $x \in [0,1]$,  
for each $\frac{3}{2} \mu_k < t < R$ we have

\begin{eqnarray*}
       \Big(\int_{B_{R}(x_k)} f_k \Big)^{\frac{1}{n-1}} 
       -
       \Big(\int_{B_t(x_k)} f_k \Big)^{\frac{1}{n-1}}
&\leq& 
   \Big(\int_{ B_{R}(x_k) \setminus B_t(x_k)} f_k  \Big)^{\frac{1}{n-1}} \\
  &\stackrel{\eqref{00954-Isolated}}{\leq}& 
 C  \Big( \int_{B_{R}(x_k)\setminus B_t(x_k)} \frac{\mu_k^\alpha}{|x-x_k|^{n+\alpha}}  \Big)^{\frac{1}{n-1}}  \\
 &=&
C \Big( \int_{t}^{R}   \frac{\mu_k^\alpha}{s^{n+\alpha}} s^{n-1} ds\Big)^{\frac{1}{n-1}}  \\
 &\leq& 
 C   \big(  \frac{\mu_k}{t} \big)^{\frac{\alpha}{n-1} }
\end{eqnarray*}
for some uniform $C >0$. Hence, for all $ \frac{3}{2}  \mu_k \leq r \leq R$ we have
$$
     \Big|\int_r^{R}  \Big\{ \big( \int_{B_t(x_k)} f_k \big)^{\frac{1}{n-1}}- \big( \int_{B_{R}(x_k)} f_k \big)^{\frac{1}{n-1}} \Big\} \frac{dt}{t} \Big|
     \, \leq \, 
     C \big( \frac{\mu_k}{r} \big) ^{\frac{\alpha}{n-1} }
     \, \leq \, C  
     $$
which gives 
\begin{equation} \label{eq:Fulbright2-Isolated}
     \int_r^{R}  \big(\int_{B_t(x_k)} f_k\big)^{\frac{1}{n-1}} \frac{dt}{t} \, = \, 
     \big( \int_{B_{R}(x_k)} f_k \big)^{\frac{1}{n-1}} \ \log \big( \frac{R}{r} \big) + O (1)  .
\end{equation}
The conclusion \eqref{eq:Fulbright-Isolated} follows from \eqref{eq:Fulbright1-Isolated}, \eqref{eq:Fulbright2-Isolated}. 
Finally,  the oscillation estimate \eqref{0942Bis-Isolated} is a direct consequence of \eqref{eq:Fulbright-Isolated}. 
\end{proof}

\bigskip

The following shows that the oscillation of $u - \tilde u$ is bounded on some ``shrinking annuli".

\begin{pro} \label{prop:Storia-Isolated}
Let $u_k$ be an isolated blow-up solutions sequence  to \eqref{758}, satisfying~\eqref{eq:TaschenTuch-Mouchoir2} and \eqref{1553}, and
consider for some fixed exponent $ 0 < \theta <  \frac{1}{2}$ the annulus 
$$ 
   A_k^{\prime} := \Big\{   \mu_k  \big( \frac{R}{\mu_k} \big)^{\theta} < |x - x_k| <   R \big( \frac{  \mu_k}{R} \big)^{\theta}  \Big\} \subset B_R (x_k),
   \quad R = \frac{1}{2}
   $$
Then, for some constant $C> 0$ (independent of $\mu_k$) we have
\begin{equation} \label{0955-Isolated}
      \osc_{A'_k} (u_k - \tilde u_k) \, \leq \, C   .
\end{equation}
\end{pro}
\begin{proof}
Consider $A_k = B_{r_2} (x_k) \setminus \overline{ B_{r_1} (x_k)}$  defined with $r_1 = \mu_k  ( \frac{R} {\mu_k})^{ \frac{\theta}{2} }$, $r_2 = R ( \frac{\mu_k}{R} \big)^{ \frac{\theta}{2}} $. 
Since $\mu_k \in (0,R)$ we have  
$A'_k \subset A_k$ with radii $r_i, r'_i$ satisfying the requirements stated in Prop. \ref{prop:Storia-General} (with an exponent $\theta_0  = \frac{\theta}{2 ( 1 - \theta)}$).

\medskip

On the one hand, from the elliptic estimate \eqref{eq:Lucevan-Isolated}, we get
\begin{equation} \label{eq:Hilo}
   \|u_k - \bar u_k \|_{L^{\infty} (A_k)} \, \leq  \, C.
\end{equation}
On the other hand, the decay estimate \eqref{eq:TrenoTunisia-Isolated}  gives  
$\| h_k e^{u_k} \|_{L^1 (A_k)} \leq C   \big( \frac{\mu_k}{R} \big)^{\alpha \frac{\theta}{2} } $, whereas on the ball $B_{r'_1} (x_k)$ we have the following lower bound
   $$   \| h_k e^{u_k} \|_{L^1 (B_{r'_1}(x_k))} \, \geq \, \int_{B_{\mu_k}(x_k) } h_k e^{u_k}=\int_{B_1(0) } h_k(\mu_ky+x_k) e^{U_k}
      \geq a\int_{B_1(0) }  e^{U_k} \geq C .$$
Hence, Prop.~\ref{prop:Storia-General} gives
\begin{eqnarray} 
   [ \osc_{A'_k} ( \bar u - \tilde u) ]^2 \, \leq \,  C \big( \mu_k | \log \mu_k | + 1 \big) .
 \label{1647-Isolated}
\end{eqnarray}
By combining  \eqref{eq:Hilo} and \eqref{1647-Isolated}  we then reach the conclusion \eqref{0955-Isolated}.
\end{proof}

\bigskip

We stress that the estimates  \eqref{0942Bis-Isolated} and \eqref{0955-Isolated} imply the following oscillation bound  
which will play a crucial role to describe precisely the asymptotic profile of $u_k$.  
\begin{cor} \label{cor:ConCalma-Isolated}
Using the notations and same assumptions of Proposition \ref{prop:Storia-Isolated}  we have that
\begin{eqnarray}
     u_k(x)-u_k(y)&=&  \tilde u_k(x)-   \tilde u_k(y)+O(1) \nonumber \\
&=&  
- \Big( \frac{1}{n \omega_n} \int_{B_{R}(x_k)} h_k e^{u_k} \Big)^{\frac{1}{n-1}} \log \big( \frac{|x - x_k|}{|y-y_k|} \big) +O(1) 
\label{0942-Isolated} 
\end{eqnarray}
holds uniformly for all $x,y \in A'_k$. 
\end{cor}

\bigskip \goodbreak
\section{Asymptotic profile of an isolated blow-up sequence} 
\label{sec:AsProfile} 
\medskip

We now turn our attention on the asymptotic profile of an isolated blow-up sequence $u_k$, i.e. a sequence satisfying \eqref{1553}-\eqref{eq:TaschenTuch-Mouchoir2}, with \eqref{eq:Blowup1} to \eqref{eq:Blowup4} holding with $N=1$. 
Let us remind that the rescaled function $U_k(y)=u_k(\mu_k y+x_k)+n\log \mu_k$ satisfies \eqref{1044-Isolated}, and that from \eqref{eq:Blowup3} we have 
\begin{equation}\label{1230-Simple}
             U_k(y)=u_k(\mu_k y+x_k)+n\log \mu_k \hbox{ bounded in }C_{loc}(\mathbb{R}^n) .
\end{equation}

\medskip

A first step in deriving the asymptotic profile of $u_k$ is given by the following global expansion: 
\begin{thm} \label{teo1135-Simple}
Let $(h_k,u_k)$ be a sequence such that \eqref{eq:TaschenTuch-Mouchoir2} and \eqref{1553} hold, with $u_k$ being an isolated blow-up sequence.
Then, in the domain $\overline{B_{\frac{1}{2}}(x_k)} \setminus B_{\mu_k}(x_k)$ we have
\begin{equation} \label{1208-Simple}
   u_k  - \inf_{\partial B_{\mu_k} (x_k)} u_k \, = \,  
-  \Big( \frac{1}{n \omega_n} \int_{B_{\frac{1}{2}}(x_k) } h_k e^{u_k} \Big)^{\frac{1}{n-1}} \log \Big( \frac{|x-x_k|}{\mu_k} \Big)  +O(1) . 
\end{equation}
\end{thm}
\begin{proof} 
Note that the estimate \eqref{eq:WhyNot-Isolated} gives
for some uniform constant $C >0$
\begin{equation} \label{eq:Spettacolo-Simple}
  \qquad \displaystyle  |u_k(x) - u_k(y) | \leq C 
  \quad
\forall  x,y \in \partial B_{\frac{1}{2}} (x_k) .
\end{equation} 
Therefore, \eqref{1230-Simple} with \eqref{eq:Spettacolo-Simple} give  the following information on the rescaled function $U_k$
\begin{equation} \label{eq:TraderJoe}
    \osc_{\partial B_1 (0) } U_k  \leq C ,
      \qquad
     \osc_{\partial B_{\frac{1}{2 \mu_k}} (0) } U_k  \leq C ,
\end{equation} 
and from the decay estimate \eqref{00954-Isolated}, we get for some $C >0$
\begin{equation} \label{1141-Simple}
           0 \leq -\Delta_n  U_k  \, \leq \, \frac{C}{|y|^{n+\alpha}} \quad \hbox{ in } B_{\frac{1}{2 \mu_k}}(0) \setminus \overline{B_1(0) }.
\end{equation}
So, from \eqref{eq:TraderJoe}--\eqref{1141-Simple} and Prop.~\ref{prop:Galois} applied with $U_k/ C^{ \frac{1}{n-1} }  $, we get
$$
    U_k -  \inf_{\partial B_{1} (x_k)} U_k  = -\alpha_k \log |y|+O(1)  \quad 
    \hbox{ in } \overline{B_{\frac{1}{2 \mu_k}}(0)} \setminus B_1(0),
    $$ 
where $\displaystyle \alpha_k :=   \big[  \inf_{\partial B_{1} (x_k)} U_k  - \displaystyle \inf_{\partial B_{\frac{1}{2 \mu_k}}(x_k)} U_k  \big] \big[ \log ( \frac{1}{2 \mu_k} )  \big]^{-1}$.
By scaling back to $u_k$,  we get  the following uniform estimate in $\overline{B_{\frac{1}{2}}(x_k) }\setminus B_{\mu_k}(x_k)$:
\begin{equation} \label{0911-Simple}
         u_k (x) -  \inf_{\partial B_{\mu_k} (x_k)} u_k  \, = \,  -\alpha_k \log \Big( \frac{|x-x_k|}{\mu_k} \Big)  +O(1) ,
\end{equation}
with  $\displaystyle \alpha_k =    \big[  \displaystyle  \inf_{\partial B_{ \mu_k}(x_k)} u_k  - \inf_{\partial B_{\frac{1}{2}}(x_k)} u_k    \big] \big[ \log ( \frac{1}{2 \mu_k} )  \big]^{-1}$.

\medskip

From \eqref{0911-Simple} it follows that for any  $x,y \in \{ \mu_k < |x - x_k| < \frac{1}{2} \}$ we have  
$$
   u_k (x) - u_k (y) = -\alpha_k \log \Big(  \frac{|x-x_k|}{ |y - x_k|} \Big) +O(1),
$$
which can be combined with  the oscillation estimate \eqref{0942-Isolated} to derive
$$
    \alpha_k  \log\big(\frac{1}{ 2 \mu_k}  \big)^{1 - 2 \theta}   = 
   \Big( \frac{1}{n \omega_n} \int_{B_{\frac{1}{2}}(x_k) } h_k e^{u_k} \Big)^{\frac{1}{n-1}} \log\big(\frac{1}{ 2 \mu_k}  \big)^{1 - 2 \theta}  + O(1) .
$$
Consequently,
\begin{equation} \label{0930-Simple}
   \alpha_k \, = \,  \Big( \frac{1}{n \omega_n} \int_{B_{\frac{1}{2}}(x_k) } h_k e^{u_k} \Big)^{\frac{1}{n-1}} \, + \, O \big(\frac{1}{|\log \mu_k|} \big).
\end{equation}
Using \eqref{0930-Simple} in \eqref{0911-Simple} we conclude the proof of the expansion as claimed in  \eqref{1208-Simple}.
\end{proof} 
\begin{oss}
The estimate \eqref{1208-Simple} also shows that $u_k$ and its nonlinear average $\tilde u_k$ have a similar expansion, namely the following 
holds uniformly in $\overline{B_{\frac{1}{2}}(x_k) } \setminus B_{\mu_k}(x_k)$
$$
      u_k \, = \,  \tilde u_k + O(1) .
$$
Indeed, by integrating \eqref{1208-Simple}  on $\partial B_{\frac{1}{2}}(x_k) $, we get
$$
  \fint_{\partial B_{\frac{1}{2}}(x_k) } u_k \, = \, -   \Big( \frac{1}{n \omega_n} \int_{B_{\frac{1}{2}}(x_k) } h_k e^{u_k} \Big)^{\frac{1}{n-1}} \log \big( \frac{1}{\mu_k} \big) +\inf_{\partial B_{\mu_k} (x_k)} u_k + O(1) ,
$$
and consequently,  from  \eqref{eq:Fulbright-Isolated}  we conclude that  $\tilde u_k = u_k + O(1) $.
\end{oss} 
The following will be helpful in connection with the Pohozaev identity that will be used for the proof of the main result of this section.
\begin{lm} \label{lem:Chien-GouBis}
Consider a sequence $(h_k,u_k)$ satisfying \eqref{eq:TaschenTuch-Mouchoir2}, \eqref{1553} with $u_k$ being an isolated blow up sequence, and consider the rescaled function $U_k$  as in  \eqref{1230-Simple}. 
Then, in any annulus $A' = \{ r'_1 = R^{\frac{1}{2} } < |x| < R^{\frac{3}{2}} = r'_2 \}$  with $R \in (1, \frac{1}{ \sqrt{2\mu_k}})$ we have
\begin{equation} \label{eq:DieSchlussfolgerung}
  \int_{A'}  \left\{  |\nabla U_k|^{n-2}(\partial_\nu U_k)^2-\frac{|\nabla U_k|^n}{n} \right\} =
  (n-1)M_k^n \omega_n   \ln \big( \frac{r'_2}{r'_1} \big) + O(1) .
 \end{equation}
\end{lm}
\begin{proof}
Setting $M_k:=\Big( \frac{1}{n \omega_n} \int_{B_{\frac{1}{2}}(x_k) } h_k e^{u_k} \Big)^{\frac{1}{n-1}}$, from Theorem \ref{teo1135-Simple}  the sequence $U_k$ satisfies
\begin{equation} \label{1738}  
     U_k  = -  M_k \log |y|  +O(1) \quad \hbox{in }\overline{B_{\frac{1}{2 \mu_k}}(0)} \setminus B_1(0). 
\end{equation}
By setting $\Psi (y) :=   U_k+M_k \log |y|$, we write
$\displaystyle  
     \nabla U_k \, = \, - M_k \frac{y}{|y|^2}  +   \nabla \Psi (y) $,
 which leads to
\begin{eqnarray}
        |\nabla U_k|^n &=& \big| M_k \frac{y}{|y|^2}  -    \nabla \Psi (y) \big|^n 
        \nonumber \\
   &=&  \frac{M_k^n}{|y|^n}  +  | \nabla \Psi (y) \big| \, O \big( \frac{1}{|y|^{n-1}} + | \nabla \Psi (y)|^{n-1}  \big) 
   \label{eq:Deer1}
 \end{eqnarray} 
 where we have used the vectorial estimate $ |a -b|^n = |a|^n + |b| \, O( |a|^{n-1} + |b|^{n-1} )$. The same argument also shows that
   \begin{equation} \label{eq:Deer2}
   \qquad 
   \,|\nabla U_k|^{n-2}(\partial_\nu U_k)^2 =  \, 
     \frac{M_k^n}{|y|^n}  +  |\nabla \Psi (y)\big| \, O \big( \frac{1}{|y|^{n-1}} + |\nabla \Psi (y)|^{n-1}  \big).
\end{equation}
Notice that in $B_{\frac{1}{2 \mu_k}}(0) \setminus \{ 0 \}  $ we have 
$$         -\Delta_n U_k =h_k(\mu_k y+x_k) e^{U_k}
         \quad
         \hbox{ and } 
         \quad
        \Delta_n ( M_k \log |y|) = 0 ,$$
and we view $u_0 := - M_k \log |y|$ as a comparison function to be used in Prop~ \ref{prop:Benutzen}. More specifically, 
after taking into account \eqref{1738} and that $\displaystyle \lim_{k \to \infty} M_k = \frac{n^2}{n-1}$ (see \eqref{eq:Blowup3} and \eqref{00954-Isolated}), the inequality obtained in Prop~ \ref{prop:Benutzen} reads:
\begin{equation} \label{eq:Bibliothek}
 C \int_{B_{\frac{1}{2\mu_k}} (0)} \eta^n  \Big\{ |\nabla \Psi|^n + \frac{|\nabla \Psi |^2}{|y|^{n-2}}  \Big\} \, \leq \,  
 \int_{B_{\frac{1}{2\mu_k}} (0)}  \eta^n e^{U_k} \, + \,  \int_{B_{\frac{1}{2\mu_k}} (0)}  \Big\{ |\nabla \eta|^n+  \frac{ |\nabla \eta|^2}{|y|^{n-2}} \Big\}
\end{equation}
for any function $\eta \in C_0^{1} (B_{\frac{1}{2 \mu_k}} (0) \setminus \overline{B_1 (0)} )$. As a special case,  we choose two annuli 
 $$
      A'=\{r_1'<|y|<r_2'\} \subset A=\{r_1<|y|<r_2\} \subset B_{\frac{1}{2 \mu_k}(0)} \setminus \overline{B_1(0)}  
 $$
with $r_1=R^{\frac{1}{4}}$, $r_1'=R^{\frac{1}{2}}$,  $r_2'=R^{\frac{3}{2}}$, $r_2=R^2$ for some fixed $R \in (1, \frac{1}{\sqrt{2\mu_k}})$, and consider a cut-off function 
$\eta \in C_0^{1} (A)$ as in \eqref{eq:Choice-Eta}. Since the bounds \eqref{1738} and the estimates \eqref{eq:Regeln1}, \eqref{eq:Regeln2} give 
$$
  \int_A e^{U_k}     \stackrel{\eqref{1738}}{=}
 O(R^{-\frac{M_k-n}{4}})
 \quad \hbox{ and } \quad
  \int_A  \Big\{  |\nabla \eta|^n  + \frac{|\nabla \eta|^2}{|y|^{n-2}} \Big\} =  O(\frac{1}{\log R}) ,
 $$
from the inequality \eqref{eq:Bibliothek} we deduce that
\begin{equation} \label{2000}
            \int_{A'} \Big\{ |\nabla \Psi|^n +  \frac{ |\nabla \Psi|^2}{|y|^{n-2}} \Big\} \leq \frac{C}{\log R}. 
\end{equation}

\medskip

Hence, from the estimates \eqref{eq:Deer1}, \eqref{eq:Deer2} and \eqref{2000} together with  
the Cauchy-Schwarz inequality, we arrive at 

\begin{eqnarray*} 
  \int_{A'}
  \Big\{ |\nabla U_k|^{n-2}(\partial_\nu U_k)^2-\frac{|\nabla U_k|^n}{n} \Big\}  &=&
  \big(1 - \frac{1}{n} \big) M_k^n \int_{A'} \frac{1}{|y|^{n}} \, + \, O \Big( \int_{A'} \Big\{ |\nabla \Psi|^n + \frac{|\nabla \Psi|}{|y|^{n-1}} \Big\}  \Big) 
  \nonumber  \\
  &=&
  (n-1)M_k^n \omega_n  \ln \big( \frac{r'_2}{r'_1} \big) + O(1) 
\end{eqnarray*} 
Therefore, the conclusion \eqref{eq:DieSchlussfolgerung} follows. 
\end{proof}

\medskip

As shown in \cite{EsLu}, the type of decay estimate \eqref{thm:VENTO-Isolated} is exploited to prove that  a blow-up solution sequence satisfies $ \int_{B_{1}(x_k) } h_k e^{u_k} \to c_n \omega_n N$ for some integer $N \in \mathbb N$, with $c_n$ 
as in~\eqref{quantization}.  When $N=1$, the following result provides a more refined estimate, and it contains the main Theorem~\ref{thm:Fundamental-Estimate} stated in the introduction. 

\begin{thm} \label{pro:Fundamental-Giessen}
Let $(h_k,u_k)$ be a sequence satisfying \eqref{eq:TaschenTuch-Mouchoir2} and \eqref{1553} with $\| \nabla h_k \|_{L^{\infty}} \leq C$.
Then, $u_k$ is an isolated blow-up solutions sequence, and we have

\begin{equation} \label{eq:Verzoegern} 
     \Big( \frac{1}{n \omega_n} \int_{B_{\frac{1}{2}}(x_k) } h_k e^{u_k} \Big)^{\frac{1}{n-1}} = \frac{n^2}{n-1} + O \big(\frac{1}{ | \log \mu_k| } \big)  .
\end{equation} 
Furthermore, the rescaled function $U_k(y) =  u_k (\mu_k x + x_k) + n \log \mu_k$ satisfies  
\begin{equation} \label{eq:Hitze-Koln}  
   \qquad
     U_k  = -  \frac{n^2}{n-1} \log |y|  +O(1) \quad \hbox{in }\overline{B_{\frac{1}{2\mu_k}}(0)} \setminus B_1(0), 
  \end{equation}
  and in particular, 
\begin{equation} \label{eq:Hitze-Tunis}  
   u_k  =U_{\mu_k,x_k} +O(1) 
   \quad 
   \hbox{ in } B_1(0) .
\end{equation}
\end{thm} 
\begin{proof}
Let us first justify that $u_k$ is an isolated blow-up sequence.  
On the one hand, from \cite{EsMo}, thanks to the $C^1$-bound on $h_k$, in the ball $B_2 (0)$ the following convergence holds in the sense of measure
\begin{equation} \label{eq:Dorr1}
   h_k e^{u_k} \xrightharpoonup{\, \, k \to \infty \, \,}
 c_n \omega_n \delta_0 .
\end{equation}
On the other hand, from \cite{EsLu} we have 
\begin{equation} \label{eq:Dorr2}
   \int_{B_2 (0) } h_k e^{u_k}  \xrightarrow{ \, \, k \to \infty \, \,} c_n \omega_n N ,
\end{equation}
where $N$ is given by the number of ``peaks" coming into play in  \eqref{eq:Blowup1}--\eqref{eq:Blowup4}. Hence, by comparing \eqref{eq:Dorr1} with \eqref{eq:Dorr2} we deduce that $N=1$. 

\medskip

Concerning \eqref{eq:Hitze-Tunis},  we already know that the estimate holds in $B_{ \mu_k} (x_k) $  thanks to \eqref{1230-Simple}. 
Furthermore, by scaling back \eqref{eq:Hitze-Koln} to $u_k$,  we easily derive that the claim  \eqref{eq:Hitze-Tunis}  holds in $\overline{B_{\frac{1}{2}}(x_k)} \setminus B_{\mu_k} (x_k)$,
and \eqref{1615-Isolated} allows to extend the validity of \eqref{eq:Hitze-Tunis}  to $B_1(0)$.

\medskip Hence we focus on \eqref{eq:Verzoegern} and~\eqref{eq:Hitze-Koln}.  
Setting $H_k(y)=h_k(\mu_k y+x_k)$, the Pohozaev identity applied to the PDE satisfied by  
$U_k(y)=u_k(\mu_k y+x_k)+n\log \mu_k$  (see \eqref{1044-Isolated}) gives

 \begin{equation} \label{eq:Pohozaev} 
  \int_{B_{r}(0)} \big\{ n H_k+  \langle \nabla H_k,y\rangle \big\} e^{U_k} \, = \,
  r  \int_{\partial B_{r}(0)} 
  \Big\{ H_k e^{U_k}+|\nabla U_k|^{n-2}(\partial_\nu U_k)^2-\frac{|\nabla U_k|^n}{n} \Big\},
\end{equation}
in each ball $B_r (0)$ with $r \in (1, \frac{1}{2\mu_k})$.

\medskip

 The estimate \eqref{1738} with the assumption  $|\nabla h_k | \leq C$ and the fact that $M_k-n  >1$ (since $M_k \to \frac{n^2}{n-1}$ by \eqref{eq:Dorr1}), 
 allow to estimate the left-hand side of \eqref{eq:Pohozaev}  as follows
 \begin{eqnarray} 
     \int_{B_{r}(0)} \big\{ n H_k+ \langle \nabla H_k,y\rangle \big\} e^{U_k}  &=& n \int_{B_{\frac{1}{2 \mu_k}}(0)} H_k e^{U_k}
     \, - \, n \int_{\{ r <  |x| < \frac{1}{2\mu_k} \}} H_k e^{U_k} \, + \, O(\mu_k) 
     \nonumber \\
     &=& n^2 \omega_n M_k^{n-1} +O(r^{n-M_k} +\mu_k) .
     \label{eq:PSC-CUNY}
 \end{eqnarray}

Concerning the first boundary term in \eqref{eq:Pohozaev}, using again \eqref{1738} we have
\begin{equation}  \label{eq:PSC-CUNY58}
             \int_{\partial B_{r}(0)}  H_k e^{U_k}   \, \stackrel{\eqref{1738}}{=} \, O (r^{n-M_k-1}).
\end{equation}

Hence, substituting \eqref{eq:PSC-CUNY} and \eqref{eq:PSC-CUNY58} in the Pohozaev identity \eqref{eq:Pohozaev} we get
\begin{equation} \label{eq:Pohozaev2} 
 \qquad   \frac{n^2 \omega_n M_k^{n-1}}{r}   \, = \, \int_{\partial B_{r}(0)} 
  \Big\{  |\nabla U_k|^{n-2}(\partial_\nu U_k)^2-\frac{|\nabla U_k|^n}{n} \Big\} +O \big(r^{n-M_k-1 } + \frac{\mu_k}{r} \big)  .
\end{equation}
By considering  $1 < r'_1 < r'_2 < \frac{1}{2 \mu_k}$ and the corresponding annulus $A' = \{ r'_1 < |x| < r'_2 \}$, 
after integrating the identity \eqref{eq:Pohozaev2} on the interval $[r'_1, r'_2]$ with respect to $r$, we obtain
$$
 \qquad 
 n^2 \omega_n M_k^{n-1}  \, = \, 
 \frac{1}{ \log (\frac{r'_2}{r'_1}) } \int_{A'}  \Big\{  |\nabla U_k|^{n-2}(\partial_\nu U_k)^2-\frac{|\nabla U_k|^n}{n} \Big\}dy + 
 O \Big( \mu_k +      \frac{(r'_1)^{n-M_k}}{ \log (\frac{r'_2}{r'_1}) }    \Big)  .
$$

Choosing the annulus $A'$ as in Lemma \ref{lem:Chien-GouBis} with $R = \frac{1}{2 \sqrt{\mu_k}}$, we then derive
$$
    n^2 \omega_n M_k^{n-1}=(n-1)\omega_n M_k^n+O(\frac{1}{|\log \mu_k|})
 $$
leading to $M_k=\frac{n^2}{n-1}+O(\frac{1}{|\log \mu_k|})$, which is the claim \eqref{eq:Verzoegern}. By using this estimate of $M_k$ into \eqref{1738} we immediately obtain the conclusion \eqref{eq:Hitze-Koln}.
\end{proof}

\bigskip
\section{Oscillation bounds for a general blow-up sequence}

In order to address the question of the optimal  $\sup+\inf$ inequality \eqref{317bis} 
without imposing $C^1$-bound on the admissible class of weight $h$, we must analyze the behavior of a general blow-up solutions sequence that might exhibit several peaks. 
To reach this goal,  the discussion undertaken in Section \ref{sec:Oscillation-Liouville} for an isolated blow-up sequence must be generalized. 

\medskip

Given a sequence $u_k$ that satisfies \eqref{eq:Blowup1}--\eqref{eq:Blowup2}, after relabeling the indices we can assume that 
$u_k (x_k^1)=\displaystyle \max_{B_1 (0)} u_k $ and the distances
$d_k^j :=|x_k^j-x_k^1|$, $j=2, \dots,N$, are ordered
$$0  < d_k^2 \leq d_k^3 \leq \dots \leq d_k^N \ll 1.$$
Denoting $x_k:=x_k^1$ and $\mu_k:=\mu_k^1$, by convention we set  
$$     d_k^1 := \mu_k \qquad d_k^{N+1}:=1.
$$

\medskip

As in the case $N=1$, in  a small ball centered at the points $x_k^i$,  thanks to the property \eqref{eq:Blowup3} we know that after proper rescaling the sequence $u_k$ 
 is well approximated by the standard bubble. However, in contrast to the case $N=1$, the analysis outside those blowup regions involve several  ``necks" given by 
 $\{ d_k^{j} + \varepsilon < |x - x^1_k| < d_k^{j+1} - \varepsilon \}$.  To get estimates in those sets, we will first establish some a priori bound for a any sequence $u_k$ 
that satisfies along $N-$sequences of points $x_k^1,\ldots, x_k^N \in B_2 (0) \subset \Omega$ the properties:
\begin{equation}  \label{1317}
   \left.
    \begin{array}{c}
     \displaystyle \lim_{k \to \infty}  x_k^{j}  = 0 ,  \quad ( j =1, \cdots, N) 
     \vspace{2mm}\\
\min \{|x-x_k^1|^n,\ldots,|x-x_k^N|^n\}  e^{u_k} \leq C \quad  \hbox{ in }B_2(0)\subset \Omega 
\end{array}
\right\}
\end{equation}
for some uniform constant $C>0$.

\medskip

For each $x \in \Omega$ let us first select an index $j_x \in \{1, \cdots, N \}$ for which 
\begin{equation} \label{eq:Slectionj}
    \displaystyle |x - x_k^{j_x} | =  \min_{1 \leq  i\leq N } \big\{ |x - x_k^{i}  |  \big\} .
\end{equation} 
Then, given $B_{3r} (p) \subset \Omega \setminus \{ x^1_k, \cdots, x_k^N \}$, we first notice that for all  $x \in B_{2r} (p)$ we have
\begin{equation} \label{eq:Frozen}
    \qquad  
    |x - x^{i}_k| \geq r , 
    \quad
    \forall i \in \{ 1 , \cdots, N \} ,
\end{equation} 
and in that case
\begin{equation} \label{eq:Frozen2}
    r^n h_k e^{u(x) }  \, \stackrel{\eqref{1317}}{\leq} \,  \frac{C r^n }{|x - x^{j_x}_k |^n}    \, \stackrel{\eqref{eq:Frozen}}{\leq} \,   C,
    \quad
    \forall x \in B_{2 r} (p) ,
\end{equation} 
a fact that will be useful in establishing some oscillation bounds. 
First, we present an analogue of Prop.~\ref{coro1135-Isolated}.
\begin{pro} \label{coro1135} 
Let $u_k$ be a sequence of solutions to \eqref{758} satisfying \eqref{eq:TaschenTuch-Mouchoir2} and \eqref{1553}, with the condition~\eqref{1317} being satisfied
along $N$-sequences $(x^{j}_k)$. Then, up to a subsequence, the following hold:
 \begin{enumerate} \itemsep=6pt
 \def\labelenumi{\rm \bf (\alph{enumi})}
 \item
On each annulus
 $A_{r, 3 \delta} \subset  B_1 (0) \setminus \{ x_k^1, \cdots, x_k^N \}$ (centered at $x_k^1$) we have 
\begin{equation} \label{eq:WhyNot}
   {\rm osc}_{A_{r, \delta} } (u_k)  \leq C  \Big\{  \log ( \frac{r}{\delta} ) + 1 \Big\} , 
\end{equation}
\begin{equation} \label{eq:Lucevan}
    \Big|  u_k (x) -\fint_{\partial B_r(x_k^1)}u_k \Big| \, \leq \, C   \Big\{  \log ( \frac{r}{\delta} ) + 1 \Big\} ,  \quad \forall x  \in A_{r, \delta}  .
\end{equation}
\item
Given $j \in \{ 2, \cdots, N \}$, define
\begin{eqnarray} 
      \mathscr{j} \! &:=& \! \min \big\{  2 \leq j' \leq  j \, : \, \limsup_{k \to \infty} \frac{d^{j}_k}{d_k^{j'}} < \infty \big\} ,
   \label{eq:Brendel}\\
   \mathcal{J}  \! &:=& \! \max  \big\{  j  \leq j' \leq  N \, : \, \limsup_{k \to \infty} \frac{d_k^{j'}}{d^{j}_k } < \infty \big\} .
   \label{eq:Brendel2} 
\end{eqnarray}
Then, by considering for some fixed $\alpha \in (0, 1)$ and $\beta \in (1, \infty)$ the spheres 
$$
{S}^{\mathscr{j}}_{k} := \{ |x - x^1_k| = \alpha d_k^{\mathscr{j} }\} ,
\quad 
S^{\mathcal{J}}_k := \{ |x - x^1_k| = \beta d_k^{\mathcal{J} }\} ,
$$ 
we have
\begin{equation} \label{prop:Neve}
  \qquad  \qquad  \qquad
  |u_k (p) - u_k(q) |  \leq C , 
  \quad
  \forall p \in  S^{\mathcal{J}}_k  , \, \, \forall q \in S_{k}^{\mathscr{j}} .
\end{equation}
\end{enumerate} 
\end{pro}    
\begin{proof}
Remember that Lemma~\ref{lem:Arrival} and the assumption \eqref{1553} allow to apply Theorem~\ref{thmBMO} and consequently 
$\sup_k [u_k]_{BMO,B_{\frac{3}{2}}(0)}<+\infty$. Therefore, 
thanks to~\eqref{eq:Frozen2}, the proof of the claims \eqref{eq:WhyNot} and \eqref{eq:Lucevan} is obtained by following the steps of Prop.~\ref{coro1135-Isolated}.
  
\medskip

{\bf Proof of \eqref{prop:Neve}:}
Consider the ball  $B_R (x^1_k) $ with $R:= 2 \beta d_k^{ \mathcal{J} } $, where we will apply the elliptic estimates stated in Proposition~\ref{eq:Faust1}. 

\medskip

Firstly, from the maximal property defining the index $\mathcal{J}$ we deduce that, up to a subsequence,
$$
   B_{R} (x_k^1) \cap \big\{ x_k^{j'} \, : \, j' \geq {\mathcal{J} +1} \big\} = \emptyset.
$$

Given then $p  \in S^{\mathcal{J} }_k$  it follows that the ball $B_{3 \delta} (p)$ with  $3 \delta :=  (\beta -1) d_k^{\mathcal{J} }$ satisfies
$$
   B_{3 \delta} (p) \subset  B_{R} (x_k^1) \setminus \overline{B_{d_k^\mathcal{J}} (x_k^1)  }  
   \subset B_1 (0)  \setminus \{ x_k^1, \cdots, x_k^N\} . 
 $$
Thus, from \eqref{eq:Frozen2} we get that
$
   \delta^n  \| h_k e^{u_k} \|_{L^{\infty} (B_{2 \delta} (p))}  \leq C 
$.
Taking into account that $\frac{R}{\delta} \leq C$ together with the uniform BMO bound, an application of the estimate \eqref{eq:Faust} yields
\begin{equation} \label{eq:FullbrightSaima2}
  \qquad
  u_k (p) - \fint_{B_{R} (x^1_k)} u_k =  O(1) ,
  \quad \hbox{ on } S^{ \mathcal{J} }_k,
\end{equation} 
in view of $B_{2 \delta} (p)\subset B_R(x_k^1)\subset B_{\frac{3}{2}}(0)$. Next, for $q  \in S^{\mathscr{j} }_k$ consider the ball $B_{3 \delta} (q)$  with $3 \delta =\min\{1 - \alpha,\frac{\alpha}{2}\} d_k^{\mathscr{j}}$.   
We  have $B_{3 \delta} (q)  \subset B_{d_k^{\mathscr{j}}} (x_k^1)\setminus B_{\frac{\alpha}{2}d_k^{\mathscr{j}}} (x_k^1)$, and from the minimal property defining $\mathscr{j}$ 
we also get, up to a subsequence,
$$
   B_{3 \delta} (q) \cap \big\{ x_k^{j'} \, : \,  j' \leq \mathscr{j} -1  \big\} = \emptyset .
$$
Hence,  $B_{3 \delta} (q)  \subset B_1 (0)  \setminus \{ x_k^1, \cdots, x_k^N\} $, and so
from \eqref{eq:Frozen2}  we have
$   \delta^n  \| h_k e^{u_k} \|_{L^{\infty} (B_{2\delta} (q))}  \leq C.$

Taking into account that $\frac{R}{\delta} \leq C \frac{d_k^{\mathcal{J}}}{d_k^{\mathscr{j}}} \leq C$ and the uniform BMO bound, with the estimate \eqref{eq:Faust} we arrive at
\begin{equation} \label{eq:FullbrightSaima1}
  \qquad
  u_k (q) - \fint_{B_{R} (x_k^1)} u_k =  O(1),
  \quad \hbox{ on } S^{\mathscr{j}}_k,
\end{equation} 
in view of $B_{2 \delta} (q)\subset B_R(x_k^1)\subset B_{\frac{3}{2}}(0)$. The claim~\eqref{prop:Neve} follows then from \eqref{eq:FullbrightSaima2} and \eqref{eq:FullbrightSaima1}.
\end{proof} 

\bigskip

The estimates \eqref{eq:WhyNot} and \eqref{prop:Neve} give oscillation bounds in annuli whose radii ratio is bounded. 
Deriving a priori bounds in  the necks $\{ d_k^j < |x - x_k^1| < d_k^{j+1} \}$ satisfying 
\begin{equation} \label{eq:Bavaria}
     \lim_{k \to \infty} \frac{d_k^{j+1}}{d_k^{j}} = \infty,
\end{equation}
will require more work.  A preliminary step will be to improve  in those annuli  the estimate \eqref{eq:Frozen2}  in order to get a better decay of the term $h_k e^{u_k}$. 
For that goal, we need a technical lemma that  uses the following simple fact: given $0 < R_0 < R$ and two balls $B_{R_0} (p_0) \subset B_{R} (p_0) \subset \Omega$, then we have 
\begin{equation} \label{eq:Ext-In}
    \qquad
    |e - b_0| \geq (\frac{R}{R_0} -1) |b_0 - p_0|,
    \quad
    \forall (b_0, e)  \in \overline{ B_{R_0} (p_0)} \times \big( \Omega \setminus B_{R} (p_0) \big) .
\end{equation} 
Indeed, starting from the reverse triangle inequality we have 
$$
    |e - b_0|  \, \geq \,  |e - p_0| -  |b_0 - p_0| \, \geq \, R   -  |b_0 - p_0| \, \geq \, \frac{R}{R_0}  |b_0 - p_0|  -  |b_0 - p_0| ,
$$
which gives \eqref{eq:Ext-In}. This is useful to prove the following.

\medskip

\begin{lm}  \label{lem:CarminaBurana}
With the hypotheses  as in  Prop. \ref{coro1135}, consider $j \in \{ 2, \cdots, N \} $ such that $d_k^{j} < d_k^{j+1}$
and  given any  $1 < s <  \big( \frac{d_k^{j+1}}{d_j^{j}} \big)^{1/2}$ set $ \tilde A_k^{(j)} := \big\{ s d_k^{j} < |x - x^1_k | <   \frac{d_k^{j+1} }{s}  \big\} $.
Then,
\begin{equation}  \label{eq:Taxi}
 \qquad \quad
  |x - x^{1}_k |^n e^{u_k} \leq  \frac{C}{(s -1)^n} 
 \quad
 \hbox{ in  }  B_{ \frac{ d^2_k}{s}} (x^1_k) \cup \tilde A_k^{(j)}  .
\end{equation} 
\end{lm} 
\begin{proof}
By applying  \eqref{eq:Ext-In} in the cases:

\smallskip

\begin{enumerate}\def\labelenumi{\rm(\roman{enumi})}
\item
$x \in  B_{ \frac{ d^2_k}{s}} (x^1_k) $ and $x^{i}_k \not \in B_{ d^2_k} (x^1_k) $ for all $i\geq 2$,
\item
$x \in  B_{ \frac{ d^{j+1}_k}{s}} (x^1_k) $ and $x^{i}_k \not \in B_{ d^{j+1}_k} (x^1_k) $ for all $i\geq j+1$,
\item
$x \not \in  B_{ s d_k^{j}} (x^1_k) $ and $x^{i}_k  \in \overline{ B_{ d^{j}_k} (x^1_k)} $ for all $i\leq j$,
\end{enumerate}

\smallskip

we obtain $ (s -1) |x - x^1_k |  \leq     |x - x_k^{i} |$
 for all $x \in B_{\frac{d_k^2}{s}} (x_k^1) \cup \tilde A_k^{(j)}$, and all indices $i$. 
 Since this remains true for the index~$j_x$ defined in \eqref{eq:Slectionj},  from the assumption~\eqref{1317}  we obtain
$$
     (s -1)^n  |x - x^1_{k}|^n e^{u_k(x)}  \, \leq \,   |x - x^{j_x}_k |^n  e^{u_k(x)} \,  \leq \,   C  
       \quad \forall x \in B_{ \frac{d^2_k}{s}} (x^1_k) \cup   \tilde A_k^{(j)}  ,
$$
which leads to the conclusion \eqref{eq:Taxi}.
\end{proof}

\medskip

We can now state the counterpart of Prop.~\ref{thm:VENTO-Isolated} for a general blowing-up solutions sequence.

\begin{pro} \label{thm:VENTO}
Under the same hypotheses  as in  Prop. \ref{coro1135},  
consider an annulus  
\begin{equation} \label{eq:InfiniteMod-Again}
   A_k^{(j)} = \Big\{ 4  d_k^{j} < |x - x^1_k | <  \frac{d_k^{j+1} }{4}  \Big\}  ,
\end{equation}
 with an index $j \in \{ 2, \cdots, N \} $ for which 
\eqref{eq:Bavaria} holds. Then,  we have 
for some uniform  constants $C, \alpha>0$
\begin{equation} \label{00954}  
    \qquad
        e^{u_k} \, \leq \,   \frac{ C \mu_k^{\alpha}}{|x-x^1_k|^{n+\alpha}}  \quad \hbox{ in } B_{\frac{d_k^2}{4}} (x_k^1) \cup A_k^{(j)}  .
\end{equation}
\end{pro} 
\begin{proof}
Using \eqref{eq:Taxi} we have $|x - x^1_k|^n e^{u_k} \leq C$ holds in $B_{\mu_k} (x^1_k)$ and so
 \eqref{00954} obviously holds in $B_{\mu_k} (x^1_k)$ for any $\alpha >0$. Furthermore by \eqref{eq:Taxi} we have that 
 $$
 \qquad \qquad 
         |x- x_k^1|^n e^{u_k} \leq C 
         \quad \hbox{ in } B_{\frac{d_k^2}{2}} (x^1_k) \cup \Big\{ 2  d_k^{j} < |x - x^1_k | <  \frac{d_k^{j+1} }{2}  \Big\} .
 $$ 
 Combining this fact with the estimate \eqref{eq:Harn-Sup-Inf}  we deduce
$$    
      \qquad
       u_k (x)  \, \leq \, C - c_2 u(x^1_k) - (1 + c_2) \log |x - x^1_k|^n ,
       \quad \forall x \in B_{\frac{d_k^2}{4} } (x_k^1) \cup A_k^{(j)},
$$
for some $c_2>0$, which readily gives \eqref{00954} with $\alpha=nc_2$.
\end{proof}

This decay estimate \eqref{00954} can be used to obtain an oscillation bound on the nonlinear spherical average $\tilde u_k$ 
for a  general blowing-up solutions sequence $u_k$. More specifically, by considering
$\tilde u_k$ is defined as in \eqref{eq:1209Bis-Isolated} with $R := \frac{d_k^{j+1}}{4}$ (in lieu of $R=\frac{1}{2}$), 
and using the estimate \eqref{00954} in each neck $A_{k}^{(j)}$, the same proof given in Prop.~\ref{lem:Tobin-Isolated} gives the following generalization:

\smallskip
\begin{pro} \label{lem:Tobin}
Consider $(h_k, u_k)$ such that \eqref{1553} holds,  and  satisfying condition~\eqref{1317} along $N$-sequences $(x^{j}_k)$. Then, given $j \in \{ 1, \cdots, N\}$ for which \eqref{eq:Bavaria} holds, we have  
  $$ \tilde u_k (x) = 
 \fint_{\partial B_{\frac{d_k^{j+1}}{4}}(x_k) } u_k
 - \Big( \frac{1}{n \omega_n} \int_{B_{\frac{d_k^{j+1}}{4}} (x_k) } h_k e^{u_k} \Big)^{\frac{1}{n-1}} \log \Big( \frac{|x-x_k|}{d_k^{j+1}} \Big)\, + \, 
 O (1) ,$$
\begin{equation} \label{0942Bis} 
   \tilde u_k (x) - \tilde u_k (y) = 
   - \Big( \frac{1}{n \omega_n} \int_{B_{\frac{d_k^{j+1}}{4}} (x_k)} h_k e^{u_k} \Big)^{\frac{1}{n-1}} \log \Big( \frac{|x - x_k|}{ |y - x_k|} \Big) \, + \,
   O (1) ,
\end{equation}
uniformly for all $x,y \in B_{\frac{d_k^{j+1}}{4}} (x_k)  \setminus B_{ 4 d_k^{j}} (x_k)$. 
\end{pro}

\bigskip

For a general blow-up sequence one can extend Proposition \ref{prop:Storia-Isolated} to derive  oscillation of $u_k - \tilde u_k$ on some annuli, and it reads as follows:

\begin{pro} \label{prop:Storia}
Under the same assumptions of Prop.~\ref{lem:Tobin} and \eqref{eq:TaschenTuch-Mouchoir2}, we have
 \begin{equation} \label{0955}
 {\rm osc}_{A^{\prime}_k} (u_k - \tilde u_k) \leq C ,
\end{equation}
where $A^{\prime}_k := \Big\{ d_k^{j} \big( \frac{ d_k^{j+1} }{ d_k^{j} }  \big)^{\theta} < |x - x_k| <  d_k^{j+1}   \big( \frac{ d_k^{j} }{ d_k^{j+1} }  \big)^{\theta}  \Big\}$
for some fixed  $0 < \theta < \frac{1}{2} $.  
\end{pro}
\begin{proof}
By considering the annulus 
$$
   A_k :=   \Big\{ r_1 = d_k^{j} \big( \frac{ d_k^{j+1} }{ d_k^{j} }  \big)^{\frac{\theta}{2} } < |x - x_k| <  r_2 = d_k^{j+1} \big( \frac{ d_k^{j} }{ d_k^{j+1} }  \big)^{\frac{\theta}{2}}   \Big\}
   \subset B_{d_k^{j+1} } (x_k)  
$$
observe that, $ A^{\prime}_k \subset A_k $ and $A_k$ is a subset of the annulus $A_{k}^{(j)}$ defined in~\eqref{eq:InfiniteMod-Again}.  Following the proof of 
Prop. \ref{prop:Storia-Isolated} the conclusion is mainly a consequence of Prop~\ref{prop:Storia-General}  (with $\theta_0 = \frac{\theta}{2(1 - \theta) }$) by also taking into account 
 $$ \| h_k e^{u_k} \|_{L^1 (A_k)}    \stackrel{\eqref{00954}}{\leq}
  C \int_{A_k}
 \frac{ \mu_k^\alpha}{|x-x_k|^{n+\alpha}} \leq C \Big( \frac{\mu_k}{d_k^j} \Big)^{\alpha} \Big( \frac{d_k^{j}}{d_k^{j+1}} \Big)^{\alpha \frac{\theta}{2}} 
 \leq 
 C  \Big( \frac{d_k^{j}}{d_k^{j+1}} \Big)^{\alpha \frac{\theta}{2}} .$$
\end{proof}
The estimates  \eqref{0942Bis} and \eqref{0955} imply the following analogue of Cor.~\ref{cor:ConCalma-Isolated} in a general situation where the possibility of 
multiple necks satisfying \eqref{eq:Bavaria} must be taken into account:
\begin{cor} \label{cor:ConCalma}
Under the same assumptions of Proposition \ref{prop:Storia} we have that
\begin{eqnarray*}
     u_k(x)-u_k(y)&=&  \tilde u_k(x)-   \tilde u_k(y)+O(1) \nonumber \\
&=&  
- \Big( \frac{1}{n \omega_n} \int_{B_{\frac{d_k^{j+1}}{4}}(x_k)} h_k e^{u_k} \Big)^{\frac{1}{n-1}} \log \big( \frac{|x - x_k|}{|y-y_k|} \big)  +O(1) 
\end{eqnarray*}
holds uniformly for all $x,y \in A'_k := \Big\{ d_k^{j} \big( \frac{ d_k^{j+1} }{ d_k^{j} }  \big)^{\theta} < |x - x_k| <  d_k^{j+1}   \big( \frac{ d_k^{j} }{ d_k^{j+1} }  \big)^{\theta}  \Big\}$.
\end{cor}

\bigskip
\section{Asymptotic profile for multi-bubbling solutions}

As discussed in Section \ref{sec:Oscillation-Liouville},
even though we are assuming the boundary oscillation  to be bounded, 
we cannot exclude the existence of several blowing-up local maximum points converging to $0$. 
To understand the profile description we must take into consideration that for each index $j \in \{ 1, \cdots, N \}$, up to a subsequence, we might have 
either $\displaystyle \lim_{k \to \infty} \frac{d^{j+1}_k}{d^j_k} = \infty$, or $\displaystyle \lim_{k \to \infty} \frac{d^{j+1}_k}{d^j_k} <  \infty$. 

\begin{pro} \label{prop:TreeTree}
Given $j \in \{1, \cdots,N\}$ for which \eqref{eq:Bavaria} holds, then
\begin{equation} \label{eq:TreeTree}
        u_k  (x) - \inf_{\partial  B_{4 d_k^j}(x_k)} u_k \, = \,  - \Big( \frac{1}{n \omega_n} \int_{B_{\frac{d_k^{j+1}}{4}}(x_k)} h_k e^{u_k} \Big)^{\frac{1}{n-1}}  \log \Big(  \frac{|x-x_k|}{d_k^j} \Big) + O(1) 
\end{equation}
holds uniformly in the set $\overline{  B_{\frac{d^{j+1}_k}{4}} (x_k) } \setminus B_{4d_k^j} (x_k)$.
\end{pro} 
\begin{proof}
Consider the rescaled function $W_k(y) :=u_k( 4d_k^j  y+x_k)+ n\log (4 d_k^j) $, which satisfies 
$$           -\Delta_n W_k \, = \, h_k( 4 d^j_k y+x_k) e^{W_k} \qquad 
             \quad  \hbox{ in } B_{\frac{d_k^{j+1}}{ 16 d_k^j}}(0) \setminus \overline{ B_1 (0) }.
$$
From the same arguments  used in the proof of Theorem~\ref{teo1135-Simple}, by \eqref{eq:WhyNot} and \eqref{eq:Bavaria} we have
\begin{equation} \label{eq:OSc-W}
   \osc_{\partial B_{ 4 d^j_k} (x_k)} u_k = O(1)
   \quad \hbox{ and } \quad 
      \osc_{\partial B_{\frac{d^{j+1}_k}{4}} (x_k)} u_k = O(1)  .
\end{equation}

From the decay estimate \eqref{00954} and the oscillations bounds \eqref{eq:OSc-W}, the sequence $W_k$ satisfies
$$\left.
  \begin{array}{c}
   \quad
   0\leq -\Delta_n W_k \leq C (\frac{\mu_k}{d_k^j})^{\alpha} |y|^{-(n+\alpha)}
             \quad  \hbox{ in } B_{\frac{d_k^{j+1}}{ 16 d_k^j}}(0) \setminus \overline{ B_1(0) } ,
  \vspace{1mm} \\
   \displaystyle \osc_{\partial B_1 (0)} W_k = O(1)
   \quad \hbox{ and } \quad 
      \osc_{\partial B_{\frac{d^{j+1}_k}{16 d_k^j} } (0)} W_k = O(1)  .
      \end{array}
      \right\}$$
Therefore, applying Prop.~\ref{prop:Galois} with the sequence $ W_k  \, [ C (\frac{\mu_k}{d_k^j})^{\alpha}]^{- \frac{1}{n-1}}$ gives 
$$
   \displaystyle 
   \qquad
   W_k -  \inf_{\partial B_{ 1 } (0)} W_k  =  - \alpha_k \log |y|  + O \big( (\frac{\mu_k}{d_k^j})^{\frac{\alpha}{n-1}}+1 \big), 
             \quad  \hbox{ in }  \overline{B_{\frac{d_k^{j+1}}{ 16 d_k^j}}(0) } \setminus B_1 (0) , 
$$
where 
$\displaystyle \alpha_k :=   \big[   \inf_{\partial B_1 (0)}  W_k   - \displaystyle  \inf_{\partial B_{\frac{d^{j+1}_k}{16 d_k^j} } (0)} W_k   \big] \big[ \log ( \frac{d_k^{j+1}}{16 d^j_k} )  \big]^{-1}$. 
Scaling back, there holds
\begin{equation} \label{959-New}
u_k  (x) - \inf_{\partial B_{4 d_k^j}(x_k)} u_k \, = \, -\alpha_k \log \Big(  \frac{|x-x_k|}{4 d_k^j} \Big) + O(1),  
\end{equation}
uniformly in $\overline{B_{\frac{d_k^{j+1} }{4} } (x_k) }\setminus B_{4 d_k^j}(x_k)$, where $\alpha_k$ can also be rephrased in term of $u_k$. Hence, 
$$
   \qquad
   u_k(x) - u_k (y) = - \alpha_k \log \big( \frac{ |x-x_k| }{ |y-y_k|} \big) + O(1) ,
   \quad
   \forall x,y \in \overline{B_{\frac{d_k^{j+1} }{4} } (x_k) }\setminus B_{4 d_k^j}(x_k) .
$$

\medskip

Comparing with the oscillation of the nonlinear average given in Corollary~\ref{cor:ConCalma} leads then to
$$
 - \alpha_k \log \big( \frac{d_k^{j+1}}{d_k^j} \big)^{(1-2\theta)}  \, =  \,
- \Big( \frac{1}{n \omega_n} \int_{B_{\frac{d_k^{j+1}}{4}}(x_k)} h_k e^{u_k} \Big)^{\frac{1}{n-1}} \log \big( \frac{d_k^{j+1}}{d_k^j} \big)^{(1-2\theta)} +O(1) .
$$
Hence, we get $ \alpha_k  = \Big( \frac{1}{n \omega_n} \int_{B_{\frac{d_k^{j+1}}{4}}(x_k)} h_k e^{u_k} \Big)^{\frac{1}{n-1}}  + O(\frac{1}{ \log \big( \frac{d_k^{j+1}}{d_k^j} \big) } ) $, which together with \eqref{959-New} 
gives the conclusion.
\end{proof}

\bigskip

\begin{pro} \label{teo1135}
Given $j \in \{ 2, \cdots, N \}$, consider the indices ${\mathscr{j}} , {\mathcal J}$ defined by \eqref{eq:Brendel} and \eqref{eq:Brendel2}. Then,  
\begin{equation} \label{eq:Venezia}
        u_k  (x) - 
        \inf_{\partial  B_{\frac{d_k^{{\mathscr{j}}}}{4} } (x_k)} u_k \, = \,  
        - \Big( \frac{1}{n \omega_n} \int_{B_{\frac{d_k^{{\mathcal J} +1} }{4}}(x_k)} h_k e^{u_k} \Big)^{\frac{1}{n-1}}  \log \Big(  \frac{|x-x_k|}{d_k^{\mathscr{j}}} \Big) + O(1) 
\end{equation}
holds uniformly in the set $\overline{  B_{\frac{d^{{\mathcal J} +1}_k}{4}} (x_k) } \setminus B_{4 d_k^{ \mathcal{J}}  } (x_k)$.
\end{pro}
\begin{proof}
By definition of ${\mathcal J}$ we have $\displaystyle \lim_{k \to \infty} \frac{d^{{\mathcal J} +1}_k}{d_k^{\mathcal J} }  = \infty$ and 
$\displaystyle  1 \leq \frac{d^{{\mathcal J}}_k}{d_k^{\mathscr{j}} } \leq C$, 
so Prop.~\ref{prop:TreeTree} yields
\begin{eqnarray}
u_k  (x) - 
        \inf_{\partial  B_{4 d_k^{{\mathcal{J}}}} (x_k)} u_k 
        &=&
        - \Big( \frac{1}{n \omega_n} \int_{B_{\frac{d_k^{{\mathcal J} +1} }{4}}(x_k)} h_k e^{u_k} \Big)^{\frac{1}{n-1}}  \log \Big(  \frac{|x-x_k|}{d_k^{\mathcal{J}}} \Big) + O(1)  
        \nonumber \\
        &=&
        - \Big( \frac{1}{n \omega_n} \int_{B_{\frac{d_k^{{\mathcal J} +1} }{4}}(x_k)} h_k e^{u_k} \Big)^{\frac{1}{n-1}}  \log \Big(  \frac{|x-x_k|}{d_k^{\mathscr{j}}} \Big) + O(1) 
        \label{eq:KopfSchmerzen}
\end{eqnarray}
uniformly on $\overline{  B_{\frac{d^{{\mathcal J} +1}_k}{4}} (x_k) } \setminus B_{4 d_k^{ \mathcal{J}}  } (x_k)$.
Furthermore,  estimate \eqref{prop:Neve} in Prop.~\ref{coro1135} gives the following
$$
    \inf_{\partial  B_{4 d_k^{{\mathcal{J}}} } (x_k)}  u_k =  \inf_{\partial  B_{\frac{d_k^{{\mathscr{j}}}}{4} } (x_k)} u_k + O(1) ,
$$
which combined with \eqref{eq:KopfSchmerzen} gives the conclusion \eqref{eq:Venezia}.
\end{proof}

\bigskip

In order to have a version of the previous result that takes into account all the ``necks", we introduce the following set
\begin{equation} \label{eq:J-infinito}
  J_{\infty} := \Big\{ j \in \{ 2, \cdots, N +1 \} \, : \,  \limsup_{k \to \infty} \frac{d^{j}_k }{d^{j-1}_k} = \infty \Big\},
\end{equation} 
and notice that we have $2, N +1 \in J_{\infty}$. 
Let us order the elements of $J_{\infty}$ using the notation
$$  
    2 =j_1 < \cdots < j_{\tilde N} <  j_{\tilde N +1} =  N+1,
$$ 
and define
\begin{equation} \label{eq:Thurston} 
    \displaystyle M (r)  := \Big( \frac{1}{n \omega_n}  \int_{B_{r} (x_k)} h_k e^{u_k}  \Big)^{\frac{1}{n-1}}.
\end{equation} 
Then, in the general case when $N \geq 1$ in \eqref{eq:Blowup1}--\eqref{eq:Blowup4}, the analogue of Theorem \ref{teo1135-Simple}  reads:

\begin{pro} \label{pro:Astma}
In the annulus $\overline{B_{\frac{1}{4} } (x_k) } \setminus B_{ 4 d_k^N  } (x_k)$ we have
\begin{eqnarray*} 
  \quad 
  u_k (x)  -\inf_{\partial B_{\frac{\mu_k}{4}} (x_k)} u_k  
  &=& 
 -     M (\frac{1}{4})  \log \Big( \frac{|x-x_k|}{d_k^{j_{\tilde N}}} \Big)-  M (\frac{d^2_k}{4}) \log \Big( \frac{d_k^{2}}{\mu_k} \Big) 
 \\
  && \qquad  - \sum_{m =2}^{\tilde N} M  (\frac{d_k^{j_{m}}}{4}  )  \log (\frac{d_k^{j_{m}}}{d_k^{j_{m-1}}}  ) +O(1)  .
  \quad
\end{eqnarray*}
\end{pro}
\begin{proof}
For each $j_m \in J_{\infty}$ notice first that the associated pair of indices $(\mathscr{j}, \mathcal{J})$ defined by  \eqref{eq:Brendel} and \eqref{eq:Brendel2} is given by 
$\mathscr{j} = j_m$ and $\mathcal{J} = j_{m+1}-1$.

\medskip

Since $\displaystyle \lim_{k \to \infty} \frac{d^2_k}{\mu_k}  = \infty$, from Proposition~\ref{prop:TreeTree} applied with the index $j=1$ and \eqref{1230-Simple},  we get
     $$u_k  - \inf_{\partial B_{ \frac{\mu_k}{4}} (x_k)    }  u_k \, = \, 
     - M (\frac{d^2_k}{4}) \log \Big( \frac{|x-x_k|}{ \mu_k} \Big) 
    +O(1) ,$$
in $\overline{B_{\frac{d_k^{2}  }{4} } (x_k) } \setminus B_{ 4\mu_k} (x_k)$. Theorem \ref{teo1135} applied then for each $j=j_m$,  $m=1,\dots,  \tilde N$, gives 
\begin{equation} \label{eq:Second-Iteration-PU}
     u_k  - \inf_{\partial B_{\frac{d_k^{j_m}  }{4} } }  u_k \, = \, 
     -  M  (\frac{d_k^{j_{m+1}}}{4}  )  \log \Big( \frac{|x-x_k|}{d_k^{j_m}} \Big) 
      +O(1) 
\end{equation}
in $\overline{B_{\frac{d_k^{j_{m+1}}  }{4} } (x_k) } \setminus B_{4 d_k^{j_{m+1}-1}} (x_k)$. Hence, by setting $j_0 := 1$ we obtain 
\begin{equation} \label{eq:Noise2} 
     \inf_{\partial B_{ \frac{d_k^{j_{m+1}}}{4}} (x_k)    }  u_k  -  \inf_{\partial B_{\frac{d_k^{j_m}  }{4} }(x_k) }  u_k \, = \, 
     -  M  (\frac{d_k^{j_{m+1}}}{4}  )  \log \Big( \frac{d_k^{j_{m+1}}    }{ 4 d_k^{j_m}     } \Big)  +O(1) 
\end{equation} 
 for all $m= 0, \dots, \tilde N$.
Hence, for  $x \in \overline{B_{\frac{1}{4} } (x_k) } \setminus B_{4 d_k^N } (x_k)$,  by using \eqref{eq:Second-Iteration-PU} and \eqref{eq:Noise2} we obtain
\begin{eqnarray*}
   u_k (x)  -\inf_{\partial B_{\frac{\mu_k}{4}} (x_k)} u_k  &=&  u_k (x)  -\inf_{\partial B_{\frac{d_k^{j_{\tilde N}}}{4}} (x_k)} u_k  
   + \sum_{m=0}^{\tilde N -1} \big( \inf_{\partial B_{\frac{d_k^{j_{m+1}}}{4}} (x_k)} u_k   -\inf_{\partial B_{\frac{d_k^{j_{m}}}{4}} (x_k)} u_k \big) \\
   &=& 
   - M(\frac{1}{4}) \log \Big( \frac{|x - x_k|}{ d_k^{j_{\tilde N}}} \Big)  
   - \sum_{m=0}^{\tilde N -1}M  (\frac{d_k^{j_{m+1}}}{4}  )  \log \Big( \frac{d_k^{j_{m+1}}    }{ 4 d_k^{j_m}     } \Big)  +O(1) 
\end{eqnarray*}
which concludes the proof. 
\end{proof}

\bigskip \goodbreak
\section{Asymptotic expansion of the mass bubbles} 
\label{sec:LowerBound}
\medskip

Given a  blowing-up sequence $u_k$ to the $n$-Liouville equation, we know that \eqref{eq:Blowup1}-\eqref{eq:Blowup4} hold. 
Furthermore, 
in each ball $B_{r_k} (x^j_k)$ with $\frac{r_k}{\mu_k^j} \to \infty$, 
after rescaling $u_k$ with respect to the parameter $\mu_k^j$ and applying Fatou's Lemma, we notice that
$$
   \liminf_{k \to \infty}  \int_{B_{ r_k}(x_k^j) }  h_k e^{u_k}
   \, = \, 
   \liminf_{k \to \infty} \int_{B_{ \frac{r_k}{\mu_k^j} }(x_k^j) }    h_k(\mu^j_k y+x^j_k) e^{U_k}  
   \, \geq \, 
   \int_{\mathbb R^n} h(0) e^U
   \, = \,
   c_n \omega_n ,
$$
where $c_n$ is defined by \eqref{quantization} and $\omega_n$ stands for the volume of the unit ball. 
In particular,
since $\frac{d^2_k}{\mu^1_k} \to \infty$ and $\frac{d^2_k}{\mu^2_k} \to \infty$ (by \eqref{eq:Blowup1}) we have
\begin{equation} \label{eq:Ambee}
  \liminf_{k \to \infty}  \int_{B_{ \frac{d^2_k}{2}}(x_k) }  h_k e^{u_k} \geq c_n \omega_n,
  \qquad
  \liminf_{k \to \infty} \int_{B_{ \frac{d^2_k}{2}}(x^2_k) }   h_k e^{u_k} \geq c_n \omega_n,
\end{equation}
with this latter case happening only when $N \geq 2$. 
In what comes next, we will derive a quantitative asymptotic estimate of such lower bounds. To do this we rely on 
the isoperimetric arguments  of~\cite{ChenLin} that can be extended to the quasilinear problem we are considering.

\medskip
\subsection{A differential inequality} 

Given $ 0 < H  \in L^{\infty} (B_L(0))$ on a ball $B_L (0)  \subset \mathbb R^n$ and a function $U \in C^1(B_L(0))$ satisfying
\begin{equation} \label{eq:ViolinoConcerto}
  - \Delta_n U = H e^U ,
  \quad \hbox{ in } B_L (0),
\end{equation} 
denote its upper level sets  by $\Omega_t := \big\{ y \in B_{L} (0) \, : \, U (y) > t \big\}$, and 
let $\Omega_t^{*}$ be the ball centered at the origin of volume equals to $\Omega_t$.  Setting $t_0:=\displaystyle \max_{ |y|=L} U$ and $t_1:=\displaystyle \max_{|y| \leq L} U$, assume hereafter $\Omega_{t_0}\not= \emptyset$, or equivalently $t_0<t_1$. Notice that $|\Omega_{t_0}|>0$, $|\Omega_t|>0$ iff $t<t_1$ and
\begin{equation} \label{eq:SafedTunisia} 
     \qquad
       \overline{\Omega_{t} } \subset B_{L} (0) ,
       \quad
      \forall  t > t_0.
\end{equation} 
Since $Z=\{y \in B_L(0) \, : \,  \nabla U(y)=0\}$ has zero Lebesgue measure \cite{CiaMaz,Lou} (see also \cite{Sci1}),  by the Lipschitz continuity of $U$ we have that
\begin{equation} \label{1315}
H_1(U(Z))=0, \quad U(Z)=\{t \in \mathbb{R}: \, \exists y \in B_L(0) \hbox{ s.t. }U(y)=t,\ \nabla U(y)=0\},
\end{equation}
where $H_k$ stands for the $k-$dimensional Hausdorff measure, and 
\begin{equation} \label{1917-Appendix}
\Big|\{y\in B_L(0)\,: \, U(y)=t\}\Big|=0, \quad \forall t \in \mathbb{R}
\end{equation}
thanks to the Implicit Function Theorem. Since $|\Omega_t|$ is a continuous decreasing function in $[t_0,t_1]$ thanks to \eqref{1917-Appendix}, 
the Schwarz symmetrization $U^*$ of $U$, defined as
 $$
    U^*:  B_{R_{*}} (0)\to \mathbb R,
    \qquad
    U^* (x) = \sup \big\{ t > t_0  \, :  \, \omega_n |x|^n < |\Omega_t| \big\},
 $$
is a radial decreasing function characterized by the property $|\Omega_{U^*(r)}|=\omega_n r^n$ for every $0\leq r < R_*:=(\frac{|\Omega_{t_0}|}{\omega_n})^\frac{1}{n}$. Since $U$ is Lipschitz continuous and satisfies \eqref{eq:SafedTunisia} , it is known that $U^*$ is also Lipschitz continuous and in particular differentiable a.e.  in  $B_{R_*}(0) $.  We then consider the function
\begin{equation} \label{eq:MassF}
   F(r) :=   \int_{\Omega_{U^* (r)}} \! H  e^U  ,
   \quad
   0 \leq r  < R_* ,
\end{equation}
which is increasing and has the following properties.

\begin{pro} \label{pro:ChenLin}
Let $U$ be a solution to \eqref{eq:ViolinoConcerto} for which $\Omega_{t_0} \not= \emptyset$,  then the function $F$ defined by~\eqref{eq:MassF} satisfies the following

 \begin{enumerate} \itemsep=6pt
 \def\labelenumi{\rm \bf (\alph{enumi})}
\item
$F$ is  Lipschitz continuous, differentiable a.e. and  
\begin{equation} \label{eq:Mice}
    \qquad a  n \omega_n  r^{n-1} e^{U^* (r) }   \, \leq \, F'(r)  \leq  b  n \omega_n   r^{n-1} e^{U^* (r) } ,
    \quad
    \hbox{ for a.e. } r \in (0,R_{*});
\end{equation} 
\item
for a.e. $r \in (0,R_{*})$ such that $F(r) \geq n^n \omega_n$ we have
\begin{equation} \label{eq:Primavera}
r F' (r)  \, \geq \,  
  - \frac{b-a}{b} n^n \omega_n + F(r) \Big(
  n  -  \frac{n-1}{n}  (\frac{F(r)}{n \omega_n})^{\frac{1}{n-1} }  \Big) ,
\end{equation}
\end{enumerate}
where $a:=\displaystyle  \inf_{B_{L} (0)} H$, $b:= \displaystyle \sup_{B_L(0)} H$. 
\end{pro} 
\begin{proof}
The proof follows the arguments of~\cite[Theorem 1.1]{ChenLin}, and we provide it for the sake of completeness. 

\medskip

{\bf Proof of \eqref{eq:Mice}:}
Since $H > 0$, for $r>s$ we have
\begin{eqnarray*}
|F(r)-F(s)|&=&   \int_{\Omega_{U^*(r)} \setminus \Omega_{U^*(s)}}H  e^{U} \,\leq \, b e^{U^*(0)}\Big(|\Omega_{U^*(r)}|-|\Omega_{U^*(s)}|\Big)\\
&=& b \omega_n  e^{U^*(0)} (r^n-s^n) \,\leq \, b n \omega_n R_{*}^{n-1} e^{U^*(0)} |r-s|
\end{eqnarray*}
and then $F$ is a Lipschitz function in $(0,R_{*})$. In particular, $F'(r)$ exists a.e. in $(0,R_{*})$ and, if $F$ is differentiable at $r$, we have that
\begin{eqnarray*}
&& a n \omega_n  r^{n-1} e^{U^*(r)} =a \omega_n \lim_{h \to 0^+}  \frac{(r+h)^n-r^n}{h} e^{U^*(r+h)} \leq \lim_{h \to 0^+} \frac{1}{h} \int_{\Omega_{U^*(r+h)} \setminus \Omega_{U^*(r)}}H  e^{U}\\
&& =F'(r)\leq b \omega_n \lim_{h \to 0^+}\frac{(r+h)^n-r^n}{h} e^{U^*(r)}
=b n \omega_n r^{n-1} e^{U^*(r)}
\end{eqnarray*}
in view of the continuity of $U^*$.

\bigskip

{\bf Proof of \eqref{eq:Primavera}:}
Since $U^*$ is differentiable a.e.  in  $B_{R^*} (0)$,  using the radial coordinate $r \in (0,R_{*})$, for some set $J_1 \subset (0,R_{*})$ with $H_1(J_1)=0$ we have
$$
   \frac{d U^*}{dr} \hbox{ exists for all } r \in (0,R_{*}) \setminus J_1.
$$
Setting
$$
      Z^* :=
      \big\{r \in (0,R_{*}) \setminus J_1\, : \,  \frac{d U^*}{dr}(r)=0 \big\} ,
 $$
we have that  
$$      t \to |\Omega_t| \hbox{ is differentiable for all } t>t_0, \, t \notin U^*( J_1 \cup Z^*),$$
with a derivative at $t = U^*(r)$,  $ r \in (0,R_*) \setminus (J_1 \cup Z^*)$,  given by 
\begin{equation} \label{1805}
   \frac{d}{dt} |\Omega_t| \, = \, \frac{d}{dr} |\Omega_{U^*(r)}| \Big|_{r=(U^*)^{-1}(t)}\frac{d}{dt}(U^*)^{-1}(t) \, = \,
   n \omega_n \frac{  [(U^*)^{-1}(t)]^{n-1} }{  \frac{d U^*}{dr}[(U^*)^{-1}(t)]  } \, = \,  
   n \omega_n \frac{ r^{n-1} }{  \frac{d U^*}{dr}(r)  } .
\end{equation}
Fix $ t>t_0$, $t \notin U(Z)$. By the co-area formula we have
$$\int_{\partial \Omega_t} \frac{d \sigma}{|\nabla U|}  \leq - \frac{d}{dt} |\Omega_t| $$
in view of \eqref{1315}. Applying the H\"older inequality with the conjugate exponents $(\frac{n}{n-1} ,\,  n)$ we obtain
\begin{equation}  \label{eq:Fisio}
 H_{n-1}(\partial \Omega_t) =   \int_{\partial \Omega_t}   |\nabla U|^{- \frac{n-1}{ n}}  |\nabla U|^{ \frac{n -1}{n} } d \sigma  
 \, \leq \,  
     \Big( \int_{\partial \Omega_t} \frac{d \sigma}{|\nabla U|}  \Big)^{\frac{n-1}{n}} 
          \Big( \int_{\partial \Omega_t}  |\nabla U|^{n-1}  d \sigma  \Big)^{\frac{1}{n}}  .
 \end{equation}
By integrating \eqref{eq:ViolinoConcerto} on $\Omega_t$ and applying the divergence theorem, by \eqref{eq:SafedTunisia} we arrive at
$$       \int_{\partial \Omega_t}   |\nabla U|^{n-1}  d \sigma   = \int_{\Omega_t} H e^U .$$
Hence, applying the classical isoperimetric inequality together with \eqref{eq:Fisio}  we obtain
\begin{equation} \label{1826}
 H_{n-1}(\partial \Omega^*_t)^{\frac{n}{n-1}}   \, \leq \,  H_{n-1}(\partial \Omega_t)^{\frac{n}{n-1}}   \, \leq \, \big(- \frac{d}{dt} |\Omega_t| \big) \big(\int_{\Omega_t} H e^U \big)^{\frac{1}{n-1}}
\end{equation}
for all $t >t_0$,  $t \notin U(Z)$.  

\medskip For $0<r<  R_{*}$ with $r \notin J_1\cup J_2 \cup Z^*$,  $J_2= (U^*)^{-1}\big[U(Z)\setminus \big(J_1 \cup U^*(Z^*)\big) \big]$, we have at the same time the validity of \eqref{1805} 
and \eqref{1826} with $t=U^*(r) \notin U(Z)$, which yield to
\begin{equation} \label{eq:Lavandino}
  (n \omega_n)^{\frac{1}{n-1} } r  \, \big( - \frac{d U^*}{dr}  \big) \,  \leq F^{\frac{1}{n-1}}(r)  
\end{equation}
in view of $H_{n-1}(\partial \Omega_{U^*(r)}^*)=n \omega_n r^{n-1}$. Since $H_1 \big(J_2) 
=0$ in view of  \eqref{1315} and $\frac{d U^*}{dr} \not=0$ in $J_2$, we have that \eqref{eq:Lavandino}, being trivially valid when $r \in Z^*$, is valid for all $r \in (0,R_*)$, $r \notin J:=J_1\cup J_2$. Therefore,  there holds
\begin{eqnarray}
    \frac{d}{dr} ( r^n e^{U^*} ) &=& n r^{n-1} e^{U^*} + r^n  e^{U^*}   \frac{d U^*}{dr} 
    \nonumber \\
    &\stackrel{\eqref{eq:Lavandino}}{\geq}& 
    r^{n-1} e^{U^*}   \big\{ n -  \big( \frac{ F}{  n \omega_n }   \big)^{\frac{1}{n-1} }   \big\} 
    \nonumber \\
    &\stackrel{\eqref{eq:Mice}}{\geq}& 
    \left\{ 
    \begin{array}{rl}
        \frac{F'}{a n \omega_n} \big\{ n - \big( \frac{ F}{  n \omega_n }   \big)^{\frac{1}{n-1} }   \big\} 
        &\hbox{ if } F(r) >  n^{n} \omega_n 
        \vspace{1mm} \\    
                \frac{F'}{b n \omega_n} \big\{ n - \big( \frac{ F}{  n \omega_n }   \big)^{\frac{1}{n-1} }     \big\} 
        &\hbox{ if }  F(r) \leq  n^n \omega_n  
        \label{eq:Pasta}  
    \end{array} 
    \right. 
\end{eqnarray}
for all $r \in (0,R_{*}) \setminus  J$ with $H_1(J)=0$.

\medskip When the set $\{r \in (0,R_{*}): F(r) \geq n^n \omega_n\}$ is non empty,  there  exists a unique $r_0 \in (0, R_{*})$ so that $F(r_0) = n^n \omega_n$.  By the Lipschitz continuity of $U^*$ and $F$ we can integrate \eqref{eq:Pasta} on the interval $(0,r)$ for $r > r_0$, to obtain
\begin{eqnarray}
  n \omega_n  r^n e^{U^*(r)}  
  &\geq&   
  \frac{1}{b } \int_0^{r_0}   F' \Big\{ n -  \big( \frac{ F}{  n \omega_n }   \big)^{\frac{1}{n-1} }    \Big\}   
  + \frac{1}{a}  \int_{r_0}^r   F' \Big\{ n -  \big( \frac{ F}{  n \omega_n }   \big)^{\frac{1}{n-1} }      \Big\}  
  \nonumber \\
  &=& 
  -\frac{b-a}{ab} F(r_0) \Big( n  -  \frac{n-1}{n}   (\frac{F(r_0)}{n \omega_n})^{\frac{1}{n-1}}    \Big) 
  + \frac{F(r)}{a} \Big( n -  \frac{n-1}{n}  (\frac{F(r)}{n \omega_n})^{\frac{1}{n-1} }   \Big) 
  \nonumber \\
  &=& 
  - \frac{b-a}{ab} n^n \omega_n + 
   \frac{F(r)}{a} \Big( n  -  \frac{n-1}{n}  (\frac{F(r)}{n \omega_n})^{\frac{1}{n-1} }   \Big).
   \label{eq:Attico}
\end{eqnarray}

\medskip

The conclusion~\eqref{eq:Primavera} follows then by estimating from above the left-hand side of \eqref{eq:Attico} using
the inequality $F'(r) \geq a n \omega_n r^{n-1} e^{U^* (r)}$ proved in~\eqref{eq:Mice}. 
\end{proof} 

\bigskip

The differential inequality~\eqref{eq:Primavera}  can be used to derive the following lower bound on the function $F: [0,R_*] \to  [0,\infty)$ defined in \eqref{eq:MassF}.
\begin{pro} \label{pro:Beethoven-Emperor}
Let $U$ be a solution to \eqref{eq:ViolinoConcerto} for which $F(R_*) := \int_{\{ U > t_0 \}}  H e^{U} > n^n \omega_n$. Then, on the interval
$[r_0, R_{*}] := \{ F(r) \geq n^n \omega_n \}$ we have 
 \begin{equation} \label{eq:Leopardi}
   \big( \frac{F (r)}{n \omega_n} \big)^{\frac{1}{n-1}}   \, \geq \, 
        \frac{n(1+\delta) (\frac{r}{r_0})^{1+\delta} }{\delta +  (\frac{r}{r_0})^{1+\delta }} ,
      \quad
      \forall r \in [r_0, R_*] ,
\end{equation}
where $\delta :=  \frac{a}{b(n-1)}>0 $.
\end{pro} 
\begin{proof} 
If $F(r) \geq n^n \omega_n$,  the inequality~\eqref{eq:Primavera}  with $- \frac{b-a}{b} n^n \omega_n \geq - \frac{b-a}{b} F(r)$ gives
$$
  r F' (r)  \, \geq \,  
  F(r) \Big(
   - \frac{b-a}{b}  + n  -  \frac{n-1}{n}  (\frac{F(r)}{n \omega_n})^{\frac{1}{n-1} }  \Big) ,
   \quad
   r \in [r_0, R_*] .
$$
Since $ F(r_0) = n^n \omega_n$ and using the notation $\delta = \frac{a}{b(n-1)}$, we get
$$r F' (r)  \, \geq \,  
   \frac{n-1}{n} F(r) \Big(
   n(1+\delta)   -  (\frac{F(r)}{n \omega_n})^{\frac{1}{n-1} }  \Big) ,
   \quad
  F(r_0) = n^n \omega_n. $$
On the interval (which is non-empty) 
\begin{equation} \label{eq:R-Tilde}
   [r_0, \tilde R ] :=  \left\{ r \in [r_0, R_*] \, : \,   n \leq \big( \frac{F(r)}{n \omega_n} \big)^{\frac{1}{n-1}} \, < \, n(1+\delta)  \right \},
\end{equation}
after setting  $f := (\frac{F}{n\omega_n})^{\frac{1}{n-1}} $ we arrive at 
\begin{equation} \label{eq:DiffIneqBIS}
   \frac{ f'(r) }{f(r) \big( n(1+\delta) - f(r) \big) }  \, \geq \, \frac{1}{ n r}   \quad  \hbox{ on }  (r_0,\tilde R), 
\qquad
      f(r_0) = n .
  \end{equation}

The inequality \eqref{eq:DiffIneqBIS} can be explicitly integrated, leading to
$$
  \log \Big(\frac{\delta f}{n(1+ \delta) - f }\Big)   \, \geq (1+\delta) \log (\frac{r}{r_0})  ,
   \quad
   \hbox{ on }(r_0, \tilde R) ,
$$
and after rearranging we obtain 
$$
         \big( \frac{F}{n \omega_n} \big)^{\frac{1}{n-1}} = f(r)    \, \geq \,  \frac{  n (1+\delta) (\frac{r}{r_0})^{1+\delta} }{\delta +  (\frac{r}{r_0})^{1+\delta} } ,
           \quad
           r \in (r_0, \tilde R).
$$
Hence, from the definition of $\tilde R$ in \eqref{eq:R-Tilde} we obtain 
\begin{equation} \label{eq:Apps}
   \big( \frac{F}{n \omega_n} \big)^{\frac{1}{n-1}} 
    \, \geq \, 
    \left\{ 
    \begin{array}{cl} 
\frac{  n (1+\delta) (\frac{r}{r_0})^{1+\delta} }{\delta +  (\frac{r}{r_0})^{1+\delta} } ,
           \quad
           &r \in (r_0, \tilde R) , 
           \vspace{3mm} \\
    n(1+\delta) ,
     \quad
     &r \in    [\tilde R, R_*) .
   \end{array}
   \right. 
\end{equation}
Since $1 \geq \frac{s}{\delta+s}$, the inequality \eqref{eq:Apps} leads immediately to the conclusion \eqref{eq:Leopardi}.
\end{proof} 

\bigskip

Since the function $s \mapsto \frac{s}{1 + s}$ is increasing on $(0,\infty)$,  notice that
\eqref{eq:Leopardi} gives 
\begin{equation} \label{eq:Leopardi2}
   \big( \frac{F (r)}{n \omega_n} \big)^{\frac{1}{n-1}}   \, \geq \, 
      \frac{n(1+\delta) r}{\delta r_0+r } \geq 
      n(1+\delta)(1- \frac{\delta r_0}{r}) ,
      \quad
      \forall r \in (r_0, R_*] .
\end{equation}

\medskip
\subsection{Lower bound for the mass of a bubbling sequence} 

The conclusion \eqref{eq:Leopardi2} will be applied to get a lower bound on the mass $\int_{\Omega} h_k e^{u_k}$ along a blowing-up solutions sequence to the $n$-Liouville equation. 
More specifically, we have

\begin{thm}  \label{thm:Schumann}
Consider a sequence $h_k \in \Lambda_a^b$ satisfying \eqref{eq:Peso1}.  Setting $D_k := \frac{d_k^2}{4 \mu_k}$, consider in the ball $B_{D_k} (0)$ the rescaled functions
\begin{equation} \label{eq:Gebuhren} 
    \qquad
    U_k  (y) := u_k (\mu_k y + x_k) +n \log \mu_k,
    \quad 
    H_k(y)=h_k(\mu_k y+x_k),
 \end{equation}
and let 
$$
    \quad
    a_k := \min_{|y| \leq  D_k^\frac{1}{4}} H_k,
    \quad
    b_k := \max_{|y| \leq D_k^\frac{1}{4}} H_k .
    $$
Then,  there exists a constant $C \geq 0$ such that (along a subsequence)
\begin{equation} \label{eq:Schumann} 
   \qquad
     \big( \frac{1}{n \omega_n} \int_{B_{D_k}(0)} H_k e^{U_k} \big)^{\frac{1}{n-1} } 
      \, \geq \,  \frac{n}{n-1}\big(n+ \frac{a_k}{b_k}-1\big)\big(1- C D_k^{-\frac{1}{4}} \big).
\end{equation} 
\end{thm}
\begin{proof}
We start by noting that  $D_k := \frac{d_k^2}{4 \mu_k}$
satisfies $\displaystyle \lim_{k \to \infty} D_k = \infty$ (by \eqref{eq:Blowup2}). Furthermore, we have $H_k (y) \to h(0)$ in $L^{\infty}_{loc} (\mathbb R^n)$ and the rescaled function 
$U_k$ solves~\eqref{1044-Isolated} with $U_k \to U$ where $U$ is the bubble \eqref{eq:Bubble}. 

\medskip

{\bf Claim 1:} There exists $R_0 >0$ such that for all $k$ large we have 
\begin{equation} \label{eq:Schumann-1856}
     \int_{B_{R_0}(0)}  H_k e^{U_k} \, > \,  n^n \omega_n 
\end{equation} 

Indeed,  since 
$$ \lim_{R \to \infty}  \int_{B_R(0)} h (0)  e^{U} = \, c_n \omega_n \, > \, n^n \omega_n ,
$$
we can choose a radius $R_0>0$ such that 
$$  \lim_{k \to \infty} \int_{B_{R_0}(0)}  H_k e^{U_k}=  \int_{B_{R_0}(0)} h (0)  e^{U}  > n^n \omega_n $$
and then \eqref{eq:Schumann-1856} holds for all $k$ large.

\medskip  

{\bf Claim 2:} Given $\sigma \in (0,\frac{1}{4}]$ we have (along a subsequence)
\begin{equation} \label{eq:Fille-Du-Regiment}
   \qquad
   B_{D_k^{\sigma}  }  (0)   \, \subset \,
  \{y \in B_{D_k^{2\sigma}}(0): U_k >  t_{0,k}\}  \,  \subset \, B_{D_k^{4 \sigma} } (0) 
   \,  \subset \, B_{D_k } (0) 
\end{equation}
where $\alpha_k := \Big( \frac{1}{n \omega_n} \int_{B_{{D_k}}(0)} H_k e^{U_k} \Big)^{\frac{1}{n-1}} = (c_n \omega_n )^{\frac{1}{n-1}} + O(1)$ and $t_{0,k}:=\displaystyle \max_{|y|=D_k^{2\sigma}} U_k$.

\medskip

The last inclusion follows from $D_k >1$ with $\sigma \leq \frac{1}{4}$, and so we focus on the remaining ones. 
On the one hand,  by Prop.~\ref{prop:TreeTree} and the fact that $\displaystyle \inf_{\partial B_4 (0)} U_k = O(1)$, we have
\begin{equation} \label{eq:Schumann-1810}
     \qquad \qquad 
      U_k  (y)   \, = \, - \alpha_k \log |y|+O(1) ,
      \, \, 
      \hbox{ in } \overline{B_{D_k} (0)} \setminus B_4(0)
\end{equation} 
and in particular 
\begin{equation} \label{1901}
t_{0,k}=-\alpha_k \log D_k^{2\sigma}+O(1).
\end{equation}
So, given $y \in \overline{B_{D_k^{\sigma}} (0)} \setminus B_4(0)$ we get 
 \begin{equation} \label{eq:YoelCorta2}
    U_k (y) = - \alpha_k \log |y| + O(1) > - \alpha_k \log D_k^{\sigma} + O(1) > t_{0,k}
\end{equation}
in view of \eqref{1901}. On the other hand, since $U_k \to U$ in $L^{\infty} (B_ 4(0))$ and $D_k \to \infty$, for all $k$ large we have
 \begin{equation} \label{eq:YoelCorta1}
          \inf_{B_4(0) } U_k     >  t_{0,k} 
\end{equation}
in view of \eqref{1901}. Hence, from \eqref{eq:YoelCorta2} and  \eqref{eq:YoelCorta1} we deduce the first inclusion in \eqref{eq:Fille-Du-Regiment}. Concerning the second one, we do have 
$B_4 (0) \subset B_{D_k^{4 \sigma} (0)}$  and for  $y \in  \overline{B_{D_k^{2 \sigma}} (0)} \setminus B_4(0)$ satisfying $U_k (y) >t_{0,k}$, from \eqref{eq:Schumann-1810} and \eqref{1901},  
we get
$$
       \log |y|  \, \leq \, \log (D_k^{2 \sigma} ) + O(1) < \log (D_k^{4 \sigma} ) ,
$$
which concludes the proof of \eqref{eq:Fille-Du-Regiment}.

\medskip
 
Based on Claim 1 and Claim 2, we can apply Proposition~\ref{pro:Beethoven-Emperor} in the ball $B_{D_k^{2\sigma} } (0) $ with the function $F_k(r) := \int_{\Omega_{U_k^* (r)}} H_k e^{U_k}$
on the interval $[r_{0,k}, R_{*,k}]$,  where $r_{0,k} \in (0, R_0]$ and $R_{*,k} $ are defined through the relations: 
$$ \quad
     \int_{B_{r_{0,k}}(0)} H_k e^{U_k} = n^n \omega_n ,
\quad
     |B_{R_{*,k}} (0)| = | \{ y \in B_{D_k^{2\sigma}}(0):U_k > t_{0,k} \} |.$$
We then deduce 
\begin{eqnarray} \label{eq:Sui-Clivi}
    \big( \frac{1}{n \omega_n} \int_{B_{D_k^{4 \sigma}}} H_k e^{U_k} \big)^{\frac{1}{n-1} }
    &\geq&
     \big( \frac{1}{n \omega_n} \int_{\{y \in B_{D_k^{2\sigma}}(0): U_k > t_{0,k} \}} H_k e^{U_k} \big)^{\frac{1}{n-1} } \, =\big( \frac{F_k (R_{*,k})}{n \omega_n}  \big)^{\frac{1}{n-1} } \nonumber\\
&\geq &     \big( \frac{F_k(D_k^{\sigma})}{n \omega_n}  \big)^{\frac{1}{n-1} }   ,
\end{eqnarray} 
where the last inequality follows from the fact that $D_k^{\sigma} \leq R_{*,k}$ (by \eqref{eq:Fille-Du-Regiment}).
Furthermore, since  $r_{0,k} \in (0, R_0]$ the inequality \eqref{eq:Leopardi2} gives for some constant $C >0$
\begin{equation} \label{eq:Sui-Clivi2}
     \big( \frac{F_k(D_k^{\sigma})}{n \omega_n}  \big)^{\frac{1}{n-1} }
      \geq   \frac{n}{n-1}(n+ \frac{a_k}{b_k}-1)(1- \frac{R_0 D_k^{-\sigma}}{n-1} ) .
 \end{equation}
By combining \eqref{eq:Sui-Clivi} and \eqref{eq:Sui-Clivi2} we arrive at
$$
      \big( \frac{1}{n \omega_n} \int_{B_{D_k^{4 \sigma}}} H_k e^{U_k} \big)^{\frac{1}{n-1} } 
      \, \geq \,  \frac{n}{n-1}(n+ \frac{a_k}{b_k}-1)(1- C D_k^{-\sigma} ) ,
$$
and by choosing $\sigma = \frac{1}{4}$ we obtain the conclusion \eqref{eq:Schumann}.
\end{proof}

\medskip For a sequence of functions $h_k$ that have a uniform logarithmic modulus of continuity, Theorem~\ref{thm:Schumann} combined with Prop.~\ref{prop:TreeTree}
lead to the following estimate: 
\begin{pro}  \label{thm:Blizzard}
Consider a sequence $h_k \in    \Lambda^{{\rm Log}}_{a,b} $.
Then, 
\begin{equation}  \label{eq:BuferaDiNeve} 
     u_k (x) - \inf_{\partial B_{4\mu_k} (x_k)} u_k  
     \, \leq \, - \frac{n^2}{n-1} \log \Big( \frac{|x-x_k|}{\mu_k} \Big)  +O(1)
 \end{equation} 
uniformly in $\overline{B_{\frac{d_k^2}{4}} (x_k)} \setminus B_{4 \mu_k}(x_k)$. 
\end{pro} 
\begin{proof}
Recalling the definitions  \eqref{eq:Thurston} and \eqref{eq:Gebuhren}, by applying  \eqref{eq:TreeTree} we get
\begin{eqnarray} \label{eq:Orecchie}
        u_k (x) - \inf_{\partial B_{4 \mu_k} (x_k)} u_k 
        = 
        - M(\frac{d_k^2}{4}) \log \Big( \frac{|x-x_k|}{\mu_k} \Big) +O(1) 
 \end{eqnarray} 
uniformly in $\overline{B_{\frac{d_k^2}{4}} (x_k)} \setminus B_{4\mu_k}(x_k)$, where $D_k=\frac{d_k^2}{4\mu_k}$ and
$$M(\frac{d_k^2}{4})=\big( \frac{1}{n \omega_n}  \int_{B_{D_k}(0)} H_k e^{U_k} \big)^{\frac{1}{n-1} }.$$
Observing that by Ascoli-Arzela $h_k \in\Lambda^{{\rm Log}}_{a,b}$  implies the uniform convergence \eqref{eq:Peso1} (up to a subsequence), Theorem~\ref{thm:Schumann} gives
\begin{equation} \label{1920}
M(\frac{d_k^2}{4}) \geq 
\frac{n}{n-1}\big(n+ \frac{a_k}{b_k} -1\big)  \big( 1 -   C D_k^{-\frac{1}{4}}  \big),
\end{equation}
where $a_k$ and $b_k$ are as in the Theorem~\ref{thm:Schumann}.  Since $h_k \in\Lambda^{{\rm Log}}_{a,b}$ , the term $ 1 - \frac{a_k}{b_k} \geq 0$ can be estimated  as follows
 \begin{eqnarray}
   1 - \frac{a_k}{b_k} &=&  \frac{1}{b_k} ( b_k - a _k )  \, \leq  \, \frac{1}{a} \: \osc_{B_{D_k^{\frac{1}{4}}} (0) } h_k (\mu_k y + x_k) 
    \leq \frac{b}{a \big| \log (\mu_k D_k^\frac{1}{4})  \big|}  
    \, \leq \, 
    \frac{C}{   \log D_k },
    \label{eq:Frazioneab}
\end{eqnarray}
where we have used that $s \mapsto \frac{1}{ | \log s|} $ is increasing in $(0,1)$ and $\mu_k D_k^\frac{1}{4}=\frac{d_k^2}{4 D_k^\frac{3}{4}} \leq D_k^{-\frac{1}{2}}$ for $k$ large.
Since , by using \eqref{eq:Frazioneab} into \eqref{1920} we then get the lower bound
\begin{equation} \label{1921}
M(\frac{d_k^2}{4}) \geq 
\frac{n^2}{n-1}+O\Big(\frac{1}{|\log \frac{d_k^2}{\mu_k}|}\Big),
\end{equation}
which inserted into \eqref{eq:Orecchie} gives the validity of the inequality \eqref{eq:BuferaDiNeve}.
\end{proof}

\medskip

\begin{oss} \label{rem:Inf-Mu} 
It follows from \eqref{eq:Blowup3} that $u_k = -n  \log \mu_k + O(1)$ in $\overline{B_{ 4 \mu_k} (x_k)}$.
\end{oss}

\bigskip
\section{An optimal $\sup+\inf$ inequality}
\label{sec:SupInf} 
\medskip

\begin{pro} \label{pro:BlowUp-Harnack}
Consider a blowing up sequence $u_k$ such that \eqref{eq:TaschenTuch-Mouchoir2} and \eqref{1553} hold.  Then,
\begin{equation} \label{eq:Carnevale}
        u_k(x_k)  + (n-1) u_k (x)  \, \leq \,  O(1) ,
        \quad
        \hbox{ on } \partial B_{\frac{1}{4}} (0) .
\end{equation}
\end{pro}
\begin{proof}
When the properties  \eqref{eq:Blowup1}--\eqref{eq:Blowup4} hold with $N=1$,  we have $d^2_k = 1$, and Prop.~\ref{thm:Blizzard}  together with Remark~\ref{rem:Inf-Mu}  give for $p \in \partial B_{\frac{1}{4}}(0)$ the estimate
$$
       u_k  (p) \, \leq \,   - \frac{n^2}{n-1} \log \big( \frac{1}{\mu_k} \big)  - n \log \mu_k  + O(1)  \, = \, \frac{1}{n-1} \log \mu_k^n + O(1).
$$
From the definition $\mu_k^n := e^{-u_k (x_k)}$ we arrive then to the conclusion \eqref{eq:Carnevale}.

\medskip

In the general case when the sequence $u_k$ develops multiple bubbles concentrating at the same point, 
namely the case when $N \geq 2$ in \eqref{eq:Blowup1}--\eqref{eq:Blowup4}, 
by applying  Prop.~\ref{pro:Astma} with \eqref{1921}  we get after reordering the terms 
the following expansion  for $p \in \partial \overline{B_{\frac{1}{4}}(x_k)}$
\begin{eqnarray*} 
  \quad 
  u_k  (p) -\inf_{\partial B_{\frac{\mu_k}{4}} (x_k)} u_k  
  &\leq &     - 
    \sum_{m=2}^{\tilde N +1}   
M ( \frac{d_k^{j_m} }{4}  )    \log \big( \frac{ d_k^{ j_{m }}    }{ d_k^{j_{m-1} }   } \big)
    -\frac{n^2}{n-1} \log \big( \frac{d_k^{2}}{\mu_k} \big) +O(1) \\
    &&
 =   \sum_{m=3}^{\tilde N +1}   \big[ M ( \frac{d_k^{j_m} }{4}  )  - M ( \frac{d_k^{j_{m-1}} }{4}  )   \big] \log d_k^{j_{m-1}}  + 
\Big[    M ( \frac{d_k^{j_2} }{4}  )   -\frac{n^2}{n-1} \Big] \log d_k^{2} \\ 
&& \qquad +\frac{n^2}{n-1} \log \mu_k  + O(1) 
\end{eqnarray*}
where the indices $j_{m}$ are the one introduced in \eqref{eq:J-infinito}, and 
$\displaystyle M (r)$ has been defined in \eqref{eq:Thurston}.  Since $  j_2  \geq 3$, $B_{\frac{d_k^2}{2}} (x_k) , B_{\frac{d_k^2}{2}}(x_k^2) \subset B_{d_k^{j_2}} (x_k)$ for $k$ large,  by \eqref{eq:Ambee} we have 
$$
     \liminf_{k \to +\infty}\Big[  \big(  \frac{1}{n \omega_n} \int_{B_{\frac{d_k^{j_2}}{4}}(x_k) } h_k e^{u_k} \big)^{\frac{1}{n-1}}-\frac{n^2}{n-1} \Big]\geq \frac{n^2}{n-1}(2^{\frac{1}{n-1}}-1)>0.
 $$

Therefore,
$$
   u_k (p) - \inf_{\partial B_{\frac{\mu_k}{4}} (x_k)} u_k   \, \leq \,  \frac{n^2}{n-1} \log \mu_k \, + \, O(1).
$$
Taking into account Remark~\ref{rem:Inf-Mu}, and the fact that $\mu_k^n := e^{-u_k (x_k)}$ we obtain the claim. 
\end{proof}

\medskip 

We have now all the ingredients to prove the second main result of our paper.

\medskip

\emph{Proof of Theorem \ref{711intro}:} 
By contradiction assume the existence of a compact set 
$K \subset \Omega$ and a sequence $u_k \in {\mathcal U}_M$ so that
\begin{equation} \label{eq:Kiwi2}
   \displaystyle \max_K u_k+(n-1) \displaystyle \inf_\Omega u_k \xrightarrow{\, k \to \infty \,} +\infty .
\end{equation}
{Letting $\bar x_k \in K$ be so that $u_k(\bar x_k)=\displaystyle   \max_K u_k$,  up to a subsequence we can assume $p=\displaystyle \lim_{k \to +\infty} \bar x_k$. Since $\displaystyle \max_K u_k \to +\infty$ as $k \to +\infty$,  we deduce that $p \in \mathcal{S}$ where $\mathcal{S}$ is given by \eqref{eq:Blowup-Set}.  We can assume $p=0$ with $B_2(0) \subset \Omega$ and  by~\cite[Theorem 3.1]{EsLu} find $N$ local maximum points $x_k^1,\dots,x_k^N \in B_2 (0)$  so that \eqref{eq:Blowup1},  \eqref{eq:Blowup3}, \eqref{eq:Blowup2} and \eqref{eq:Blowup4} hold,  up to a subsequence.  Up to relabeling, $x_k:=x_k^1$ satisfies $u_k(x_k)=\displaystyle \max_{B_1(0)} u_k$ and then $u_k(x_k)\geq u_k(\bar x_k)=\displaystyle \max_K u_k$.  Hence, assuming by contradiction that \eqref{eq:Kiwi2} holds, we have
\begin{equation} \label{eq:Kiwi}
   \displaystyle  u_k (x_k) +(n-1) \displaystyle \inf_\Omega u_k \, \to \, +\infty,
   \quad   \hbox{ osc}_{\partial \Omega} (u_k)+   \int_{\Omega} e^{u_k} \leq C .
\end{equation}

\medskip 

Therefore, from Prop.~\ref{pro:BlowUp-Harnack} we have 
$$
  u_k (x_k) + (n-1) \inf_{\Omega} u_k \,\leq \,  C,
$$
which is in contradiction with the assumption \eqref{eq:Kiwi}. Hence, given a compact set $K$,  there is a constant $C(K)$ for which 
$$
    \max_{K}  u + (n-1) \inf_{\Omega} u  \,\leq \,  C (K),
    \quad
    \forall u \in {\mathcal U}_M.
$$


\bibliographystyle{plain}

\end{document}